\documentclass[10pt]{article}

\usepackage{amssymb}

\usepackage{amsthm,mathtools, amsmath}
\usepackage{hyperref}

\usepackage[cmtip,all]{xy}

\newcommand{\BA}{{\mathbb{A}}}

\newcommand{\BC}{{\mathbb{C}}}

\newcommand{\BF}{{\mathbb{F}}\,\!{}}
\newcommand{\BG}{{\mathbb{G}}}

\newcommand{\BN}{{\mathbb{N}}}

\newcommand{\BQ}{{\mathbb{Q}}}

\newcommand{\BZ}{{\mathbb{Z}}}

\newcommand{\Fg}{{\mathfrak{g}}}

\newcommand{\Fm}{{\mathfrak{m}}}

\newcommand{\Fp}{{\mathfrak{p}}}

\newcommand{\Fu}{{\mathfrak{u}}}

\newcommand{\CG}{{\mathcal G}}
\newcommand{\CH}{{\mathcal H}}

\newcommand{\CL}{{\mathcal L}}
\newcommand{\CM}{{\mathcal M}}

\newcommand{\CO}{{\mathcal O}}
\newcommand{\CP}{{\mathcal P}}

\newcommand{\CR}{{\mathcal R}}
\newcommand{\CS}{{\mathcal S}}
\newcommand{\CT}{{\mathcal T}}
\newcommand{\CU}{{\mathcal U}}
\newcommand{\CV}{{\mathcal V}}

\newcommand{\CX}{{\mathcal X}}

\newcommand{\CZ}{{\mathcal Z}}

\DeclareMathOperator{\id}{{id}}
\DeclareMathOperator{\Jor}{{Jor}}
\newcommand{\wt}[1]{\widetilde{#1}}
\newcommand{\wh}[1]{\widehat{#1}}
\DeclareMathOperator{\Emb}{Emb}
\newcommand{\reg}{\mathrm{reg}}
\newcommand{\kernel}{\mathop{\rm Ker}\nolimits}
\DeclareMathOperator{\charpol}{charpol}

\DeclareMathOperator{\rank}{rank}
\DeclareMathOperator{\Hom}{Hom}
\DeclareMathOperator{\End}{End}

\DeclareMathOperator{\Spec}{Spec}
\DeclareMathOperator{\supp}{Supp}
\DeclareMathOperator{\pr}{pr}

\newcommand{\GL}{{\rm GL}}
\newcommand{\SL}{{\rm SL}}
\DeclareMathOperator{\Ad}{Ad}
\DeclareMathOperator{\Lie}{Lie}

\DeclareMathOperator{\Ind}{Ind}
\DeclareMathOperator{\Res}{Res}
\newcommand{\trace}{\mathop{\rm Tr}\nolimits}
\DeclareMathOperator{\Char}{char}

\newcommand{\Gal}{\mathop{\rm Gal}\nolimits}
\newcommand{\Frob}{{\rm Frob}}
\DeclareMathOperator{\Aut}{Aut}
\DeclareMathOperator{\Out}{Out}

\newcommand{\tame}{\mathrm{tame}}
\newcommand{\et}{{\rm et}}
\newcommand{\der}{{\rm der}}
\newcommand{\red}{{\rm red}}
\newcommand{\sep}{{\rm sep}}
\newcommand{\sat}{{\rm sat}}
\newcommand{\ab}{{\rm{ab}}}
\newcommand{\ad}{\mathrm{ad}}
\newcommand{\out}{{\mathrm{out}}}
\newcommand{\nilp}{{\mathrm{nil}}}
\newcommand{\unip}{{\mathrm{uni}}}
\newcommand{\geo}{\mathrm{geo}}
\newcommand\ssi{{\rm ss}}

\newcommand{\dash}{{\textrm{-}}}

\newcommand{\blank}{{\phantom{N}}}
\newcommand{\ublank}{\underline\blank}

\def\longto{\longrightarrow}
\def\into{\hookrightarrow}

\def\eps{\varepsilon}
\def\phi{\varphi}
\def\setminus{\smallsetminus}

\let\oldbullet\bullet
\def\bullet{{\mathchoice{\oldbullet}%
                        {\oldbullet}%
                        {\scriptscriptstyle\oldbullet}%
                        {\oldbullet}}}

\newcommand{\can}{\mathrm{can}}
\DeclareMathOperator{\diag}{diag}
\DeclareMathOperator{\ch}{ch}
\newcommand{\utau}{\underline{\tau}}
\newcommand{\height}{{\mathrm{ht}}}
\makeatletter
\@ifdefinable\equationname{\let\equationname\equationautorefname}
\def\equationautorefname~#1\@empty\@empty\null{(#1\@empty\@empty\null)}%
\@ifdefinable\AMSname{\let\AMSname\AMSautorefname}
\def\AMSautorefname~#1\@empty\@empty\null{(#1\@empty\@empty\null)}%
\@ifdefinable\itemname{\let\itemname\itemautorefname}
\def\itemautorefname~#1\@empty\@empty\null{(#1\@empty\@empty\null)%
}%
\makeatother

\makeatletter
\renewcommand{\theenumi}{\alph{enumi}}

\renewcommand{\theenumii}{\roman{enumii}}

\renewcommand{\p@enumii}{\theenumi$\m@th\vert$}

\renewcommand{\p@enumiii}{\theenumi.\theenumii.}

\renewcommand{\labelitemi}{$\m@th\circ$}
\renewcommand{\labelitemii}{$\m@th\diamond$}
\renewcommand{\labelitemiii}{$\m@th\star$}
\renewcommand{\labelitemiv}{$\m@th\cdot$}
\makeatother

\RequirePackage{aliascnt}
\newcommand{\basetheorem}[3]{%
    \newtheorem{#1}{#2}[#3]
    \newtheorem*{#1*}{#2}
    \expandafter\def\csname #1autorefname\endcsname{#2}
}%
\newcommand{\maketheorem}[3]{%
    \newaliascnt{#1}{#3}
    \newtheorem{#1}[#1]{#2}
    \aliascntresetthe{#1}
    \expandafter\def\csname #1autorefname\endcsname{#2}
    \newtheorem*{#1*}{#2}
}%
\theoremstyle{plain}   %-------------------standard Style-------------------------

\basetheorem{theorem}{Theorem}{section}

\maketheorem{Thm}{Theorem}{theorem}
\maketheorem{Prop}{Proposition}{theorem}
\maketheorem{Prop-Def}{Proposition-Definition}{theorem}
\maketheorem{Cor}{Corollary}{theorem}
\maketheorem{Lem}{Lemma}{theorem}
\maketheorem{SubLem}{Sublemma}{theorem}
\maketheorem{Conj}{Conjecture}{theorem}

\theoremstyle{definition}    %------------text not italic style------------------

\maketheorem{Def}{Definition}{theorem}
\maketheorem{DefRem}{Definition-Remark}{theorem}
\maketheorem{Claim}{Claim}{theorem}
\maketheorem{Alg}{Algorithm}{theorem}

\theoremstyle{remark}    %----------------also text not italic ------------------

\maketheorem{Rem}{Remark}{theorem}
\maketheorem{Ex}{Example}{theorem}
\maketheorem{Ques}{Question}{theorem}
\maketheorem{QuesRem}{Question-Remark}{theorem}
\maketheorem{Ass}{Assumption}{theorem}
\maketheorem{Con}{Convention}{theorem}

\numberwithin{equation}{section}

\advance\textwidth by 1in
\advance\textheight by 2.3cm
\advance\topmargin by -1.1cm
\advance\oddsidemargin by -.5in

\title{On the semisimplicity of reductions and adelic openness for $E$-rational compatible systems over global function fields}
\author{Gebhard B\"ockle, Wojciech Gajda and Sebastian Petersen}

\makeatletter
\def\blfootnote{\xdef\@thefnmark{}\@footnotetext}
\makeatother

\begin{document}
\parindent0em
\parskip.3em

\maketitle

%\begin{document}

%\title[Compatible systems over global function fields]{On the semisimplicity of reductions and adelic openness for $E$-rational compatible systems over global function fields}

%%    Only \author and \address are required; other information is
%%    optional.  Remove any unused author tags.
%
%%    author one information
%% \author[short version for running head]{name for top of paper}
%\author{Gebhard B{\" o}ckle}
%\address{IWR, University of Heidelberg, Im Neuenheimer Feld 368, 69120 Heidelberg, Germany}
%%\curraddr{}
%\email{gebhard.boeckle@iwr.uni-heidelberg.de}
%\thanks{G.B.~received support from the DFG within the FG1920 and the SPP1489}
%
%%    author two information
%\author{Wojciech Gajda}
%\address{Faculty of Mathematics and Computer Science, Adam Mickiewicz University, Umultowska 87, 61614 Pozna\'{n}, Poland}
%%\curraddr{}
%\email{gajda@amu.edu.pl}
%\thanks{W.G.~was partially supported by the NCN grant no.~UMO-2014/15/B/ST1/00128 and the Alexander von Humboldt Foundation. %It funded in particular joint research stays of all authors at Adam Mickiewicz University, Pozna\'n
%}
%
%%    author two information
%\author{Sebastian Petersen}
%\address{Universit\"at Kassel, Fachbereich 10, Wilhelmsh\"oher Allee 71-73, 34121 Kassel, Germany}
%%\curraddr{}
%\email{petersen@mathematik.uni-kassel.de}
%%\thanks{}

%    \subjclass is required.

%\date{\today}

%\dedicatory{}

%    Abstract is required.

\begin{abstract} 
Let $X$ be a normal geometrically connected variety over a finite field $\kappa$ of characteristic~$p$. Let $(\rho_\lambda\colon\pi_1(X)\to \GL_n(E_\lambda)  )_\lambda$ be any semisimple $E$-rational compatible system where $E$ is a number field and $\lambda$ ranges over the finite places of $E$ not above $p$. We derive new properties on the monodromy groups of such systems for almost all $\lambda$ and give natural criteria for the corresponding geometric adelic representation to have open image in an appropriate sense. A key input to our results are automorphic methods and the Langlands correspondence over global function fields proved in \cite{Lafforgue} by L.~Lafforgue.

To say more, let $(\overline{\rho}_\lambda\colon\pi_1(X)\to \GL_n(k_\lambda)  )_\lambda$ be the corresponding mod-$\lambda$ system, where for every $\lambda$ by $\CO_\lambda$ and $k_\lambda$ we denote the valuation ring and the residue field of $E_\lambda$, and where the reduction is done with respect to some $\pi_1(X)$-stable $\CO_\lambda$-lattice $\Lambda_\lambda$ of $E_\lambda^n$. Let also $G_\lambda^{\mathrm{geo}}$ be the Zariski closure of $\rho_\lambda(\pi_1(X_{\overline{\kappa}}))$ in $\GL_{n, E}$ and let $\mathcal{G}_\lambda^{\mathrm{geo}}$ be its schematic closure in $\Aut_{\CO_\lambda}(\Lambda_\lambda)$. Assume in the following that the algebraic groups $G_\lambda^{\mathrm{geo}}$ are connected. 

We prove that for almost all $\lambda$ the group scheme $\mathcal{G}_\lambda^{\mathrm{geo}}$ is semisimple over $\CO_\lambda$ and its special fiber agrees with the Nori envelope of $\overline{\rho}_\lambda(\pi_1(X_{\overline{\kappa}}))$. A comparable result under different hypotheses was proved in \cite{CHT} by Cadoret, Hui and Tamagawa using other methods. As an intermediate result, we show for $X$ a curve that any potentially tame compatible system of mod-$\lambda$ representations can be lifted to a compatible system over a number field, cf.~\cite{Drinfeld-ProSemisimple};  this implies for almost all $\lambda$ the semisimplicity of  the restriction $\overline{\rho}_\lambda|_{\pi_1(X_{\overline{\kappa}})}$. Finally we establish adelic openness  for $(\rho_\lambda|_{\pi_1(X_{\overline \kappa })} )_\lambda$ in the sense of Hui-Larsen \cite{Hui-Larsen}, for $E=\BQ$ in general, and for $E\supsetneq\BQ$ under additional~hypotheses.
\end{abstract}

\maketitle


\tableofcontents

\parindent0em
\parskip.3em

\section{Introduction}
\label{Sec-Intro}

Let $\kappa$ be a finite field of characteristic $p$ with an algebraic closure $\overline\kappa$ and absolute Galois group $\Gamma_\kappa=\Gal(\overline\kappa/\kappa)$. Let $X$ be a normal geometrically connected variety over $\kappa$ with arithmetic fundamental group  $\pi_1(X)$, omitting a base point in the notation. Its base change under $\kappa\to\overline\kappa$ will be $X_{\overline\kappa}$, and we write $\pi_1^\geo(X):=\pi_1(X_{\overline\kappa})$ for the geometric fundamental group of $X$, so that one has the short exact sequence $1\to \pi_1^\geo(X)\to \pi_1(X)\to\Gamma_\kappa\to1$. 

Let further $E$ be a number field and let $\CP_E'$ denote its set of finite places not above $p$. In the following we denote by $\rho_\bullet=(\rho_\lambda)_{\lambda\in\CP_E'}$ an $E$-rational compatible system consisting of a continuous homomorphism $\rho_\lambda\colon\pi_1(X)\to \GL_n(E_\lambda)$ for each $\lambda\in\CP_E'$, where $E_\lambda$ denotes the completion of $E$ at $\lambda$, subject to the usual compatibility condition, fully recalled in \autoref{Sec-Notation}. Throughout the introduction, we assume that all $\rho_\lambda$ are semisimple, possibly by semisimplifying an initially given system. Then the Zariski closure $G_\lambda$ of $\rho_\lambda(\pi_1(X))$ in $\GL_{n,E_\lambda}$ is a reductive subgroup. We also denote by $G_\lambda^\geo$ the Zariski closure of $\rho_\lambda(\pi^\geo_1(X))$ in $\GL_{n,E_\lambda}$. The group  $G_\lambda^\geo$ is semisimple cf. \autoref{Prop-GenChinTriple}. In the pure case this is a result of Deligne cf. \autoref{Thm-GeomSemisimplicity}.

There are two main sources of such systems. First, let $f\colon Y\to X$ be any smooth proper morphism. Then by Deligne, \cite{Deligne-Weil2}, for any $i\in\BZ$, the family or higher direct images $(R^i_\et f_*\BQ_\ell)_{\ell\in\CP_\BQ'}$, provides a $\BQ$-rational compatible system (pure of weight $i$). Second, suppose $X$ is a curve over $\kappa$ and $\Pi$ is an cuspidal automorphic representations for $\GL_n$ over the adele ring of the function field $\kappa(X)$ of $X$. Then by L.\ Lafforgue,  \cite{Lafforgue}, and by Drinfeld if $n=2$, \cite{Drinfeld-OnGL2}, to $\Pi$ one can attach an $E$-rational compatible system, where $E$ depends on~$\Pi$. We shall refer to systems of the former kind as {\em cohomological} and of the latter kind as {\em automorphic}.

Roughly, the motivation for the present work is to study under what conditions on $\rho_\bullet$ the homomorphism
\[\rho_{\BA}:=\prod_{\lambda\in\CP_E'} \rho_\lambda\colon\pi_1(X)\longto G(\BA_E^p):=\mathop{\prod\nolimits^\prime}\limits_{\lambda\in\CP_E'} G_\lambda(E_\lambda)\]
has open image, where $\prod'$ denotes the restricted product for a suitable choice of compact open subgroups of $G_\lambda(E_\lambda)$. If there would be a Mumford-Tate like group, as is often the case if $\pi_1(X)$ is replaced with the absolute Galois group of a number field, then one would expect the $G_\lambda^o$ to arise from a single group over $E$. Therefore we think of the above as an adelic openness question. Unlike in the number field case just mentioned, there is a well-known obstruction to adelic openness if the groups $G_\lambda$ have non-trivial torus quotients, that stems from the smallness of $\Gal (\overline\kappa/\kappa)\cong\hat\BZ$. Therefore one has to study the above question either under the hypothesis that all $G_\lambda$ are semisimple, or for $\pi_1^\geo(X)$ when the groups $G_\lambda^\geo$ are semisimple. One can show that the geometric case does imply the arithmetic case, and so, for simplicity of exposition, for the rest of the introduction we shall focus on the geometric case and always have superscripts~$\geo$.

Unless $G^\geo_\lambda$ is simply connected for almost all $\lambda$, it has been well-known for a while that the above is not the right question, even for $E=\BQ$. A good formulation of adelic openness is given in \cite{Hui-Larsen} by Hui and Larsen. For this, denote by $\wt G^\geo_\lambda$ the universal cover of $G^\geo_\lambda$, and by $ \langle\rho_{\BA}(\pi^\geo_1(X)),\rho_{\BA}(\pi^\geo_1(X))\rangle$ the group generated by the commutators of $\rho_{\BA}(\pi^\geo_1(X))$. The latter can, in a natural way, be regarded as a subgroup of $ \wt G^\geo(\BA_E^p):=\prod{}^\prime_{\lambda\in\CP{}_E^\prime}\wt G^\geo_\lambda(E_\lambda)$, where we assume that for almost all $\lambda$ we have a model of $\wt G_\lambda$ defined over $\CO_\lambda$ so that the restricted product makes sense. Then it is suggested in \cite{Hui-Larsen} to study the openness of 
\begin{equation}\label{Eqn-One}
 \langle\rho_{\BA}(\pi^\geo_1(X)),\rho_{\BA}(\pi^\geo_1(X))\rangle\subset  \wt G^\geo(\BA_E^p).
\end{equation}

For cohomological $\BQ$-rational compatible systems, the openness of the inclusion \autoref{Eqn-One} is an immediate consequence of the recent article \cite{CHT} by Cadoret, Hui, Tamagawa. In the present work, in \autoref{Cor-OpennessAboveEll0} and \autoref{Cor-Reduction-IsSaturated}, we shall give a proof for arbitrary $\BQ$-rational compatible systems that is independent of \cite{CHT}. We also indicate a second proof of this fact based on \cite{CHT} and on the work of L.\ Lafforgue, using some reduction techniques of the present work. For general $E$ there are obvious further conditions, as can be deduced from \cite{Pink-Compact}. Namely one needs that $E$ is the subfield of $\overline\BQ$ that is generated over $\BQ$ by $\{\trace(\rho_\lambda(\Frob_x))\mid x\in X\}$. We shall deduce explicit conditions from \cite{Pink-Compact}, under which this holds. We shall also single out one class of $E$-rational compatible system where we can prove fully the openness of the inclusion \autoref{Eqn-One}; see \autoref{Thm-OpennessAboveEll2}.
 
A main problem in studying adelic openness is to understand the image of the reduction $\bar\rho_\lambda$ of $\rho_\lambda$ for almost all $\lambda\in\CP_E'$, and in fact most of the present work concerns precisely this question. Let us choose an $\CO_\lambda$-lattice $\Lambda_\lambda$ in $E_\lambda^n$ that is stable under the action of $\pi_1(X)$ via $\rho_\lambda$.\footnote{In \autoref{Cor-Reduction-IsSaturated} we show that it suffices to choose $\Lambda_\lambda$ stable under the action of $\pi_1^\geo(X)$.} Following \cite[Prop.~1.3]{LarsenPink95}, the groups $\CG^\geo_\lambda$, defined as the Zariski closure of $G_\lambda^\geo$ in $\Aut_{\CO_\lambda}(\Lambda_\lambda)$ endowed with the unique structure of reduced closed subscheme, are smooth group schemes over $\CO_\lambda$ for almost all $\lambda$. Let $\CG^\geo_{k_\lambda}$ denote the special fiber of $\CG_\lambda^\geo$, regard it as a subgroup of $\GL_{n,k_\lambda}$ via a choice of basis of $\Lambda_\lambda$, and let $\bar\rho_\lambda$ be the reduction of $\rho_\lambda$ with respect to the lattices $\Lambda_\lambda$. Then we have the following inclusions of groups:
\[ \bar\rho_\lambda(\pi_1^\geo(X))\subset \CG^\geo_{k_\lambda}(k_\lambda)\subset \GL_n(k_\lambda).\]
Unlike $G^\geo_\lambda$, the group $\CG_{k_\lambda}^\geo$ is a priori not known to be reductive. Now for subgroups of $\GL_n$ over finite fields Nori has defined an algebraic hull in \cite{Nori}, often called the Nori envelope. We shall write $\bar\rho_\lambda(\pi_1^\geo(X))^\sat_{k_\lambda}$ for the closed subgroup of $\GL_{n,k_\lambda}$ that is the Nori envelope of $\bar\rho_\lambda(\pi_1^\geo(X))$. Nori's construction has been generalized by Serre to the concept of saturation, which is the notion we will in fact use. We shall prove the following main result:
\begin{Thm}[{\autoref{Thm-RedIsSaturated}, \autoref{Thm-OnIrred} and \autoref{Prop-GConnSemisimple}, \autoref{Thm-AI-Reduction-IsSaturated} and \autoref{Thm-Reduction-IsSaturated}}] \label{Thm-Intro} Suppose $G^\geo_\lambda$ is connected for all $\lambda$. Then for almost all $\lambda\in\CP_E'$ the following hold:
\begin{enumerate}
\item The group $\CG^\geo_{k_\lambda}$ is saturated.
\item The group $\bar\rho_\lambda(\pi_1^\geo(X))^\sat_{k_\lambda}$ is semisimple, and thus independent of the choice of $\Lambda_\lambda$.
\item The inclusion $\bar\rho_\lambda(\pi_1^\geo(X))^\sat_{k_\lambda}\subset \CG^\geo_{k_\lambda}$ is an equality and $\CG^\geo_\lambda$ is semisimple. 
\end{enumerate}
\end{Thm}
The first major result in this direction is \cite{Larsen95} due to Larsen, who proved the above for $\BQ$-rational compatible systems not for almost all $\ell$, but for all $\ell$ in a set of density one.\footnote{The formulation in \cite{Larsen95} is somewhat different, but it is not difficult to carry out the transition.} Very recently the result as stated was proved for cohomological $\BQ$-rational compatible systems in the important work \cite[Thm.~7.3 and Cor.~7.5]{CHT}. We shall give proofs of the three statements above with methods independent of those of either \cite{Larsen95} or  \cite{CHT}. 

As an immediate consequence of \autoref{Thm-Intro}~(c) and results of Serre, we obtain the following result, which for geometric $\BQ$-rational systems is \cite[Thm.~1.1]{CHT}.
\begin{Cor}
Suppose $G^\geo_\lambda$ is connected for all $\lambda$. Then for almost all $\lambda\in\CP_E'$, the representation $\bar\rho_\lambda|_{\pi_1^\geo(X)}$ is semisimple and independent of the choice of~$\Lambda_\lambda$.
\end{Cor}

Several of our key steps rely on automorphic methods, mainly \cite{Lafforgue} and \cite{Chin04}, and only indirectly on~\cite{Deligne-Weil2}. In fact we shall first prove the theorem for automorphic $E$-rational compatible systems. Moreover several of our arguments benefit from the freedom of enlarging an initially given coefficient field. 

Part~(a) of \autoref{Thm-Intro} is obtained from combining the work of Larsen-Pink on the groups $\CG_\lambda$ with work of Nori, Serre and others on saturation. We exploit the fact that the root system of $\CG^\geo_{k_\lambda}$ looks like that of a reductive group, even when it is not known to be the case. The result of (a) is stated in \autoref{Thm-RedIsSaturated}. Parts~(b) and~(c) are first proved in the case where $\rho_\bullet$ is absolutely irreducible, and using a result of Drinfeld, we also reduce to the case of a curve $X$, see~\autoref{Cor-ReductionToCurve}. In that case, using automorphic results, we use a further result of Drinfeld who showed that $\bar\rho_\lambda$ is absolutely irreducible for almost all $\lambda$, see \autoref{Thm-OnIrred}. By standard results on saturation, given in \autoref{Prop-GConnSemisimple}, (b) is then an immediate consequence. 

From (a) and (b) it follows that one has the inclusion indicated in (c). We then use that the inclusion $\bar\rho_\lambda(\pi_1^\geo(X))^\sat_{k_\lambda}\subset \GL_{n,k_\lambda}$ is of low $\ell_\lambda$-height (for $\ell_\lambda\gg0$), and that the smaller group is saturated, to lift this representation over the Witt vectors, still in the absolutely irreducible case. Whenever the inclusion of connected groups in (c) is proper, and $\Char k_\lambda$ is large enough, this leads to a congruence of cuspidal automorphic forms for $\GL_n$ with trivial determinant and uniformly bounded conductor. An argument similar to (b) implies that only finitely many congruences are possible, i.e., that $\bar\rho_\lambda(\pi_1^\geo(X))^\sat_{k_\lambda}\subset \CG_{k_\lambda}^\geo$ can be proper only for finitely many $\lambda$, and this completes (c) in the absolutely irreducible case; see \autoref{Thm-AI-Reduction-IsSaturated}. We note that we include a proof for (b) also to stress that the method used here allows one to prove the following result -- for the notion of (tame) mod $\lambda$ compatible system and for the notation $\bar\rho_\bullet\otimes_{\CO_E}\CO_{E'}$ used below we refer to \autoref{Def-ModLambdaCS}.
\begin{Thm}\label{Thm-Intro2}
Suppose that $X$ is a curve and that $\bar\rho_\bullet$ is a potentially tame $E$-rational semisimple mod $\lambda$ compatible system. Then there exists a finite extension $E'$ of $E$ and an $E'$-rational compatible system $\rho_\bullet$, unique up to isomorphism, whose semisimplified reduction is $\bar\rho_\bullet\otimes_{\CO_E}\CO_{E'}$.
\end{Thm}

The general case of \autoref{Thm-Intro} is obtained by a reduction procedure to the absolutely irreducible case, explained in the proof of \autoref{Thm-Reduction-IsSaturated}. Given any compatible system $\rho_\bullet$, one can in a systematic way find an absolutely irreducible compatible system $\rho'_\bullet$, of possible different dimension, such that there is a map $G^\geo_\lambda\to G^{\prime,\geo}_\lambda$ with finite kernel, for all $\lambda$. We show in the rather technical \autoref{Lem-IndepOfLattice} that this can be done, by the use of an intermediate system that allows one to move integral properties also, in such a way that if (b) and (c) hold for $\rho'_\bullet$, then they hold for $\rho_\bullet$. \autoref{Lem-IndepOfLattice} is at the heart of the independence statement in \autoref{Thm-Intro}(b).

We find the above an interesting reduction method in its own right. It relies on work of Chin who in turn builds on work of L.\ Lafforgue. If $E$ is sufficiently large and if the $G_\lambda$ are semisimple, then Chin essentially shows that the compatible system $\rho_\bullet$ in fact arises from an $M$-compatible system $r_\bullet=(r_\lambda\colon\pi_1(X)\longto M(E_\lambda))_{\lambda\in\CP_E'}$, with Zariski dense image, where $M$ is split semisimple over $E$, by composing $r_\bullet$ with a representation $\alpha\colon M\to\GL_n$, defined over $E$ that is independent of $\lambda$. It is stated in \autoref{Cor-ChinRepOverM}, following  \cite{Chin04} and \cite[Prop.6.6]{BHKT}. Choosing an irreducible (almost faithful) representation $\alpha'$ instead, and using $\alpha\oplus\alpha'$ as an intermediate step, makes our argument possible. For instance one can take for $\alpha'$ the tensor product of the adjoint representations of the simple quotients of~$M$; see also \autoref{Rem-OnChoicesOfSubstitutes}. We call the group $M$ the split motivic group of $\rho_\bullet$ over $E$ and $\alpha$ the split motivic representation, following~\cite{Serre-Motives}; see \autoref{Def-ChinTriple}. The motivic group satisfies $M\otimes_EE_\lambda\cong G_\lambda$ (if $\rho_\bullet$ itself is $E$-rational).

This decomposition into tensor factors also allows one to reduce the original question about the adelic openness of \autoref{Eqn-One} to the case where $M$ is absolutely simple, cf.~\autoref{Lem-OpennessInOtimes} and \autoref{Prop-GaloisConjugates}. For such $M$ and for $\alpha$ the adjoint representation, we shall prove adelic openness.

\begin{Cor}[{\autoref{Cor-HuiLarsen}; see also \autoref{Thm-OpennessAboveEll2}}]\label{Cor-IntroCor1} Let $\rho_\bullet$ be an $E$-rational semisimple compatible system such that all $G_\lambda^\geo$ are connected. Suppose either that (i) $E=\BQ$ or that (ii) the split motivic group $M$ is simple of adjoint type, the split motivic representation $\alpha$ is the adjoint representation\footnote{The pair $(M,\alpha)$ is possibly defined only over some finite extension of $E$.}  and $E$ is the subfield of $\overline\BQ$ that is generated over $\BQ$ by $\{\trace(\rho_\lambda(\Frob_x))\mid x\in X\}$.\footnote{Given the stated properties on $(M,\alpha)$, by \cite{Pink-Compact} the representation $\rho_\bullet$ may always be defined over this $E$.} Then 
\[ \langle\rho_{\BA}(\pi^\geo_1(X)),\rho_{\BA}(\pi^\geo_1(X))\rangle\subset  \wt G^\geo(\BA_E^p)\]
is an inclusion of an open subgroup.
\end{Cor}

Furthermore, in our setting we answer Conjecture~5.4 of \cite{LarsenPink95} in the affirmative.
\begin{Cor}[{see~\autoref{Cor-MonodromyIsUnramified}}]\label{Cor-IntroCor2}
Suppose that $\rho_\bullet$ is an $E$-rational compatible system. Then the groups $G^{\geo,o}_{\lambda}$ are unramified for almost all $\lambda$.
\end{Cor}
For cohomological $\BQ$-rational compatible systems the above two corollaries are also straightforward consequences of~\cite{CHT}. \autoref{Cor-IntroCor2} was proved in \cite{LarsenPink95} for compatible systems attached to abelian varieties over global function fields.

\medskip

Let us end the introduction with a brief survey of the individual sections, noting that each section will begin with an outline of its contents. \autoref{Sec-Notation} collects the main notations used throughout the article. \autoref{BasicNew} recalls some standard results on compatible systems of Galois representations, tailored to the setting of this work, where the Galois group is the arithmetic fundamental group $\pi_1(X)$ of a normal variety $X$ over a finite field $\kappa$, or its geometric counterpart $\pi_1^\geo(X)$. In doing so, we also fix some further notation, and we present extensions of existing results or give proofs of presumably known results that we could not locate in the literature. Themes are the Weil restriction for compatible systems, Serre's result on the independence of $\lambda$ of $G_\lambda/G_\lambda^o$, the motivic group and representation, the reduction of general $X$ to the curve case, an extension of the independence of images result of \cite{BGP15}, and a reminder on Deligne's result on the semisimplicity of the representations $\rho_\lambda|_{\pi_1^\geo(X)}$.

The subsequent \autoref{Sec-AI} provides various results on absolutely irreducible compatible systems. We begin by collecting a number of basic and mostly elementary results. In \autoref{Prop-GenChinTriple}, we generalize a result of Chin from pure to arbitrary semisimple compatible systems. In \autoref{Subsec-Gcompat}, building on \cite{Chin04} and \cite{BHKT}, we explain how to use tools from representation theory to study compatible systems. We recall the notion of a $G$-compatible system for $G$ a reductive group. Then we show that any semisimple compatible system with motivic group $M$ can be recovered from an $M$-compatible system (with Zariski dense image) together with the motivic representation $\alpha$ -- if $E$ is large. By passing to other (almost faithful) representations of $M$, one can obtain the same group $M$, up to a uniform finite kernel, from various choices of absolutely irreducible compatible systems; see \autoref{Cor-ChinAndAI2-New}.  This corollary serves as a crucial reduction step in the later \autoref{Thm-Reduction-IsSaturated}, but is also important in \autoref{Subsec-ChinAdjoint}.

\autoref{Sec-Saturation} focusses on saturation. After recalling the foundational correspondence of Nori between exponentially generated subgroups and nilpotently generated Lie algebras, the second subsection provides many results on saturation in the sense of Serre that we could not locate in the literature, and we link saturation with Nori's results. Our main interest is in saturation within $\GL_n$ but in the third subsection, in relation to the Weil restriction and representations of low $\ell$-height, we also recall the notion of saturation within a general reductive group. In \autoref{Subsec-NonReduSat} we show for a new class of groups that they are saturated. In the last subsection, this is applied to the special fibers of the smooth group schemes $\CG_\lambda$ introduced in \cite{LarsenPink95} and recalled above \autoref{Thm-Intro}. This will prove \autoref{Thm-Intro}(a). We also discuss that the relevant results from \cite[\S~1]{LarsenPink95} hold for $E$- and not only $\BQ$-rational compatible systems.

In \autoref{Sec-ResidualSaturation}, we prove three main results on $E$-rational (mod $\lambda$) compatible systems $\rho_\bullet$. After some preliminary observations in \autoref{Subsec-CharAndLifting},  in \autoref{Thm-MainGenlCase} in \autoref{SubSect-ModLambdaSystems} we shall prove that any (suitably defined) potentially tame $E$-rational mod $\lambda$ compatible system (over an infinite index set) is the reduction of a compatible system. From this we shall deduce \autoref{Thm-OnIrred} in \autoref{Subsect-ModLambdaAI}, which says that any absolutely irreducible compatible system has absolutely irreducible reduction for almost all $\lambda$. The proof of this was first published by Drinfeld in \cite{Drinfeld-ProSemisimple}, and we basically give a more detailed version of his proof. The final \autoref{Subsec-ModLmabdaSaturated} contains the proof of \autoref{Thm-AI-Reduction-IsSaturated}, which is part (c) of \autoref{Thm-Intro} in the absolutely irreducible case. 

\autoref{Sec-IndepOfLat} is concerned with the proof of \autoref{Thm-Intro}(c) in general. In the first subsection we shall explain that the geometric monodromy groups $G^\geo_\lambda$ of a compatible system $\rho_\bullet$, even integrally, arise from an absolutely irreducible compatible system, up to a uniform finite kernel. In the second subsection, we will deduce from this \autoref{Thm-Intro}(c), and moreover \autoref{Cor-IntroCor2} as well as a precursor of \autoref{Cor-IntroCor1} for $\BQ$-rational compatible systems. 

\autoref{Sec-OpenImage} begins by recalling important results of Pink from \cite{Pink-Compact} in \autoref{Subsec-OnPink}  on the structure of compact subgroups of $G(F)$, where $G$ is a reductive group and $F$ is a product of $\ell$-adic fields. From this we shall in \autoref{Subsec-ChinAdjoint} deduce \autoref{Cor-IntroCor1} -- mainly the assertion under hypothesis (ii). In \autoref{Rem-AdelicOpen} we shall explain with a sketch of proof how far one could weaken the hypothesis in (ii) of \autoref{Cor-IntroCor1} to still have an adelic openness result. The subsection also explains how, up to some left-open details in \autoref{Rem-AbsIrrThmFromCHToverQ}, one can directly from \cite{CHT} obtain \autoref{Thm-Intro}(c) without using our results from \autoref{Sec-ResidualSaturation}.

{\bf Acknowledgements:} 
Part of this work was done during a stay of the three authors at Adam Mickiewicz university in Pozna\'n financed by NCN grant no. UMO-2014/15/B/ST1/00128. G.B. was supported by the DFG in the FG1920 and within the SPP1489. W.G. was partially supported by the NCN grant no. UMO-2014/15/B/ST1/00128 and the Alexander von Humboldt Foundation. The authors thank Prof.\ G. McNinch, for valuable feedback, in particular on \autoref{Sec-Saturation}, and the anonymous referee for a careful reading and many useful suggestions.

%{\bf Acknowledgements:} 
%%Part of this work was done during a stay of the three authors at Adam Mickiewicz university in Pozna\'n financed by NCN grant no. UMO-2014/15/B/ST1/00128. G.B. was supported by the DFG in the FG1920 and within the SPP1489. W.G. was partially supported by the NCN grant no. UMO-2014/15/B/ST1/00128 and the Alexander von Humboldt Foundation. 
%The authors thank Prof.\ G. McNinch, for valuable feedback, in particular on \autoref{Sec-Saturation}, and the anonymous referee for a careful reading and many useful suggestions.

\section{Notation}
\label{Sec-Notation}
Notation and conventions.
\begin{itemize}
\item We fix a prime number $p$ throughout and write $\ell$ for some prime number different from~$p$.
\item $\kappa$ is a finite field with prime field $\BF_p$ and cardinality $q$.
\item $X$ is a normal geometrically irreducible variety over $\kappa$ with function field $\kappa(X)$ of dimension at least $1$; its set of closed points is denoted by $|X|$ and the residue field of $x\in |X|$ is $\kappa_x$.
\item $X_\reg$ denotes the set of regular points of $X$; the subset $X_\reg$ is dense open in~$X$.
\item For any field $K$, $\overline K$ denotes an algebraic closure of $K$, $K^\sep\subset\overline K$ the separable closure of $K$ in $\overline K$, and $\Gamma_K=\Gal(K^\sep/K)$ its absolute Galois group. 
\item By $\pi_1(X)$ we denote the \'etale fundamental group of $X$; its geometric fundamental group is $\pi_1^\geo(X):=\pi_1(X\times_\kappa\overline\kappa)$; one has a short exact sequence $1\to \pi_1^\geo(X)\to \pi_1(X)\to \Gamma_\kappa\to1$.
\item If $K$ is a finite field, we denote by $\Frob_K\colon \overline K\to\overline K,\alpha\mapsto \alpha^{1/\#K}$ the geometric Frobenius of~$\Gamma_K$. For $x\in |X|$ and a choice $\phi\colon\pi_1(x)\to\pi_1(X)$ we write $\Frob_x$ for $\phi(\Frob_{\kappa_x})$.
\item For a (pro-)finite group $\Gamma$ denote by $\Gamma_\ell^+$ the (topological closure of the) normal subgroup of $\Gamma$ generated by all (pro-)$\ell$-Sylow subgroups of $\Gamma$. If $\ell$ is clear from context, we simply write $\Gamma^+$.
\item For a linear algebraic group $G$ over a field $K$, let $G^o$ denote its identity component, $G^\der$ the derived group of $G^o$, and $R_u(G)$ the unipotent radical of $G$. The adjoint representation will be $\Ad_G\colon G\to\Aut(\Lie G)$.
\item In this text semisimple and reductive groups need not be connected.
\item $E$ denotes a number field; any number field will carry a fixed embedding $E\into \overline\BQ$ into a fixed algebraic closure of $\BQ$; the set of finite places of $E$ not above $p$ is $\CP_E'$. For $\lambda\in\CP_E'$ we denote by $\ell_\lambda$ the rational prime under $\lambda$, by  $E_\lambda$ the completion at $\lambda$, by $\varpi_\lambda$ a uniformizer, by $\CO_\lambda$ its ring of integers and by $k_\lambda$ its residue field; for $E'$ a finite extension field of $E$ and $\lambda'\in\CP_{E'}'$, we write $E'_{\lambda'}$, $\CO'_{\lambda'}$ and $k'_{\lambda'}$, and we write $\lambda'|E$ for the restriction of $\lambda$ to $E$. We also use this notation for $E=\overline\BQ$, so that for $\lambda\in \CP_{\overline\BQ}'$ we have the objects $\overline\BQ_\lambda$, $\overline\BZ_\lambda$ and $k_\lambda$.
\item If $K$ is a local field or a number field, by $\CO_K$ we denote its ring of integers.
\item We call $\alpha\in\overline\BQ$ {\em plain of characteristic $p$}, cf.~\cite[Sect.~2]{Chin04}, if $\alpha\in(\overline\BZ[1/p])^\times$, and we call $P\in\overline\BQ[T]$ {\em plain of characteristic $p$} if its roots are all plain of characteristic~$p$.
\item For $Q\in p^\BN$ and $w\in\BZ$, a {\em $Q$-Weil number of weight $w$} is an element $\alpha\in\overline\BQ$ such that (i) $\alpha$ is plain of characteristic $p$ and (ii) for any embedding $\iota\colon\overline\BQ\to\BC$ one has $|\iota\alpha|_\BC=Q^{w/2}$. We call $P\in\overline\BQ[T]$ {\em pure of weight $w$ for $Q$} if every root of $P$ is a $Q$-Weil number of weight~$w$.
\item If $\rho$ is a representation of a group $G$ and if $\phi\colon H\to G$ is a (natural) group homomorphism (for instance the inclusion of a subgroup), we write $\rho|_H$ for the restriction $\rho\circ\phi$ of $\rho$ to $H$. We use the same notation for the restriction of families of representations.
\item An $n$-dimensional $E$-rational compatible system $((\rho_\lambda)_{\lambda\in\CP_E'},(P_x)_{x\in|X|})$, or for short $\rho_\bullet$, consists of a homomorphism $\rho_\lambda\colon \pi_1(X)\to\GL_n(E_\lambda)$ for every $\lambda\in\CP_E'$ and a monic degree $n$ polynomial $P_x\in E[T]$, such that for each $(\lambda,x)\in \CP_E'\times |X|$ one has
\[\charpol_{\rho_\lambda(\Frob_x)}(T)=P_x(T)\]
with respect to the embedding $E\into E_\lambda$. An {\em $E$-rational compatible system} is $n$-dimensional for some $n$.
The system is called {\em pure of weight $w$}, if for every $x\in |X|$ the polynomial $P_x$ is pure of weight $w$ for $\#\kappa_x$. The system is called {\em tame}, if for any discrete rank $1$ valuation $v$ of $K=\kappa(X)$ and for any $\lambda\in\CP_E'$ the extension $(K^\sep)^{\kernel \rho_\lambda|_{\Gamma_K}}$ is tame at $v$. It is called {\em potentially tame} if for some finite cover $X'\to X$ the restriction $\rho_\bullet|_{\pi_1(X')}$ is tame.
\item For an $E$-rational $n$-dimensional compatible system $\rho_\bullet$, and $\lambda$ in $\CP_E'$ we denote by $G_{\rho_\bullet,\lambda}$, or simply $G_\lambda$, the Zariski closure of $\rho_\lambda(\pi_1(X))$ in $\GL_{n/E_\lambda}$, and we write $\alpha^{(o)}_\lambda$, respectively $\alpha^{(o)}_{\rho_\bullet,\lambda}$, for the inclusion of $G^{(o)}_\lambda$ into $\GL_{n/E_\lambda}$. If we consider $\pi_1^\geo(X)$ in place of $\pi_1(X)$, we add a superscript $\geo$ to the notation, e.g. $G^\geo_\lambda$ instead of $G_\lambda$ etc.
\item $M_{\rho_\bullet}$ (or simply $M$ if $\rho_\bullet$ is clear from context) is a motivic group for $\rho_\bullet$ defined over a finite extension~$F/E$; it is a connected reductive algebraic group, that is usually also split; it comes with a split motivic representation $\alpha\colon M\into \GL_{n/F}$; see \autoref{Subsect-Chin}.
\item For $\alpha\in\overline\BQ$ plain of characteristic $p$, and $L=\BQ(\alpha)$, we define the $L$-rational compatible system $\rho_{\alpha,\bullet}=(\rho_{\alpha,\lambda}\colon \Gamma_\kappa\to\GL_1(L_\lambda))_{\lambda\in\CP_L'}$ by $\rho_{\alpha,\lambda}(\Frob_\kappa)=\alpha$. If $\alpha$ is a $\#\kappa$-Weil number of weight $w$, then $\rho_{\alpha,\bullet}$ is pure of weight $w$. Via $\pi_1(X)\to\Gamma_\kappa$ we regard $\rho_{\alpha,\bullet}$ as a compatible system of representations of $\pi_1(X)$.
\item If $\tau\colon \pi_1(X)\to \GL_n(E)$ is a continuous representation (with the discrete topology on $E$), we denote by $\rho_{\tau,\bullet}$ the compatible $E$-rational system (pure of weight $0$) defined by $\rho_{\tau,\lambda}=((\GL_n(E)\into\GL_n(E_\lambda))\circ\tau)_{\lambda\in\CP_E'}$.
\end{itemize}
\begin{Rem}\label{Rem-OnNotation}
It is an elementary observation, cf.~\cite[last parag.\ in Sect.~2]{Chin04}, that for a compatible system $\rho_\bullet$, for all $x\in |X|$ the polynomials $P_x$ are plain of characteristic~$p$.
\end{Rem}
\begin{Con}\label{Conv-Semisimple} From \autoref{Sec-AI} on, all compatible systems are assumed to be semisimple.
\end{Con}

\section{Basic results on compatible systems over function fields}
\label{BasicNew}

This section recalls some standard results on compatible systems of Galois representations, tailored to the setting relevant to this work, where the domain is $\pi_1(X)$ or $\pi_1^\geo(X)$, respectively. In doing so, we fix some further notation, and we present extensions of existing results or give proofs of presumably known results that we could not locate in the literature.

\autoref{Subsec-CompSys} recalls change of coefficients and a kind of Weil restriction for compatible systems, and it recalls Serre's results on independence of $\lambda$ of $\pi_0$ of the monodromy groups. It also classifies $1$-dimensional compatible systems. \autoref{Subsect-Chin} recalls a result of Chin that gives an independence of $\lambda$ of the identity component of the monodromy group after a finite coefficient change. \autoref{Subsection-ReductionToCurve} explains a reduction from $X$ to a curve that preserves the monodromy groups (and the Galois images up to a uniform finite amount). Moreover it extends results from \cite{BGP15} on independence of images in the sense of Serre to arbitrary compatible systems. The last subsection recalls a semisimplicity result of Deligne and reduces investigations on monodromy groups attached to $\pi_1^\geo(X)$ to the case of semisimple compatible systems~$\rho_\bullet$ of $\pi_1(X)$.

\subsection{Compatible systems}
\label{Subsec-CompSys}
We collect some basic results for later use.
\begin{Lem}[{\cite[I.2.3, Theorem]{Serre-Abelian}}]\label{Lem-IsomOfCS}
Suppose that $\rho_\lambda$ and $\rho_\lambda'$ are semisimple representations $\pi_1(X)\to\GL_n(E_\lambda)$ with $\charpol_{\rho_\lambda(\Frob_x)}=\charpol_{\rho_\lambda(\Frob_x)}$ for all $x\in|X|$. Then $\rho_\lambda\cong\rho_\lambda'$. 

In particular, if two semisimple compatible systems $\rho_\bullet$ and $\rho_\bullet'$ have the same Frobenius polynomials $P_x(T),P'_x(T)\in E[T]$ for all $x\in|X|$, then $\rho_\bullet\cong\rho_\bullet'$, i.e, for all $\lambda\in \CP_E'$ one has $\rho_\lambda\cong\rho_\lambda'$.
\end{Lem}

\begin{Def}\label{Def-CoeffExtension}
If $E'$ is a finite extension of $E$, we define the {\em coefficient extension} $\rho_\bullet\otimes_EE'$ of $\rho_\bullet$ to $E'$ by 
\[(\rho_\bullet\otimes_EE')_{\lambda'} := (\GL_n(E_\lambda)\into \GL_n(E'_{\lambda'}))\circ \rho_\lambda\]
for all $\lambda'\in\CP_{E'}'$ with contraction $\lambda$ to $\CP_E'$.
\end{Def}

\begin{Lem}\label{Lem-CoeffExtOfCompSys} 
For a finite extension $E'/E$ and any $\lambda'\in\CP_{E'}'$ above $\lambda\in\CP_E'$ one has \[G_{\rho_\bullet\otimes_E E',\lambda'}=G_{\rho_\bullet}\otimes_{E_\lambda}E'_{\lambda'}.\]
\end{Lem}
\begin{proof}
The assertion of the lemma follows from the following well-known result: Let $L\supset K$ be an extension of fields field, let $T$ be a $K$-variety and $T_L$ be its base change to $L$. Then the Zariski closure in $T_L$ of any subset $S\subset T(K)$ is the base change under $K\to L$ of the Zariski closure of $S$ in $T$, i.e., the operations base change and Zariski closure commute for varieties.
\end{proof}

Let $E'$ be a subfield of $E$. Then the Weil restriction $\Res_{E/E'}\GL_{n/E}$ can be regarded as a closed subgroup of $\GL_{n[E:E']/E'}$, via the natural transformation functor from $E'$-algebras $R'$ given by
\[\Res_{E/E'}\GL_{n/E}(R') = \GL_n(E\otimes_{E'} R') \hookrightarrow \GL_{n[E:E']}(R'),  \]
where the equality on the right is given by choosing a basis for $E$ over $E'$. 
\begin{Def}\label{Def-WeilRestriction}
The {\em Weil restriction $\Res_{E/E'} \rho_\bullet$} of $\rho_\bullet$ to $E'$ is defined by assigning to any $\lambda'\in\CP_{E'}'$ the representation
\[(\Res_{E/E'}\rho_\bullet)_{\lambda'} := \bigoplus_{\lambda|\lambda'} \rho_\lambda  \colon \pi_1(X) \to \big(\Res_{E/E'}\GL_{n/E}\big)(E'_{\lambda'}) = \prod_{\lambda|\lambda'} \GL_n({E_{\lambda}}),\]
where the product is over all $\lambda\in\CP_E'$ above $\lambda'$. It is an $E'$-rational compatible system. 
\end{Def}
An important example of a Weil restriction is $\Res_{E/\BQ} \rho_\bullet$; cf.~\autoref{Thm-SerreOpen}. For fields $E,E'\supset E_0$, the set of field homomorphism $\sigma\colon E\to E'$ with $\sigma|_{E_0}=\id_{E_0}$ is $\Hom_{E_0}(E,E')$.
\begin{Lem}\label{Lem-WeilResAndMonodromy}
Suppose that $\rho_\bullet$ is $E$-rational. Let $E'\subset E$ be a subfield, and $F/E$ an extension such that $F/E'$ is finite Galois. Then the following hold:
\begin{enumerate}
\item If we denote by $\ublank \otimes_E^\sigma F$ the tensor product of $F$ over $E$ with respect to $\sigma\colon E\to F$, then
\[ \Res_{E/E'} \rho_\bullet\otimes_{E'} F\cong\bigoplus_{\sigma\in\Hom_{E'}(E,F)} \rho_\bullet\otimes_E^\sigma F.\]
\item If $\rho_\bullet$ is semisimple, then so is $\Res_{E/E'} \rho_\bullet$.
\item If $G_{\rho_\bullet,\lambda}$ is semisimple for all $\lambda\in\CP_E'$, then $G_{\Res_{E/E'} \rho_\bullet,\lambda'}$ is semisimple for all $\lambda'\in\CP_{E'}'$.
\end{enumerate}
\end{Lem}
\begin{proof}
The proof of (a) is immediate from the definitions. We deduce (b) from (a). For this recall that a representation $V$ of a profinite group $\Gamma$ over $E'$ is semisimple, if $V\otimes_{E'}F$ is semisimple as a representation of $\Gamma$ over $F$: To see this, consider a $\Gamma$-stable submodule $W\subset V$ and a $\Gamma$-equivariant splitting $s\colon V\otimes_{E'}F\to W\otimes_{E'}F$ of $W\otimes_{E'}F\into V\otimes_{E'}F$. Then $\frac1{[F:E']}\sum_{\tau\in\Gal(F/E')}\tau (s)\colon W\to V$ is a $\Gamma$-equivariant splitting of $W\into V$. Thus it suffices to show the semisimplicity of the right hand side of (a). For this note that since $\rho_\bullet$ is semisimple, then so is $\rho_\bullet\otimes^\sigma_{E}F$ (cf.~\cite[Cor.~69.8]{CurtisReinerOld}), and thus also the direct sum over all $\sigma\in\Hom_{E'}(E,F)$.

For (c) let $\mu$ be in $\CP_F'$, let $\lambda'\in\CP_{E'}'$ be below $\mu$, and denote for $\sigma\in\Hom_{E'}(E,F)$ by $\mu|\sigma$ the place of $\CP_E'$ under $\mu$ with respect to $\sigma$. By (a) the group $H':=G_{\Res_{E/E'}\rho_\bullet,\lambda}\otimes_{E'_\lambda} F_\mu$ is a closed subgroup of $H:=\prod_{\sigma\in\Hom_{E'}(E,F)}G_{\rho_\bullet,\mu|\sigma}\otimes^\sigma_{E_{\mu|\sigma}}F_\mu$. Let $\pr_\sigma \colon H\to H_\sigma:=G_{\rho_\bullet,\mu|\sigma}\otimes^\sigma_{E_{\mu|\sigma}}F_\mu$ be the canonical projection. Let $R$ be the radical of $H'$. By  \autoref{Lem-CoeffExtOfCompSys}, it will suffice to show that $R$ is trivial. By our hypothesis and \autoref{Lem-CoeffExtOfCompSys}, the groups $H_\sigma$ are semisimple. Since $\pr_\sigma(H')$ is a closed subgroup of $H_\sigma$, see~\cite[2.2.5]{Springer}, it must be all of $H_\sigma$, by the definition of $H'$ and $H_\sigma$. Now because $\pr_\sigma\colon H'\to H_\sigma$ is surjective and $H_\sigma$ is semisimple, we have $R\subset \kernel\pr_\sigma$ for all $\sigma$. But then $R$ lies in the kernel of the injection $\prod_\sigma\pr_\sigma$, and hence $R$ must be trivial.
\end{proof}
\begin{Rem}
For $H$ and $H'$ as in the previous proof, without any further hypotheses there is not more we can say than what we explained above. \autoref{Prop-GaloisConjugates} gives examples with $H'=H$. But from  \autoref{Prop-GaloisConjugates} and its proof one can also construct examples where in fact $H'\cong G_{\rho_\bullet,\mu|\sigma_0}\otimes_{E_{\mu|\sigma_0}}F_\mu$ for any $\sigma_0\in \Hom_{E'}(E,F)$. \autoref{Sec-OpenImage} will give further results about possible~$H'$.
\end{Rem}
\begin{Thm}[{\cite[Prop.~2 and its Cor.]{Serre-Driebergen}}]\label{Thm-SerreOpen}
Let $\Gamma$ be a compact subgroup of $\GL_n(\BQ_\ell)$ with Zariski closure $G$. If $G^\circ$ is semisimple, then $\Gamma$ is open in $G(\BQ_\ell)$. In particular, if $\rho_\bullet$ is a $\BQ$-rational compatible system such that $G_{\rho_\bullet,\ell}$ is semisimple, then $\rho_\ell(\pi_1(X))$ is open in $G_{\rho_\bullet,\ell}(\BQ_\ell)$.
\end{Thm}

\begin{Ex}\label{Ex-E-CS-WithoutOpenImage}
Clearly, in \autoref{Thm-SerreOpen} one cannot replace $\BQ_\ell$ by a proper finite extension $L$: Let $\Gamma$ and $G$ be as in \autoref{Thm-SerreOpen} and suppose that $G^o$ is semisimple and non-trivial. Then $\Gamma$ is Zariski dense in $G\otimes_{\BQ_\ell}L$ by  \autoref{Lem-CoeffExtOfCompSys}. But $G(\BQ_\ell)$ is not open in $G(L)$, and so $\Gamma$ is not open in $G(L)$. For positive results, see~\autoref{Thm-OpennessAboveEll2}(e) and \autoref{Cor-HuiLarsen}.
\end{Ex}

The following apparently well-known result we learned from \cite[Thm.~1.6]{Chin04}:
\begin{Lem}\label{Lem-GDer=Ggeo}
Let $\rho_\lambda\colon \pi_1(X)\to \GL_n(E_\lambda)$ be semisimple. Denote by $G\subset \GL_{n,E_\lambda}$ its Zariski closure and by $G^\geo\subset G$ that of $\rho_\lambda(\pi^\geo_1(X))$. Then $G^\der=G^{\geo,o}$ and $G^\der$ is semisimple.
\end{Lem}
\begin{proof}
For lack of a reference, we include a proof. Recall that by our conventions from \autoref{Sec-Notation}, we denote by $G^\der$ the derived group of the identity component of $G$, hence it is connected. The group $\pi_1^\geo(X)$ is normal in $\pi_1(X)$, and thus so is $G^{\geo}\subset G$, and then also $G^{\geo,o}$ in $G^o$. The latter groups are connected reductive, because $\rho_\lambda$ is semisimple by hypothesis, and $\rho_\lambda|_{\pi^\geo_1(X)}$ is so because of a result of Deligne -- see \autoref{Thm-GeomSemisimplicity}. Hence the quotient $\bar G:=G^o/G^{\geo,o}$ is reductive and connected. It is also abelian because the image of $\pi_1(X)/\pi_1^\geo(X)\cong\hat\BZ$ is dense in it. It follows that $\bar G$ is a torus. This proves the first assertion; cf.~\cite[Sec.~2]{LarsenPink95}. The second follows again from \autoref{Thm-GeomSemisimplicity}.
\end{proof}
For later use, we record the following consequences of \autoref{Lem-GDer=Ggeo} and \autoref{Thm-SerreOpen}.
\begin{Cor}\label{Cor-Q-CS-OpenImage}
If $\rho_\bullet$ is $\BQ$-rational, then $\rho_\ell(\pi_1^\geo(X))\subset G^{\geo}_{\rho_\bullet,\ell}(\BQ_\ell)$ is open for all $\ell\in\CP_\BQ'$.
\end{Cor}

\begin{Cor}\label{Cor-GeomAlgFiniteIndex}
Suppose $\rho_\bullet$ is $E$-rational semisimple and the groups $G_{\rho_\bullet,\lambda}$ are semisimple for all $\lambda\in\CP_E'$. Then for all $\ell\in\CP_\BQ'$, the subgroup 
$\big(\prod_{\lambda|\ell} \rho_\lambda\big)(\pi_1^\geo(X))\subset\big(\prod_{\lambda|\ell}\rho_\lambda\big)(\pi_1(X))$ is open.
\end{Cor}
\begin{proof}
Using \autoref{Lem-WeilResAndMonodromy}(c) and the Weil restriction $\Res_{E/\BQ}\rho_\bullet$, we may assume $E=\BQ$. Because $G_{\rho_\bullet,\ell}$ is semisimple, we deduce $G^{\geo}_{\rho_\bullet,\ell} =G_{\rho_\bullet,\ell}$ from \autoref{Lem-GDer=Ggeo}. The result now follows from \autoref{Thm-SerreOpen}.
\end{proof}
\begin{Rem}\label{Rem-AdelicImprovement}
In \autoref{Cor-GeomImageOpenInArithImage0} we shall give an adelic refinement of \autoref{Cor-GeomAlgFiniteIndex}.
\end{Rem}

Regarding the discrepancy between $G_\lambda$ and $G_\lambda^o$ one has the following result due to Serre:
\begin{Thm}[{\cite[p.~15--20]{Serre-OeuvresIV} or \cite[Prop.~6.14]{LarsenPink92}\footnote{The references only consider $E$-rational compatible systems for $E=\BQ$. The adaptation to general $E$ is minor.}}] \label{Thm-SerreOnComponents}
The kernel of \[\rho_\lambda\mod{G_\lambda^o}\colon\pi_1(X)\to G_\lambda/ G_\lambda^o\] is finite and independent of $\lambda\in\CP_E'$.
\end{Thm}
\begin{Def}
We say that $\rho_\bullet$ is {\em connected} or has {\em connected monodromy} if $G_\lambda^o=G_\lambda$ for one, and hence, by \autoref{Thm-SerreOnComponents}, for all~$\lambda\in\CP_E'$.
\end{Def}
By \autoref{Thm-SerreOnComponents} for any compatible system $\rho_\bullet$ there exists a finite \'etale cover $X'\to X$ such that $\rho_\bullet|_{\pi_1(X')}$ has connected monodromy. By the following result, whose proof is a simple exercise, restriction to $\pi_1(X')$ leaves $G_\lambda^o$ invariant:
\begin{Lem}\label{Prop-ChinUnderRestr}
Let $X'\to X$ be any finite cover. Then $G^o_{\rho_\bullet|_{\pi_1(X')},\lambda}=G^o_{\rho_\bullet,\lambda}$ for all $\lambda\in\CP_E'$.
\end{Lem}

Below we  need the following immediate consequence of \cite[Thm.~1]{KatzLang} by Katz and Lang:
\begin{Thm}\label{Thm-KatzLang}
Denote by $\pi_1(X)^{\ab,p}$ the maximal abelian quotient of $\pi_1(X)$ modulo its pro-$p$ subgroup. Then the following hold:
\begin{enumerate}
\item For any $x\in X$, the cokernel of the canonical map $\pi_1(x)\to \pi_1(X)^{\ab,p}$ is finite.
\item There exists a finite extension $X'\to X$, say with field of constants $\kappa'$, such that the induced homomorphism $\pi_1(X')\to\pi_1(X)^{\ab,p}$ factors via the canonical map $\pi_1(X')\to\Gamma_{\kappa'}$.
\end{enumerate}
\end{Thm}

The following results cover $1$-dimensional representations. We first need some notation. For $\alpha\in\overline\BZ_\ell^\times$ and $\can\colon\pi_1(X)\to\Gamma_\kappa$ the canonical map, define $\rho_\alpha\colon \pi_1(X)\to\GL_1(\overline\BQ_\ell)$ as the composite continuous homomorphism
\begin{equation} \label{eq:RhoAlpha}
\pi_1(X)\,\stackrel\can\longto \, \Gamma_\kappa\stackrel{r_\alpha}\longto \GL_1(\BQ_\ell)\hbox{ determined by }r_\alpha(\Frob_\kappa)=\alpha.\end{equation}
\begin{Lem}\label{Lem-AbelianCS}
Let $\rho\colon\pi_1(X)\to\GL_1(\overline \BQ_\ell)$ be a continuous representation and define $\alpha_x=\rho(\Frob_x)$ for all $x\in |X|$. Then the following hold:
\begin{enumerate}
\item There exists $\alpha\in\overline\BZ_\ell^\times$ such that $\rho':=\rho\otimes(\rho_{\alpha})^{-1}$ has finite image.
\item If all $\alpha_x$ lie in $\overline \BQ$ and are plain of characteristic $p$, then $\alpha$ in (a) lies in $\overline\BQ$ and is plain of characteristic $p$. In this case, if $L=\BQ(\alpha)$ and $\mu\in\CP_L'$ is the place induced from $\overline\BQ\into\overline\BQ_\ell$, then $\rho_\alpha=(\rho_{\alpha,\mu})\otimes_{L_\mu}\overline\BQ_\ell$ is a member of a compatible system.
\item  If all $\alpha_x$ are $\#\kappa_x$-Weil numbers of weight $w$, then $\alpha$ in (a) is a $q$-Weil number of weight~$w$.
\item If all $\alpha_x$ lie in a number field $E$, then $\alpha$ in (a) can be chosen to lie in $E$.
\end{enumerate}
\end{Lem}
\begin{proof}
For $x\in |X|$ let $d_x$ be the degree of $\kappa_x$ over $\kappa$. Because $X$ is geometrically irreducible over the finite field $\kappa$, we can find $x_1,x_2\in |X|$ such that $\gcd(d_{x_1},d_{x_2})=1$.\footnote{\label{FootnoteDegreeOne}Lacking a reference, we give a short proof: Let $\kappa_r$ be the finite extension of $\kappa$ of degree $r$. Fix two rational prime numbers $\ell_1\neq \ell_2$. Because $X$ is geometrically irreducible, suitable estimates based on the Weil conjectures, see \cite[Thm.~2.1]{zywina-israel}, imply $\lim_{j\to \infty} X(\kappa_{\ell_i^j})=\infty$ for $i\in \{1, 2\}$. Thus we can find $x_i\in |X|$ with $d_{x_i}$ a power~of~$\ell_i$.} Let $m_1, m_2\in\mathbb{Z}$ be such that $1=m_1 d_{x_1}+m_2 d_{x_2}$ and set $\alpha=\alpha_{x_1}^{m_1}\alpha_{x_2}^{m_2}$. Since $\rho$ is continuous and $\pi_1(X)$ is compact all $\alpha_x$ lie in $\overline\BZ{}_\ell^\times$, and hence so does $\alpha$. This choice also proves (d) once the proof of (a) is complete.

To complete (a), note first that the natural map $\Gamma_{\kappa_x}\to \Gamma_{\kappa}$ is given by $\Frob_x\mapsto \Frob_{\kappa}^{d_x}$. By \autoref{Thm-KatzLang} there exists $D\in\BN$ such that $\alpha_x^{Dd_{x'}}=\alpha_{x'}^{Dd_x}$ for any $x,x'\in |X|$. A direct computation now shows that $(\rho')^D$ is the trivial character (on all $\Frob_x$), and hence the image is finite, proving~(a) (and (d)). 

For (b) and (c) we only have to show that, independently of the construction of $\alpha$ in (a), it will have the properties asserted there. The remaining part of (b) is clear from the definition of $\rho_{\alpha,\bullet}$. Let therefore $D$ denote the order of $\rho'$. Evaluating $(\rho')^D$ at $\Frob_x$ shows $\alpha^{d_xD}=\alpha_x^D$, i.e., that $\alpha$ is one of the $Dd_x$-th roots of $\alpha_x$. This implies (b) and (c) since the properties considered there are stable under taking roots.
\end{proof}

\begin{Rem}
The proof also shows that if one $\alpha_x$ is plain of characteristic $p$, or a $\#\kappa_x$-Weil number pure of weight $w$, then the same holds for all $\alpha_x$.
\end{Rem}
\begin{Cor}[{cf.~\cite[4.2, 4.3]{Chin04}}]\label{Prop-AbelianCS}
Let $\rho_\bullet$ be a $1$-dimensional $E$-rational compatible system. Then the following~hold:
\begin{enumerate}
\item There exists $\alpha\in E^\times$ plain of characteristic $p$ and a continuous representation $\tau\colon\pi_1(X)\to\GL_1(E)$ such that $\rho_\bullet=\rho_{\alpha,\bullet}\otimes\rho_{\tau,\bullet}$.
\item If $\rho_\bullet$ is pure of weight $w$, then $\alpha$ in (a) is a $q$-Weil number of weight~$w$.
\item There exists a finite \'etale cover $X'\to X$ and some $\alpha\in E^\times$ which is plain of characteristic $p$ (and a $\#\kappa'$-Weil number of weight $w$ if $\rho_\bullet$ is pure of weight $w$), such that $\rho_\bullet|_{\pi_1(X')}$ is isomorphic to $\rho_{\alpha,\bullet}|_{\pi_1(X')}$, and in particular factors via $\Gamma_{\kappa'}$ for $\kappa'$ the constant field of $X'$.
\end{enumerate}
\end{Cor}

\begin{proof} 
Because $\rho_\bullet$ is compatible $1$-dimensional, $\rho_\lambda(\Frob_x)$ can be identified with some $\alpha_x\in E^\times$ that is independent of $\lambda$. In particular, $\alpha_x$ lies in $\CO_\lambda$ for all $\lambda\in\CP_E'$ and is thus plain of characteristic~$p$; if $\rho_\bullet$ is pure of weight $w$, then moreover $\alpha_x$ is a $\#\kappa_x$-Weil number of weight $w$. Fix one $\lambda_0\in\CP_E'$ and an embedding $E_{\lambda_0}\into\overline\BQ_{\ell_\lambda}$. Let $\alpha\in E$ be as in \autoref{Lem-AbelianCS} for $\rho_{\lambda_0}\otimes_{E_{\lambda_0}}\overline\BQ_{\ell_{\lambda_0}}$. It follows that $\rho_{\alpha,\bullet}^{-1}\otimes\rho_\bullet$ is a compatible system that is of finite order at $\lambda_0$.

Then by \autoref{Thm-KatzLang} (or by \autoref{Thm-SerreOnComponents}) there exists a finite cover $X''\to X$ such that $\rho'_\bullet|_{\pi_1(X'')}$ is trivial. Now one can simply apply complex representation theory of finite groups to deduce the existence of $\tau\colon \pi_1(X)\to\GL_1(E)$ such that (a) holds. Part (b) follows from \autoref{Lem-AbelianCS}(c). Moreover (c) follows by taking for $X'$ the cover of $X$ that corresponds to the kernel of $\tau$ from~(a).
\end{proof}

\subsection{The motivic group}
\label{Subsect-Chin}

For repeated later use, we first recall two basic facts on the representation theory of split semisimple connected groups in characteristic zero, and state an immediate corollary:
\begin{Thm}[{\cite[IV.3.3.3]{Demazure-Gabriel}, \cite[Thm. 2.5]{tits1971}}]   \label{Thm-RepTheory}
Suppose that $M$ is a split reductive connected group over a field $F$ of characteristic $0$. Then the following hold.
\begin{enumerate}
\item The category of finite dimensional representations of $M$ is semisimple.
\item The irreducible representations of $M$ are in bijection with highest weight representations, parameterized by the dominant weights of $M$.
\end{enumerate}
\end{Thm}
\begin{Cor}\label{Cor-RepTheory}
Let $M$ and $F$ be as in \autoref{Thm-RepTheory}, and let $F'$ be any extension of $F$. Then $V\mapsto V\otimes_FF'$ sets up a bijection between the isomorphism classes of irreducible, finite-dimensional representations of $M$ and those of $M\otimes_FF'$.
\end{Cor}

\begin{Def}\label{Def-ChinTriple}
Let $\rho_\bullet$ be a semisimple $E$-rational $n$-dimensional compatible system. We call a pair $(M,\alpha)$ consisting of
\begin{enumerate}
\item[(i)] a connected reductive algebraic group $M$ over $E$,
\item[(ii)] and a faithful $E$-rational representation $\alpha\colon M\to \GL_n$
\end{enumerate}
a {\em motivic pair for $\rho_\bullet$}, $M$ the {\em motivic group of $\rho_\bullet$} and $\alpha$ the {\em motivic representation of $\rho_\bullet$}, if the following conditions hold: For any $\lambda\in\CP_E'$ one has an isomorphism of algebraic groups
\begin{equation}\label{Eqn-ChinGroups1}
\iota_\lambda\colon  M\otimes_EE_\lambda\stackrel\simeq\longto G^o_{\rho_\bullet,\lambda},
\end{equation}
such that the following diagram is commutative up to conjugation
\begin{equation}\label{Eqn-ChinReps1}
\xymatrix@R-.5pc@C-.6pc{M\otimes_{E} E_\lambda\ar[rr]^-{i_\mu}\ar[dr]_{\alpha\otimes_E E_\lambda} && G_{\rho_\bullet, \lambda}^\circ \ar[dl]^{\mathrm{incl}}\\
&GL_{n, E_\lambda} \rlap{;}&  }
\end{equation}
i.e., there is an isomorphism of representations $M\to \GL_n$ under the identification in \autoref{Eqn-ChinGroups1}.

We call $(F,M,\alpha)$ a {\em motivic triple for $\rho_\bullet$} if $F$ is a finite extension of $E$ and if $(M,\alpha)$ is a motivic pair for $\rho_\bullet\otimes_EF$. The motivic triple is called {\em split} if $M$ is split.
\end{Def}

\begin{Rem}\label{Rem-OnChin}
Let $\rho_\bullet$ be a semisimple $E$-rational $n$-dimensional compatible system. Then for split motivic triples $(F,M,\alpha)$ for $\rho_\bullet$, the root datum and the weights defining $(M,\alpha)$ are independent of $F$. Therefore we shall sometimes simply speak of the {\em split motivic group} and the {\em split motivic representation}, assuming that a sufficiently large $F$ has been chosen.

For the same reason, if $M$ is split, then the pair $(M,\alpha)$ is uniquely determined by the properties \autoref{Eqn-ChinGroups1} and  \autoref{Eqn-ChinReps1}  for a single $\mu$.
\end{Rem}
\begin{Rem}\label{Rem-OnChin2} 
\begin{enumerate}
\item The definition of the motivic group $M$ should also include conditions at the places of $E$ above $p$; this is possible using the work \cite{Abe} of Abe that attaches certain isocrystals at places above $p$ and thus extends the work of L.\ Lafforgue; cf.~also \cite{Drinfeld-ProSemisimple}. For the present work these places are irrelevant, and so we omit them.
\item It is expected that the motivic group exists (over $E$). The existence of $\alpha$ over $E$ depends on a Brauer obstruction; cf.~\cite[E.9]{Drinfeld-ProSemisimple}. But ignoring the places above $p$, it might well be possible that $\alpha$ can be realized over $E$ also. See \cite{Hui-Motivic} for recent progress.
\end{enumerate}
\end{Rem}

We shall use of the following result of Chin, generalizing \cite[Thm.~2.4]{LarsenPink95}, at various places:
\begin{Thm}[{\cite[Thm.~1.4 and Thm.~1.6]{Chin04}}]\label{Thm-MainThmChin}
For any pure semisimple $E$-rational $n$-dimensional compatible system $\rho_\bullet$ there exists a split motivic triple $(F,M,\alpha)$ for $\rho_\bullet$. Moreover (cf.~\autoref{Lem-GDer=Ggeo}) for any $\mu\in\CP_F'$ with contraction $\lambda\in\CP_E'$ the isomorphism \autoref{Eqn-ChinGroups1} (for $\rho_\bullet\otimes_EF$) induces an isomorphism of algebraic groups
\begin{equation}\label{Eqn-ChinGroups2}
G^o_{\rho_\bullet|_{\pi_1^\geo(X)},\lambda}\otimes_{E_\lambda}F_\mu\cong M^\der\otimes_FF_\mu,
\end{equation}
and the diagram \autoref{Eqn-ChinReps1} an isomorphisms of representations $M^\der\to \GL_n$
\begin{equation}\label{Eqn-ChinReps2}
\alpha|_{M^\der}\otimes_FF_\mu \cong \alpha^o_{\rho_\bullet|_{\pi_1^\geo(X)},\lambda}\otimes_{E_\lambda}F_\mu .
\end{equation}
\end{Thm}

\begin{Rem}\label{Rem-OnChin3}
It cannot be expected that in \autoref{Thm-MainThmChin} there is a smallest or optimal choice for $F$, since $M$ is split over $F$. This is analogous to the well-known fact that for representations of finite groups over $\overline\BQ$ there is in general no smallest field of definition -- because of the Schur index. A simple example in the present case would be a system $\rho_\bullet$ of dimension $4$ with $G_{\lambda_0}$ a quaternion division algebra for some $\lambda_0\in \CP_E'$ and $G_\lambda\cong\GL_{2,E_\lambda}$ for all other $\lambda\in\CP_E$. Then there could be many quadratic extensions of $E$ over which a split motivic group exists, but over $E$ this is not possible.
\end{Rem}

\begin{Lem}\label{Lem-FactsOnChin}
Suppose $\rho_\bullet$ is $E$-rational and semisimple with motivic pair $(M,\alpha)$. Then:
\begin{enumerate}
\item If $E'/E$ is a finite extension, then $(M\otimes_EE',\alpha\otimes_EE')$ is a motivic pair for $\rho_\bullet\otimes_EE'$.
\item If $E'\subset E$ is a subfield and if $\rho_\bullet$ is pure, 
then $\Res_{E/E'}\rho_\bullet$ has a motivic triple $(F',M',\alpha')$, and if $M$ is semisimple, then so is~$M'$.
\item If $X'\to X$ is a finite cover, then $(M,\alpha)$ is a motivic pair for $\rho_\bullet|_{\pi_1(X')}$.
\end{enumerate}
\end{Lem}
\begin{proof} Part (a) is immediate from \autoref{Lem-CoeffExtOfCompSys}, and Part (c) is immediate from \autoref{Prop-ChinUnderRestr}. We 
shall now prove Part (b). $\Res_{E/E'}\rho_\bullet$ is a semisimple compatible system by \autoref{Lem-WeilResAndMonodromy}(b), and it is easy to
see that $\Res_{E/E'}\rho_\bullet$ is pure (of the same weight as $\rho_\bullet$). By \autoref{Thm-MainThmChin} it follows that 
$\Res_{E/E'}\rho_\bullet$ has a motivic triple $(F',M',\alpha')$. If $M$ is semisimple in addition, then all $G_{\rho_\bullet, \lambda}$ are semisimple, thus
all $G_{\Res_{E/E'}\rho_\bullet, \lambda'}$ are semisimple as well by \autoref{Lem-WeilResAndMonodromy}(c), and this implies that $M'$ is semisimple. 
\end{proof}

\begin{Rem}\label{Rem-OnChin4}
In \autoref{Prop-GenChinTriple} and \autoref{Cor2-GenChinTriple} we shall remove the purity hypothesis in \autoref{Thm-MainThmChin} and \autoref{Lem-FactsOnChin}(b), respectively.
\end{Rem}

\begin{Rem}
In \autoref{Lem-FactsOnChin}(b), for $\lambda\in\CP_E'$ and $\lambda'\in\CP_{E'}'$ the contraction of $\lambda$, the group $G_{\Res_{E/E'}\rho_\bullet,\lambda'}$ is a closed subgroup of $\Res_{E_\lambda/E'_{\lambda'}}G_{\rho_\bullet,\lambda}$. In general the containment is proper; for instance for $\rho_\bullet=\rho'_\bullet\otimes_{E'}E$ with $\rho'_\bullet$ an $E'$-rational compatible system.
\end{Rem}

The proof of \autoref{Thm-MainThmChin} crucially relies on the work of L.\ Lafforgue on the global Langlands correspondence. Since in later parts we shall need to directly apply his work, we now state it in a form to be used later; our formulation also requires \cite[Thm.~4.1]{Chin04}.
\begin{Thm}[{\cite[Thme.~(p.~2), Thme.~VI.9, Prop.~VII.4]{Lafforgue}}]\label{Thm-Lafforgue}
Let $X$ be a smooth projective curve over $\kappa$. Then the following holds:
\begin{enumerate}
\item Suppose $\Pi$ is a cuspidal automorphic representation of $\GL_{n/\BA_{\kappa(X)}}$ with finite order central character $\tau\colon \pi_1(X)\to\GL_1(\overline \BQ)$, conductor $N$ which is a divisor on $X$, and Hecke field $E$. For $x\in |X|\setminus\supp N$ let $p_{\Pi,x}\in E[T]$ be the Hecke polynomial at $x$. Then there exists a finite extension field $E'$ of $E$ and an $E'$-rational compatible system
\[(\rho_{\Pi,\lambda}\colon\pi_1(X\setminus\supp N)\to\GL_n(E'_\lambda))_{\lambda\in\CP_{E'}'},\] 
pure of weight $0$ with determinant $\rho_{\tau,\bullet}$, such that for all $\lambda\in\CP_{E'}'$ and $x\in |X|\setminus\supp N$ the Frobenius polynomial of $\rho_{\Pi,\bullet}$ is equal to $p_{\Pi,x}$. Moreover each $\rho_{\Pi,\lambda}$ is absolutely irreducible.
\item If $\rho\colon \pi_1(X\setminus S)\to \GL_n(\overline{\BQ_\ell})$ for some finite $S\subset|X|$ is continuous and absolutely irreducible and if $\det\rho(\pi_1(X))$ is finite, then there exists a cuspidal automorphic representation $\Pi$ of $\GL_{n/\BA_{\kappa(X)}}$ with finite central character such that one has
\[\rho=\iota\circ\rho_{\Pi,\lambda'}\]
for some $\lambda'\in\CP_{E'}'$ and embedding $E'_{\lambda'}\to\overline{\BQ_\ell}$, with $E'$ from part (a). In particular, $\rho_{\Pi,\lambda'}$ is a member of a compatible system that is pure of weight~$0$.
\end{enumerate}
Moreover the bijection set up above is also compatible with the local Langlands correspondence at all places dividing the conductor. In particular, the conductors of $\rho_{\Pi,\lambda}$ and of $\Pi$ coincide.
\end{Thm}

\subsection{Reduction to the curve case and almost independence}
\label{Subsection-ReductionToCurve}
This subsection collects two deeper results on compatible systems needed in the later parts of this work. The first result explains why for a given compatible system $\rho_\bullet$, after a finite base change $X'\to X$ all ramification of $\rho_\bullet|_{\pi_1(X')}$ is tame, and in fact pro-$\ell_\lambda$ for any $\rho_\lambda$, and why for tame systems one can find a curve $C$ embedded into $X$, say via $\iota$, such that for all $\lambda\in\CP_E'$ one has $\rho_\lambda(\pi_1(X))=\rho_\lambda(\iota_*(\pi_1(C)))$. The second result is an extension of the main result of \cite{BGP15}. We show the almost independence in the sense of \cite{Serre-Independence} of arbitrary $\BQ$-rational compatible systems of certain representations of $\pi^\geo_1(X)$. 

\medskip

We say that a profinite Galois extension $L$ of $K=\kappa(X)$ has at most pro-$\ell$ ramification if for any discrete rank $1$ valuation $v$ of $K$ with extension $w$ to $L$ and completions $K_v$ and $L_w$, the inertia subgroup in $\Gal(L_w/K_v)$ is a pro-$\ell$ group.

\begin{Lem}\label{Lem-OnTamingCover}
Let $\rho_\bullet$ be a compatible system. Then there exists a finite \'etale cover $X'\to X$ such that for any $\lambda\in\CP_E'$ the ramification of $\rho_\lambda|_{\pi_1(X')}$ is pro-$\ell_\lambda$. In particular, $\rho_\bullet$ restricted to $\pi_1(X')$ is tame.
\end{Lem}
\begin{proof}
The argument is a slight modification of the argument given in \cite[\S~4 and \S~6]{BGP15}. Let $\rho_\bullet$ be the $E$-rational compatible system. Choose a finite cover $X'\to X$, such that for two places $\lambda_0,\lambda_1\in\CP_E'$ with $\ell_{\lambda_0}\neq\ell_{\lambda_1}$ the image of $\rho_{\lambda_0}(\pi_1(X'))$ is pro-$\ell_{\lambda_0}$ and that of $\rho_{\lambda_1}(\pi_1(X'))$ is pro-$\ell_{\lambda_1}$. Consider any embedding $\iota\colon C\to X'$ of a curve into $X'$, and let $\iota_*\colon\pi_1(C)\to\pi_1(X)$ be the corresponding homomorphism of fundamental groups. Then $\rho_\bullet\circ\iota_*$ is an $E$-rational compatible system of representations of $\pi_1(C)$. For any place $c$ of $K=\kappa(C)$ denote by $(r_{\lambda,c},N_{\lambda,c})$ the associated Weil-Deligne representation, and by $I_c$ the inertia subgroup of $\Gamma_{K_c}$. By \cite[Thm.~9.8]{DeligneEquations} the family of representations $r_{\lambda,c}|_{I_c}$ is independent of $\lambda$. By our hypotheses, the image is an $\ell_{\lambda_0}$ and an $\ell_{\lambda_1}$-group, and hence trivial. It follows that all ramification of $\rho_\lambda\circ\iota_*$ is unipotent, and hence at most pro-$\ell_\lambda$. Then \cite[Prop.~4.6]{BGP15}, which builds on work of Kerz-Schmidt-Wiesend, cf.~\cite[Thm.~4.4]{KerzSchmidt}, completes the proof of the lemma.
\end{proof}

We also need the following result due to Drinfeld:
\begin{Thm}[{\cite[Prop.~2.17 and~2.18]{Drinfeld12}}]\label{Drinfeld-OnPi3}
Let $Y$ be a normal geometrically connected variety over $\kappa$ and let $U\subset Y$ be a dense regular open subvariety. Then there exists a smooth irreducible curve $C$ on $U$ such that the canonical map $\pi_1(C)\to \pi_1(U)/H$ is surjective, where $H \subset\pi_1(U)$ is defined as follows:
\begin{enumerate}
\item $H$ is a closed normal subgroup such that $\pi_1(U)/H$ contains an open pro-$\ell$ subgroup, or
\item $Y$ is projective and $H$ is the kernel of $\pi_1(U)\to \pi_1^\tame(U)$, where $\pi_1^\tame(U)$ is the quotient of $\pi_1(U)$ that classifies covers of $U$ which are at most tamely ramified at codimension $1$ divisors of $Y\setminus U$. 
\end{enumerate}
Moreover one can assume that $C$ passes through any finite number of points of $U$, and in particular one may, cf.~Footnote~\ref{FootnoteDegreeOne}, assume that $C$ is geometrically irreducible.
\end{Thm}

The following corollary is an immediate consequence of \autoref{Lem-OnTamingCover} and \autoref{Drinfeld-OnPi3}(b).
\begin{Cor}\label{Cor-ReductionToCurve}
For any compatible system $\rho_\bullet$, there exists a finite \'etale cover $\pi\colon X'\to X$ and a curve $C$ embedding into $X'$ via some $\iota$ such that 
 $\rho_\lambda(\pi_1(X'))=\rho_\lambda(\iota_*(\pi_1(C)))$ for all $\lambda\in\CP_E'$. If $\rho_\bullet$ is tame, then one can take $\pi=\id_X$.
 \end{Cor}

As an application of \cite{BGP15} and \autoref{Lem-OnTamingCover}, we also record the following result that is a partial strengthening of the main result of  \cite{BGP15}.
\begin{Thm}\label{Thm-AlmostIndependence}
Let $\rho_\bullet$ be a $\BQ$-rational compatible system for $\pi_1(X)$. Then  its restriction $\rho_\bullet|_{\pi_1^\geo(X)}$ is almost independent in the sense of Serre, cf.~\cite[\S~1]{Serre-Independence}. Moreover there is a finite cover $X''\to X$ such that for almost all $\ell\in\CP_\BQ'$ the group $\rho_\bullet(\pi_1^\geo(X''))$ is $\ell$-generated.
\end{Thm}
\begin{proof}
By \autoref{Lem-OnTamingCover} there exists a finite (possibly ramified) Galois cover $X'\to X$ such that $\rho_\ell$ is $\ell$-tame for each $\ell\in\CP_\BQ'$. In particular $\rho_\bullet$ satisfies Condition $\CS(X,\kappa)$ of \cite[Def.~5.1]{BGP15}. Using the notation from \cite[\S~3]{BGP15}, by \cite[Thm.~3.6]{BGP15}, for each $\ell\in\CP_\BQ'$ there exists a short exact sequence
\[1\to M_\ell\to\rho_\ell(\pi_1(X))\stackrel{g_\ell}\to H_\ell\to1\]
of profinite groups with $H_\ell$ finite in $\Jor_\ell(J'(n))$ and $M_\ell$ in $\Sigma_\ell(2^n)$. Now apply \cite[Prop.~5.7]{BGP15} to the family $(g_\ell\circ\rho_\ell\colon \pi_1(X)\to H_\ell)_{\ell\in\CP_\BQ'}$  (compatibility makes no sense here), to find a finite (possibly ramified) cover $X''\to X'$ such that $X''\to X$ is Galois and $\rho_\ell(\pi_1^\geo(X''))\subset M_\ell(2^n)$ for all $\ell\in\CP_\BQ'$. Now apply \cite[Thm.~5.8]{BGP15} to $\rho_\bullet|_{\pi_1^\geo(X'')}$ to deduce that this family, and hence also $\rho_\bullet|_{\pi_1^\geo(X)}$ is almost independent.
\end{proof}
\begin{Rem}\label{Rem-AIndep-AndEtoQ}
For an $E$-rational compatible system $\rho_\bullet$, \autoref{Thm-AlmostIndependence} can be applied to $\Res_{E/\BQ}\rho_\bullet$.
\end{Rem}

\subsection{Geometric monodromy}

Our final result here shows that the geometric monodromy is unaffected by semisimplification. We begin with an important result due to Deligne:
\begin{Thm}[{\cite[Cor.~1.3.9, Thm.~3.4.1(iii)]{Deligne-Weil2}}]\label{Thm-GeomSemisimplicity}
Let $K$ be a finite extension of $\BQ_\ell$ for some $\ell\in\CP_\BQ'$, and let $\rho\colon \pi_1(X)\to \GL_n(K)$ be pure. Then $\rho|_{\pi_1^\geo(X)}$ is semisimple, and the Zariski closure of $\rho|_{\pi_1^\geo(X)}$ is an extension of a finite group by a semisimple group.
\end{Thm}

\begin{Cor}\label{Cor-SSandGeomMonodromy}
Let $K$ and $\rho$ be as in \autoref{Thm-GeomSemisimplicity}, and denote by $\rho^\ssi$ the semisimplification of $\rho$ as a representation of $\pi_1(X)$. Then one has an isomorphism of representations 
\[\rho^\ssi|_{\pi_1^\geo(X)}\cong\rho|_{\pi_1^\geo(X)}\colon \pi_1^\geo(X)\to \GL_n(K).\]
\end{Cor}

\autoref{Cor-SSandGeomMonodromy} is an immediate consequence of \autoref{Thm-GeomSemisimplicity} and the following lemma, whose simple proof we leave as an exercise. 
\begin{Lem}\label{lin-alg} Let $V$ be a finite dimensional representation of a group $G$ and let $H$ be a subgroup of $G$. Denote by $V^{\Gamma\dash\ssi}$ the semisimplification of $V$ as a representation of $\Gamma\in\{G,H\}$. If the action of $H$ on all simple $G$-representations is semisimple, then as $H$-representations one has
\[V^{G\dash\ssi}\cong V^{H\dash\ssi}.\] 
If in addition $V$ is semisimple as an $H$-representation, then $V^{G\dash\ssi}\cong V$ as $H$-representations.  
\end{Lem}

\section{Absolute irreducibility}
\label{Sec-AI}
\autoref{Subsec-BasicsAI} provides various results on absolutely irreducible compatible systems. \autoref{Lem-OnAbsIrredAndTwists} shows that if one member of a compatible system is absolutely irreducible, then all of its members have this property, and also that connected absolutely irreducible compatible systems can be twisted by a character so that their monodromy is semisimple; then $\pi_1(X)$ and $\pi^\geo_1(X)$ have the same monodromy groups. If the monodromy is connected we show that absolute irreducibility is equivalent to the irreducibility of its split motivic representation. \autoref{Prop-GenChinTriple} in \autoref{Subsec-ChinTriple} generalizes \autoref{Thm-MainThmChin} of Chin. It shows the existence of a motivic triple for arbitrary and not only pure compatible systems. In \autoref{Subsec-Gcompat} we explain how to use tools from representation theory to study compatible systems. We recall the notion of a $G$-compatible system for $G$ a reductive group from \cite{BHKT}. Then we show, using \cite{Chin04} and \cite{BHKT}, that any compatible system with split motivic group $M$ can be recovered from a $M$-compatible system (with Zariski dense image) together with its split motivic representation. We show in \autoref{Cor-ChinAndAI2-New} that any connected compatible system with semisimple split motivic group has, up to a (uniform) finite kernel, the same monodromy groups as a suitable connected absolutely irreducible compatible system. This corollary serves as a crucial reduction step in the later~\autoref{Thm-Reduction-IsSaturated}.

In accordance with \autoref{Conv-Semisimple}, from now on all compatible systems are {\em semisimple}.

\subsection{Basic results on absolute irreducibility}
\label{Subsec-BasicsAI}
Recall that a representation $\rho\colon \Gamma\to\GL_n(K)$ for a group $\Gamma$ and a field $K$ is called irreducible if the vector space $K^n$ underlying the representation contains no non-trivial $\Gamma$-stable subspace. It is called absolutely irreducible, if $\rho\otimes_K\overline K$ is irreducible. The latter is equivalent to $\rho\otimes_KL$ being irreducible for all finite field extensions $K\to L$, or for all field extensions $K\to L$.

\begin{Lem}\label{Lem-OnAbsIrredAndTwists}
Let $\rho'\colon \pi_1(X)\to\GL_n(\overline\BQ_\ell)$ be absolutely irreducible. Then:
\begin{enumerate}
\item There exists $\alpha\in\overline\BZ{}_\ell^\times$ such that $\det(\rho'\otimes\rho_\alpha^{-1})$ has finite order for $\rho_\alpha$ defined as in \autoref{eq:RhoAlpha}.
\item If all $\alpha_x:=\det(\rho'(\Frob_x))$ lie in $\overline\BQ$ and are plain of characteristic $p$, then $\alpha$ from (a) has this property. In this case $\rho_\alpha=\rho_{\alpha,\bullet}\otimes_{L_\mu}\overline\BQ_\ell$ where $L=\BQ(\alpha)$ and $\mu\in\CP'_L$ is the place defined from $\BQ(\alpha)\into\overline\BQ\into\overline\BQ_\ell$. If furthermore all $\alpha_x$ are $\#\kappa_x$-Weil numbers of weight $w$, then $\alpha$ is $q$-Weil number of weight~$w$.
\item 
Suppose that all $\alpha_x$ are plain of characteristic $p$. Then through any finite set of points of $X_\reg$ there is a smooth geometrically irreducible curve $C\subset X_\reg$, a number field $E'$ and an $E'$-rational compatible system $\rho'_\bullet$ for $\pi_1(C)$ such that
\begin{enumerate}
\item each $\rho'_\lambda$ is absolutely irreducible, and 
\item for some $\lambda_0\in\CP_{E'}'$ above $\ell$ and embedding $E'_{\lambda_0}\to\overline\BQ_\ell$ one has $\rho'_{\lambda_0}\otimes_{E'_{\lambda_0}}\overline\BQ_\ell\cong \rho'|_{\pi_1(C)}$.
\end{enumerate}
\item Suppose that $\rho_\bullet$ is $E$-rational compatible system of $\pi_1(X)$ such that for some $\lambda_0\in\CP_{E}'$ above $\ell$ and embedding $E_{\lambda_0}\to\overline\BQ_\ell$ one has $\rho_{\lambda_0}\otimes_{E_{\lambda_0}}\overline\BQ_\ell\cong \rho'$. Let $\alpha$ be as in (a). Then:
\begin{enumerate}
\item $\rho_\lambda$ is absolutely irreducible for all $\lambda\in\CP_E$.
\item $\alpha\in \overline\BQ$, and, assuming $\alpha\in E$, $\rho_\bullet\otimes\rho_{\alpha,\bullet}^{-1}$ has finite order determinant and is pure of weight $0$ when restricted to $\pi_1(X_\reg)$, and one has $G_{\rho_\bullet|\pi^\geo_1(X),\lambda}=G_{\rho_\bullet\otimes\rho^{-1}_{\alpha,\bullet}|\pi^\geo_1(X),\lambda}$.
\end{enumerate}
\end{enumerate}
\end{Lem}
\begin{proof}
For (a) we apply \autoref{Lem-AbelianCS}(a) to find $\beta\in\overline\BZ_\ell^\times$ such that $(\det\rho')\otimes \rho_{\beta}^{-1}$ has finite image, where $\rho_\beta$ is the representation $\pi_1(X)\to \GL_1(\overline \BQ_\ell)$ that factors via $\Gamma_\kappa$ and sends $\Frob_\kappa$ to~$\beta$. Let $\alpha$ be an $n$-th root of $\beta$. Then (a) is clear since $\det(\rho'\otimes\rho_\alpha^{-1})=(\det\rho')\otimes\rho_\beta^{-1}$. Part (b) follows now from \autoref{Lem-AbelianCS}(b) and~(c).

To prove (c), let $\rho'':=\rho'\otimes\rho_\alpha^{-1}$ for $\alpha$ from (b). Clearly $\rho''$ is absolutely irreducible and by (a) it's determinant has finite order. Note also that by the \v{C}ebotarov density theorem $\rho''$ and $\rho''|_{\pi_1(X_\reg)}$ have the same image. Let $H:=\kernel\rho''$. From the compactness of $\pi_1(X)$ it follows that $\pi_1(X)/H=\rho''(\pi_1(X))$ contains an open pro-$\ell$ subgroup; see for instance~\cite[page~1]{KLR}. By \autoref{Drinfeld-OnPi3}(a) we can find a smooth geometrically irreducible curve $C$ on $X_\reg$, through any finite number of points on $X_\reg$, such that under the natural map $\pi_1(C)\to\pi_1(X)$ we have $\rho''(\pi_1(X))=\rho''(\pi_1(C))$. From \autoref{Thm-Lafforgue} by L.\ Lafforgue, we see that $\rho''|_{\pi_1(C)}$ is a member of a compatible system $\rho''_\bullet$ (which amounts to (ii) for $\rho''$) and such that all $\rho''_\lambda$ are absolutely irreducible. By taking for $E'$ a sufficiently large number field, we may assume that $\rho''_\bullet$ is $E'$-rational and $L\subset E'$ for $L$ from (b). Part~(c) now follow by defining $\rho'_\bullet:=\rho''_\bullet\otimes \rho_{\alpha,\bullet}$.

Finally, we prove (d). We apply (c) to $\rho_{\lambda_0}\otimes_{E_{\lambda_0}}\overline\BQ_\ell$ to obtain a curve $C$ on $X_\reg$ and an $E'$-compatible system $\rho'_\bullet$ for $\pi_1(C)$ with the properties in (c). Let us replace both $E$ and $E'$ by $EE'$ and call this field $E$ again. Observe that now $\rho_\bullet|_{\pi_1(C)}\cong\rho'_\bullet$ by \autoref{Lem-IsomOfCS}. Hence (d)(i) is implied by (c)(i). The assertion on $\det(\rho_\bullet\otimes\rho^{-1}_{\alpha,\bullet})$ follows from (a) and the classification in \autoref{Prop-AbelianCS}. The purity of $\rho_\bullet\otimes\rho_{\alpha,\bullet}^{-1}$ of weight $0$ on $\pi_1(X_\reg)$ follows by varying the curves $C$ in (c) and using that all $\rho'_\bullet\otimes\rho_{\alpha,\bullet}^{-1}|_{\pi_1(C)}$ are pure of weight zero from \autoref{Thm-Lafforgue}. The last assertion is obvious since $\rho_{\alpha,\bullet}|_{\pi_1^\geo(X)}$ is trivial.
\end{proof}

\begin{Def}
We call an $E$-compatible system $\rho_\bullet$ {\em absolutely irreducible} if for all $\lambda\in \CP_E'$, or equivalently, by \autoref{Lem-OnAbsIrredAndTwists}(c), for one $\lambda\in\CP_E'$, the representation $\rho_\lambda$ is absolutely irreducible.

We call $\rho_\bullet$ {\em absolutely completely reducible} if it is isomorphic to the (finite) direct sum $\bigoplus_i\rho_{i,\bullet}$ of absolutely irreducible compatible systems $\rho_{i,\bullet}$.
\end{Def}

The following result is well-known for single representations.
\begin{Prop}\label{Prop-OnAbsoluteCR}
For every compatible system $\rho_\bullet$ there exists a finite extension $E'$ of $E$ such that $\rho_\bullet\otimes_EE'$ is absolutely completely reducible.
\end{Prop}
\begin{proof}
Fix a place $\lambda\in\CP_E'$ and choose a finite extension $E''/E$ and a place $\lambda''\in\CP_{E''}'$ above $\lambda$ such that $(\rho_\bullet\otimes_EE'')_{\lambda''}$ is a direct sum $\oplus_{i\in I}\rho_{i,\lambda''}$ of absolutely irreducible representations. By \autoref{Lem-OnAbsIrredAndTwists}(a), (b) there exists a finite extension $E'/E''$ such that $\rho_{i,\lambda''}\otimes_{E_{\lambda''}}E_{\lambda'}$ is a member of an $E'$-rational absolutely irreducible compatible system $\rho_{i,\bullet}$ for some $\lambda'\in\CP_{E'}'$ above $\lambda''$ (for all $i$ we can arrange for the same $\lambda'$). Then $\rho_\bullet\otimes_EE'$ and $\oplus_{i\in I}\rho_{i,\bullet}$ are isomorphic compatible systems by \autoref{Lem-IsomOfCS}, and,  by construction, $\oplus_{i\in I}\rho_{i,\bullet}$ is absolutely completely reducible.
\end{proof}

In order to investigate the geometric monodromy, the following two results will be useful:
\begin{Lem}\label{Lem-AIandFiniteDet}
Let $\rho_\bullet$ be absolutely irreducible and connected. Then the following are equivalent:
\begin{enumerate}
\item $\det\rho_\bullet$ has finite order,
\item $\det\rho_\bullet$ is trivial,
\item $G_{\rho_\bullet,\lambda}=G_{\rho_\bullet|_{\pi^\geo_1(X)},\lambda}$ for all $\lambda\in\CP'_E$,
\item $G_{\rho_\bullet,\lambda}$ is semisimple  for all $\lambda\in\CP'_E$,
\item $G_{\rho_\bullet,\lambda}$ has finite center for all $\lambda\in\CP'_E$.
\end{enumerate}
\end{Lem}
\begin{proof}
Since by the absolute irreducibility of the $\rho_\lambda$ and the connectivity of $\rho_\bullet $ the group $G_\lambda$ is connected reductive for all $\lambda\in\CP_E'$, the equivalence (d)$\Longleftrightarrow$(e) is clear. The equivalence (c)$\Longleftrightarrow$(d) follows from \autoref{Lem-GDer=Ggeo}. Also, if $\rho_\bullet$ is connected then so is $\det\rho_\bullet$, and this proves (a)$\Longleftrightarrow$(b). By (d) and Galois cohomology, the maximal abelian quotient of the abstract group $G_\lambda(E_\lambda)$ is finite, and this implies (a). Finally, because $\rho_\bullet $ is absolutely irreducible, the theorem of Burnside shows that the $E_\lambda$-linear hull of $\rho_\lambda(\pi_1(X))$ is $M_{n\times n}(E_\lambda)$, and hence the center of $G_\lambda$ is equal to the center of $\GL_n$ intersected with $G_\lambda$. This gives (b)$\Rightarrow$(e).
\end{proof}

\begin{Prop}\label{Prop-WeightZeroArithMonodromy}
Suppose $\rho_\bullet$ is connected and a direct sum of absolutely irreducible compatible systems $\rho_{i,\bullet}$ each of which has finite order determinant. Then for all $\lambda\in\CP_E$ the group $G_{\rho_\bullet,\lambda}$ is semisimple connected, and one has an equality of monodromy groups $G_{\rho_\bullet,\lambda}=G_{\rho_\bullet|_{\pi^\geo_1(X)},\lambda}$. 
\end{Prop}
\begin{proof}
By our hypotheses the $G_\lambda$ are connected reductive. Since $G_\lambda$ surjects onto each $G_{\rho_{i,\bullet},\lambda}$, the latter are connected as well.\footnote{The converse is not true as can be seen by considering $\SL_2(K)\times\{\pm1\}\to\SL_2(K)\times\SL_2(K),(g,\eps)\mapsto(g,g\cdot \eps)$, where we regard $\{\pm 1\}$ as the center of $\SL_2(K)$ for a field $K$ of characteristic zero.} By \autoref{Lem-AIandFiniteDet} the image of the center $Z_\lambda$ of $G_\lambda$ in $G_{\rho_{i,\bullet},\lambda}$ is finite for all $i$. Since $G_\lambda$ is a subgroup of $\prod_i G_{\rho_{i,\bullet},\lambda}$, we deduce that $Z_\lambda$ is finite, and thus $G_{\rho_\bullet,\lambda}$ is semisimple. The asserted equality of monodromy groups now follows from \autoref{Lem-GDer=Ggeo}.
\end{proof}

\subsection{Chin's theorem in the non-pure case}\label{Subsec-ChinTriple}

In this subsection we shall prove that the purity assumption in \autoref{Thm-MainThmChin} is unnecessary.

Below we write  $\diag(a_1,\cdots, a_t)$ for the diagonal matrix in $\GL_t$ with diagonal entries $a_1,\cdots, a_t$. We shall also use this notation when $a_1$,\ldots,$a_t$ are themselves diagonal matrices with $a_i$ of size $r_i$ to obtain a matrix in $\GL_r$ with $r=\sum_ir_i$. For the following result, recall the notion of a (split) motivic triple from \autoref{Def-ChinTriple}.

\begin{Prop}\label{Prop-GenChinTriple} 
Let $\rho_\bullet$ be an $E$-rational semisimple compatible system. Then there exists a split motivic triple $(F, M, \alpha)$ for $\rho_\bullet$, and it satisfies \autoref{Eqn-ChinGroups2} and~\autoref{Eqn-ChinReps2} of \autoref{Thm-MainThmChin}.
\end{Prop}

\begin{proof}
By \autoref{Lem-FactsOnChin}(a) and (c), in the following we may and will (multiple times) replace $X$ by a connected finite \'etale cover and $E$ by a finite extension.

After replacing $E$ by a finite extension we may assume that $\rho_\bullet$ is absolutely completely reducible (cf. \autoref{Prop-OnAbsoluteCR}) and that there exists a
$1$-dimensional $E$-rational compatible system $\delta_{i, \bullet}$ and an absolutely irreducible $E$-rational compatible systems $\sigma_{i, \bullet}$ pure of weight zero and with finite order determinant such that 
\[\rho_\bullet=\bigoplus_{i=1}^t \sigma_{i, \bullet} \otimes \delta_{i, \bullet}\]
(cf. \autoref{Lem-OnAbsIrredAndTwists}). 
After replacing $X$ by a finite \'etale cover we can assume that all the compatible systems involved so far are connected and that there exists for every $i$ an element $\beta_i\in E^\times$ such that $\delta_{i, \bullet}=\rho_{\beta_i, \bullet}$. In particular the system $\delta_{i, \bullet}$ is trivial on $\pi_1^{\geo}(X)$. We define $\sigma_\bullet:=\bigoplus_{i=1}^t \sigma_{i, \bullet}$, so that $\sigma_\bullet$ is pure of weight $0$ and $\rho_\bullet|_{\pi_1^{\geo}(X)}=\sigma_\bullet|_{\pi_1^{\geo}(X)}$. Denote by $(F, M^{\geo}, \alpha^{\geo})$ a split motivic  triple for $\sigma_\bullet$ that exists by \autoref{Thm-MainThmChin}. Note that for every $\lambda\in P_E'$ the group $G_{\sigma_{\bullet}, \lambda}$ is semisimple and we have $G_{\sigma_{\bullet}, \lambda}=G_{\sigma_{\bullet}|_{\pi_1^{\geo}(X)}, \lambda}$ by \autoref{Prop-WeightZeroArithMonodromy}. If $r_i$ is the rank of $\sigma_{i, \bullet}$ and $r=\sum_{i=1}^t r_i$, then it is also clear that $\alpha^{\geo}(M^{\geo})\subset \GL_{r_1, F}\times\cdots\times \GL_{r_t, F}$. After enlarging $E$ again we can also assume $E=F$. 

Define $\delta_{\bullet}'=\oplus_{i=1}^t \,\rho_{\beta_i, \bullet}^{\oplus r_i}$. Then $\delta'_{\bullet}$ is a compatible system of rank $r$ with values in the center of  $\GL_{r_1, F}\times\cdots\times \GL_{r_t, F}$ and $\rho_\lambda(g)=\sigma_\lambda(g)\delta'_\lambda(g)$ for every $g\in \pi_1(X)$. By possibly replacing $X$ again by a finite cover $X'\to X$ (see~\autoref{Thm-SerreOnComponents}), we can assume that the groups $G_{\delta'_\bullet,\lambda}$, $\lambda\in\CP_E'$, are all connected. Next we construct a split motivic triple for $\delta_\bullet'$: Let 
\[B_i=\diag(\underbrace{\beta_i, \cdots, \beta_i}_{r_i \  \mathrm{times}}),\]
$B=\diag(B_1,\cdots, B_t)\in \GL_r(E)$, and let $\langle B\rangle$ be the cyclic group of the integral powers of $B$. Note that $B$ is a diagonal matrix which lies in the center of  $\GL_{r_1}\times\cdots\times \GL_{r_t}$. Let  $Z$ be the Zariski closure of  $\langle B\rangle$ in $\GL_{r,E}$. Let $Z_\lambda\subset\GL_{r,E_\lambda}$ be the Zariski closure of $\delta'_\lambda(\pi_1(X))$. Because $\langle B\rangle$ is dense in $\delta'_\lambda(\pi_1(X))$ in the profinite topology, using \autoref{Lem-CoeffExtOfCompSys}  it follows that $Z_\lambda\cong Z\otimes_EE_\lambda$. Since the groups $G_{\delta_\bullet',\lambda}=Z_\lambda$ are connected, they are tori, and hence $Z$ is a torus. This shows that 
$(E, Z, \mathrm{incl})$ is a split motivic triple for $\delta_\bullet'$.

We now  show that $G_{\rho_\bullet, \lambda}=G_{\sigma_\bullet, \lambda} G_{\delta_\bullet', \lambda}$. Clearly the left side is contained in the right side.  To show the opposite inclusion, we claim that some open subgroup of $\delta_\lambda'(\pi_1(X))$ lies in $\rho_\lambda(\pi_1(X))$. By \autoref{Cor-GeomAlgFiniteIndex}, there exists a connected finite \'etale cover $X'\to X$, a pure constant field extension that depends on $\lambda$, such that $\sigma_\lambda(\pi_1(X'))=\sigma_\lambda(\pi_1^{\geo}(X'))$. Let $d$ be the degree of $X'\to X$. Choose a section $s: \Gamma_\kappa\to \pi_1(X)$ of the canonical map $p: \pi_1(X)\to \Gamma_\kappa$ and choose $g\in \pi_1(X)$ such that
$p(g)$ is a topological generator of $\Gamma_\kappa$. Then $g^{d}\in \pi_1(X')$ and there exists $h\in \pi_1^{\geo}(X')$ such that $\sigma_\lambda(h)=\sigma_\lambda(g^{d})$. We define $a:=g^{d}h^{-1}$. Then $\sigma_\lambda( a)=1$ and $p( a)$ (topologically) generates an open subgroup of $\Gamma_\kappa$, and thus $\rho_\lambda( a)=\delta'_\lambda( a)$ (topologically) generates an open subgroup of $\delta_\lambda'(\pi_1(X))$. This finishes the proof of the claim.

Because $G_{\delta_\bullet', \lambda}$ is connected, it follows that the Zariski closures of the two groups $\rho_\lambda(\pi_1(X))$ and $\sigma_\lambda(\pi_1(X))\delta_\lambda'(\pi_1(X))$ agree, i.e., that 
\[G_{\rho_\bullet, \lambda}=G_{\sigma_\bullet, \lambda} G_{\delta_\bullet', \lambda},\]
where  $G_{\sigma_\bullet, \lambda}$ is semisimple and $G_{\delta_\bullet', \lambda}$ is a torus commuting with $G_{\sigma_\bullet, \lambda}$. 

The torus $Z$ commutes with the semisimple algebraic group $\alpha^{\geo}(M^{\geo})$ and thus we may define
\[M:=\alpha^{\geo}(M^{\geo})\cdot Z\]
inside $\GL_{r, E}$. It follows that $(E, M, \mathrm{incl})$ is a split motivic triple for $\rho_\bullet$, and that
\autoref{Eqn-ChinGroups2} and~\autoref{Eqn-ChinReps2} of \autoref{Thm-MainThmChin} hold. 
\end{proof}

We have the following immediate consequence (of the above proposition and notion of Chin triple).
\footnote{\cite{LarsenPink92} proves the first part for a density 1 set, Serre proves the second part in general; both for $\BQ$-rational systems. See \cite{LarsenPink92}, the first two pages of the introduction.}
\begin{Cor}\label{Cor-GenChinTriple} 
Let $\rho_\bullet$ be an $E$-rational compatible system. Then
\begin{enumerate}
\item the root datum of the groups $G_{\rho_\bullet,\lambda}$ is independent of $\lambda\in\CP_E'$,
\item the formal character of the representation $G_{\rho_\bullet,\lambda}\subset\GL_{n,E_\lambda}$ is independent of $\lambda\in\CP_E'$.
\end{enumerate}
\end{Cor}

\begin{Cor}\label{Cor2-GenChinTriple}
 The statement of \autoref{Lem-FactsOnChin}(b) remains true without the purity assumption. 
\end{Cor}

\begin{Prop}\label{Prop-ChinAndAI1}
Suppose $\rho_\bullet$ is $E$-rational connected, and let $(F,M,\alpha)$ be a split motivic triple for $\rho_\bullet$. Then $\alpha$ is irreducible if and only if $\rho_\bullet\otimes_EF$ is absolutely irreducible.
\end{Prop}
\begin{proof}
It is clear that if $\rho_\bullet\otimes_EF$ is absolutely irreducible, then $\alpha$ is irreducible. For the converse, suppose that $\alpha$ is irreducible. Let $F'/F$ be a finite extension and $\mu'\in \CP'_{F'}$. Now $(F', M\otimes_F F', \alpha\otimes_F F')$ is a split motivic triple for $\rho_\bullet\otimes_E F'$ by \autoref{Lem-FactsOnChin}(a). As $M$ is connected split reductive we can conclude that $\alpha\otimes_F F'_{\mu'}$ is irreducible  (cf. \autoref{Cor-RepTheory}). From \autoref{Def-ChinTriple} we now see that 
the action of $G_{\rho_\bullet \otimes_E F', \mu'}^\circ$ is irreducible, and hence $(\rho_\bullet\otimes_E F')_{\mu'}$ is irreducible, as desired.
\end{proof}

\subsection{{$G$}-compatible systems}
\label{Subsec-Gcompat}
Let $G$ be a reductive group over $\overline\BQ$. Recall that a representation $\rho\colon \Gamma\to G(\overline \BQ)$ is called $G$-completely reducible if whenever $P\subset G$ is a parabolic subgroup with $\rho(\Gamma)\subset P(\overline \BQ)$, then there exists a Levi subgroup $L\subset P$ such that $\rho(\Gamma)\subset L(\overline \BQ)$. It is called $G$-completely irreducible if there is no parabolic  $P\subset G$ such that $\rho(\Gamma)\subset P(\overline \BQ)$. For $L\supset\BQ$ an algebraically closed field and $g\in G(L)$, we write $g_s\in G(L)$ for the semisimple part of~$g$ in its Jordan decomposition $g=g_sg_u$.

\begin{Def}[{\cite[Def.~6.1]{BHKT}}]\label{Def-GCompSys}
A {\em $G$-compatible system} of $\pi_1(X)$ is a datum \[ ((\rho_\mu)_{\mu\in\CP_{\overline\BQ}'},[g_x]_{x\in|X|})\] consisting of
\begin{itemize}
\item a continuous $G$-completely reducible representation $\rho_\mu \colon \pi_1(X) \to G(\overline\BQ_\mu)$ for every $\mu\in\CP_{\overline\BQ}'$, 
\item a semisimple $G(\overline\BQ)$-conjugacy class $[g_x]$ for every $x\in |X|$
\end{itemize}
such that for any $x\in |X|$ and any $\mu\in\CP_{\overline\BQ}'$, we have \[[g_x]=[\rho_\mu(\Frob_x)_s]\]
as an equality of semisimple conjugacy classes in $G(\overline \BQ_{\mu})$.

We call homomorphisms $\rho,\rho'\colon \pi_1(X)\to G(\overline{\BQ_\ell})$ {\em equivalent}, and write $\rho\cong\rho'$, if they are conjugate. 

We call $G$-compatible systems $\rho_\bullet,\rho'_\bullet$ {\em equivalent}, and write $\rho_\bullet\cong\rho'_\bullet$, if we have $\rho_\mu\cong\rho'_\mu$ for all $\mu\in\CP_{\overline\BQ}'$.

The system $\rho_\bullet$ is called {\em $E$-rational}\ if $G$ is defined over $E$ and if for all $\lambda\in\CP_E'$ there is a~representation $\rho'_\lambda\colon \pi_1(X)\to G(E_\lambda)$, such that for all $\mu\in\CP_{\overline\BQ}'$ one has $ \rho_\mu=\rho'_{\mu|E}\otimes_{E_\lambda}\overline \BQ_{\mu}$.

We say that a representation $\rho\colon\pi_1(X)\to G(\overline{\BQ_\ell})$ is {\em a member} of a $G$-compatible system $\rho_\bullet$, if there exists $\mu\in\CP_{\overline\BQ}'$, and an isomorphism $\iota\colon\overline \BQ_{\mu}\to\overline{\BQ_\ell}$ such that $\rho\cong\iota\circ\rho_\mu$.

We say that $\rho_\bullet$ {\em has Zariski dense image} if $\rho_\mu(\pi_1(X))\subset G(\overline\BQ_\mu)$ is Zariski dense for all $\mu\in\CP_{\overline\BQ}'$.
\end{Def}
\begin{Rem}
Since for $\GL_n$ there is a bijection between semisimple conjugacy classes and monic polynomials of degree $n$ given by the characteristic polynomial (and defined over any field),  
$\GL_n$-compatible systems are the compatible systems we consider in the rest of this work.
\end{Rem}

{\em From now on, for the remainder of this section, we assume that $X$ is a curve.} The following two results from \cite{BHKT} make strong use of results from \cite{Chin04}.
\begin{Thm}[{\cite[Thm.~6.5, Prop.~6.6]{BHKT}}]\label{Thm-GCompSysExist}
Suppose $\rho\colon \pi_1(X)\to G(\overline{\BQ_\ell})$ is a continuous homomorphism with Zariski dense image for some prime $\ell\neq p$ and $G$ is semisimple. Then:
\begin{enumerate}
\item The homomorphism $\rho$ is a member of a $G$-compatible system $\rho_\bullet$.
\item The system $\rho_\bullet$ from (a) has Zariski dense image.
\item The system $\rho_\bullet$ from (a) is unique up to equivalence.
\item Suppose $\rho_\bullet$ is a $G$-compatible system with Zariski dense image. Then there exists a number field $E$ and an $E$-rational $G$-compatible system $\rho'_\bullet$ that is equivalent to $\rho_\bullet$.
\end{enumerate}
\end{Thm}

\begin{Cor}\label{Cor-ChinRepOverM}
Let $\rho_\bullet$ be a $\GL_n$-compatible connected system with semisimple groups $G_\lambda$. Then there exists a finite extension $F$ of $E$ such that $\rho_\bullet\otimes_EF$ has a motivic pair $(M,\alpha)$, and such that there exists an $F$-rational $M$-compatible system $\rho_\bullet^M$ with Zariski dense image, unique up to equivalence, such that $\alpha\circ\rho_\bullet^M\cong \rho_\bullet\otimes_EF$.
\end{Cor}
\begin{proof}
By passing from $E$ to a finite extension, we may assume that $\rho_\bullet$ itself possesses a motivic pair $(M,\alpha)$. Write $M_F$ and $\alpha_F$ for $M\otimes_EF$ and $\alpha\times_EF$ for any field extension $F$ of $E$. Choose $\lambda\in\CP_E$. Then, by the definition of motivic pair, there is a continuous homomorphism $\rho'\colon \pi_1(X)\to M(E_\lambda)$ with Zariski dense image, such that $\rho_\lambda=\alpha_{E_\lambda}\circ \rho'$. By \autoref{Thm-GCompSysExist} there is a finite extension $F$ of $E$ and an $F$-rational $M_F$-compatible system $\rho'_\bullet$, which contains $\rho_\lambda\times_{E_\lambda}F_\mu$ as a member for some $\mu\in\CP_F'$ above $\lambda$. Then $\alpha_F\circ\rho'_\bullet$ is an $F$-rational $\GL_n$-compatible system (with semisimple conjugacy classes $[\alpha(g_x)]_{x\in|X|}$). From $\alpha_F\circ\rho'_\mu\cong\rho_\lambda\times_{E_\lambda}F_\mu$ and \autoref{Lem-IsomOfCS}, the corollary follows.
\end{proof}

Recall that a representation $\beta\colon G\to\GL_n$ of a reductive connected group $G$ is called {\em almost faithful} if $\kernel\beta$ is a finite subgroup scheme. 
\begin{Cor}\label{Cor-ChinAndAI2-New}\label{Cor-TensorDecompForAjdoint}\label{Cor-ExistsAIRep}\label{Cor-RedToAdjoint}
Let $\rho_\bullet$ be a $\GL_n$-compatible connected system with split motivic triple $(F,M,\alpha)$ such that $M$ is semisimple, and suppose that $\rho_\bullet^M$ from \autoref{Cor-ChinRepOverM} is defined over $F$. Then the following hold:
\begin{enumerate}
\item For any representation $\beta$ of $M$ (defined over $F$), we have an $F$-rational compatible connected system $\rho_\bullet^\beta:=\beta\circ\rho_\bullet^M$.
\item If $\alpha=\oplus_i\alpha_i$, then $\rho_\bullet=\oplus_i\rho_\bullet^{\alpha_i}$. In particular, $\rho_\bullet$ is absolutely completely reducible.
\item If $\phi\colon M\to M'$ is a surjective homomorphism, then $\phi\circ\rho^M_\bullet$ is an $M'$-compatible system with Zariski dense image, and if $\alpha'$ is a faithful representation of $M'$, then $(F,M',\alpha')$ is a motivic triple of the compatible system $\rho'_\bullet:=\alpha'\circ\phi\circ\rho^M_\bullet$ and $\rho^{\prime,M'}_\bullet=\phi\circ\rho_\bullet^M$.
\item There exist $\phi$, $\alpha'$ in (c) such that $\phi$ is a central isogeny and $\alpha'$ is faithful and irreducible (and so $\alpha'\circ\phi$ is almost faithful), and then the compatible system $\rho'_\bullet$ is absolutely irreducible.
\item Suppose that $M\cong \prod_i M_i$ for simple groups $M_i$ and denote the projections $M\to M_i$ by $\pr_i$ (such an isomorphism always exists if $M$ is simply connected or of adjoint type). Then $\rho^M_\bullet$ is isomorphic to the product $\prod_i (\pr_i\circ\rho_\bullet^M)$ of the $M_i$-compatible systems $\pr_i\circ\rho_\bullet^M$.
\item Suppose that $M$ is as in (e) and $\alpha$ is irreducible. Then $\alpha=\bigotimes_i\alpha_i$ for faithful irreducible representations $\alpha_i$ of $M_i$ (considered as representations of $M$ via $\pr_i$), the systems $\rho_{i,\bullet}:=\alpha_i\circ\pr_i\circ\rho_\bullet^M$ are absolutely irreducible with motivic triple $(F,M_i,\alpha_i)$, and we have $\rho_\bullet=\bigotimes_i\rho_{i,\bullet}$.
\end{enumerate}
\end{Cor}
It is sometimes useful to depict (d) in the following diagram, where $d=\dim \alpha'$:
\[\xymatrix@C+2pc@R-1pc{
&&\GL_n(F_\mu)\\
\pi_1(X)\ar@/^.8pc/[urr]^{\rho_\mu} \ar@/_.8pc/[drr]_{\rho'_\mu}  \ar[r]^{\rho^M_\mu}&M(F_\mu)\ar[ur]_{\alpha\otimes_FF_\mu}\ar[dr]^{\alpha'\otimes_FF_\mu}&\\
&&\GL_d(F_\mu)\rlap{,}
} \]
\begin{proof} 
Part (a) is immediate from \autoref{Cor-ChinRepOverM}. The first part of (b) is clear from the definitions of the $\rho_\bullet^{\alpha_i}$. For the second part, let the $\alpha_i$ be irreducible summands of $\alpha$. Then the $ \rho_\bullet^{\alpha_i}$ are absolutely irreducible by \autoref{Prop-ChinAndAI1}. Part (c) is straightforward from the previous results. For the construction of $\alpha'$ and $\phi$ in (d), see \autoref{Rem-OnChoicesOfSubstitutes}. Part (e) is obvious. The assertion in (f) is a well-known fact in the representation theory of reductive groups. For some further details see \autoref{Lem-TensorProdOfReps}(c).
\end{proof}
\begin{Rem}\label{Rem-OnChoicesOfSubstitutes}
Either of the following two constructions shows that in Part (d) of \autoref{Cor-ChinAndAI2-New} one may realize $\rho_\bullet'$ as a direct summand of compatible system constructed out of $\rho_\bullet$.
\begin{enumerate}
\item Let $\phi$ be the central isogeny $M\to M^\ad$ to the adjoint quotient of $M$ and write $M^\ad=\prod_{i\in I} M_i$ with $M_i$ simple. Let $\alpha'$ be the representation $\otimes_{i\in I} \Ad_{M_i}$. Because each $\Ad_{M_i}$ is irreducible, so is $\alpha'$. Now $\alpha$ defines a homomorphism $M\to \GL_d$ for $d=\dim\alpha$, and hence a monomorphism $\Ad_M\into \Ad_{\GL_d}\cong \alpha\otimes\alpha^\vee$ of $M$-representations. It follows that $\alpha'$ occurs as a direct summand of $(\alpha\otimes\alpha^\vee)^{\otimes \#I}$.
\item Write $\alpha$ as a sum $\bigoplus_{i\in I} \alpha_i$ with irreducible representations $\alpha_i$. Then $\alpha^{\otimes \#I}$ contains a summand isomorphic to $\gamma:=\bigotimes_{i\in I}\alpha_i$. Let $\beta$ be an irreducible summand of $\gamma$ whose highest weight is the sum of the highest weights of the $\alpha_i$; cf.~\autoref{Thm-RepTheory}. The restriction of $\beta$ to the simply connected cover of $M$ is of the form $\otimes_j\beta_j$ corresponding to the simple factors of that cover, cf.~\autoref{Lem-TensorProdOfReps}(c), and it can easily be verified that the highest weight of each $\beta_j$ is non-trivial. Hence each $\beta_j$, and hence also $\beta$ are almost faithful. Take now $\phi$ as the canonical map $M\to M':=M/\kernel(\beta)$ and $\alpha'$ as the representation on $M'$ induced~from~$\beta$.
\end{enumerate}
\end{Rem}
\begin{Rem}
Suppose the split motivic group $M$ of a connected compatible system is of adjoint type. If $M$ is also simple, then we shall prove in \autoref{Sec-OpenImage}, using \cite{Pink-Compact} by Pink, that there is a natural field of definition $E_\lambda^0\subset E_\lambda$ for each $\rho_\lambda$ over which $\rho_\lambda$ has open image. If $M$ is not simple, this cannot be expected. Take for instance representations $\rho_i$, $i=1,2$, with open image in $\SL_2(\BZ_\ell)$ and $\SL_2(\BZ_{\ell^2})$, respectively, and form their tensor product.
\end{Rem}

\section{Saturation}
\label{Sec-Saturation}
This section begins by recalling the foundational correspondence of Nori between exponentially generated subgroups and nilpotently generated Lie algebras. The second subsection concerns the more general notion of saturation introduced by Serre. The literature seems to focus mainly on algebraically closed base fields while our main interest is in finite base fields. The adaptations are minor, and we also investigate the relation between being exponentially generated and between being saturated. One more aspect of saturation is that it can be developed within general reductive groups and not just within $\GL_n$, which was Nori's original setting. Our treatment focusses on the $\GL_n$ case but in \autoref{Subsec-GenlSat} we also present some results in the general case, useful when combining Weil restriction and saturation. \autoref{Subsec-OnLifting} provides a result on lifting a mod $\ell$-representation of small $\ell$-height of a saturated group to a representation over the Witt vectors of a lift of that group. In \autoref{Subsec-NonReduSat} we show for a new class of groups that they are saturated. For $k=\BF_p$ our findings can be deduced from results in \cite[\S~7]{CHT}. Our motivation for this class of groups stems from some results of Larsen-Pink. Their results attach smooth group schemes over $\BZ_\ell$ to the monodromy at $\ell$ of a $\BQ$-rational compatible systems, if $\ell$ is large. In the last subsection, \autoref{Subsec-LPschemesAndSat}, we first extend their results to $E$-rational systems, and then we show that the smooth reductions of these subschemes of $\GL_n$ are saturated for $\ell_\lambda\gg0$.

We remind the reader once more that from \autoref{Sec-AI} on all compatible systems are {\em semisimple}.

\subsection{Reminders from Nori {\cite[\S~2]{Nori}}}
\label{Subsection-Nori}
We fix a natural number $M\ge3$. By $\ell$ we denote a prime and by $n$ a natural number such that $n\le M < \ell$. We let $k$ be a field of characteristic $\ell$ with an algebraic closure $\overline k\supset k$. As in \cite[Sec.~1]{Nori}, we define
\[\GL_{n}^\unip(k):=\{u\in\GL_n(k)\mid u^\ell=1\}\hbox{ and }M_n^\nilp(k):=\{X\in M_{n}(k)\mid X^\ell=0\} ,\]
and further
\begin{equation}\label{Eqn-ExpLogNDef}
\left.
\begin{array}{l}
\exp_n\!:\!
M_n^\nilp(k)\!\to\! \GL_{n}^\unip(k), X\mapsto \sum_{i=0}^{\ell-1}\frac1{i!}X^i,\ \ \\
\log_n\!:\!
\GL_n^\unip(k)\!\to\! M_{n}^\nilp(k), u\mapsto -\sum_{i=1}^{\ell-1}\frac1{i}(1-u)^i.\end{array}\right.
\end{equation}
Then $\exp_n$ and $\log_n$ are well-defined mutually inverse bijections. For $X\in M_n^\nilp(k)$ let $f_X\colon \BG_a\to \GL_n$ denote the morphism $ t\mapsto \exp_n(tX)$ of algebraic groups; it is defined over~$k$.

\begin{Def}[{\cite[\S~1]{Nori}}]
A Lie subalgebra $L$ of $M_n(k)$ is called {\em nilpotently generated} if $L$ is the $k$-span of $L\cap M_n^\nilp(k)$. 

If $L$ is nilpotently generated, the algebraic subgroup of $\GL_{n,k}$ generated by the one-parameter subgroups $f_X$ for all $X\in L\cap M_n^\nilp(k)$ will be denoted by $\exp_nL$. It is a geometrically connected subgroup scheme.

If $A$ is a closed subgroup-scheme of $\GL_{n,k}$, one defines $\langle\log_n A(k)\rangle$ as the $\BF_\ell$-span of $\{\log_nu\mid u\in A(k)\cap \GL_n^\unip(k)\}\subset M_n^\nilp(k)$. It is an $\BF_\ell$-sub Lie algebra of $M_n(k)$.
\end{Def}

\begin{Def}[{\cite[Def.~2.3, Rem.~2.14]{Nori}}]\label{Def-acceptable}
Let $L$ be a $k$-Lie subalgebra of $M_n(k)$ and let $A$ be a closed subgroup-scheme of $\GL_{n,k}$. Then $(L,A)$ is an {\em acceptable pair} if and only if 
$A$ is smooth and geometrically connected, $\Lie A=L$, and for all $X\in M_n^\nilp(\overline k)$ one has \[X\in L\otimes_k\overline k\Longleftrightarrow \exp_n(X)\in A(\overline k).\]
One calls $L$ {\em acceptable} if there exists $A$ so that $(L,A)$ is an acceptable pair.\footnote{Nori develops this notion more generally over any commutative ring and not only over a field.}
\end{Def}

\begin{Def}[{\cite[Def.~2.4]{Nori}}]\label{Def-ExpoGen}
A closed subgroup-scheme $A$ of $\GL_{n,k}$ is {\em exponentially generated} if there exists $S\subset M_n^\nilp(k)$ such that $A$ is generated by the groups $f_X$ for all $X\in S$.\footnote{There is a related notion by Suslin, Friedlander and Bendel of subgroup schemes of exponential types, cf.~\cite{SFB} or \cite{McNinch-Abelian}.}
\end{Def}

\begin{Thm}[{\cite[Thm.~A]{Nori}}]\label{Thm-Nori-A}
There is a natural number $c_1(M)\ge 2M-1$ so that whenever $k$ is a field of characteristic $\ell>c_1(M)$, and if $1\le n\le M$, then (1) and (2) below are true:
\begin{enumerate}
\item[(1)] If $L$ is a nilpotently generated Lie subalgebra of $M_n(k)$, then $(L,\exp_n L)$ is an acceptable pair, and in particular $L = \Lie (\exp_n L)$.
\item[(2)] If $A\subset \GL_{n,k}$ is exponentially generated, then $A$ is smooth and $\Lie A = \langle\log_n A(k)\rangle$ and furthermore $A = \exp_n (\Lie A)$.
\end{enumerate}
\end{Thm}
Note that if $A$ is exponentially generated, then by (2) the $\BF_\ell$ Lie subalgebra $ \langle\log_n A(k)\rangle$ of $M_n(k)$ is in fact a $k$ Lie subalgebra. By definition it is nilpotently generated, and thus from $A=\exp_n(\Lie A)$ and $\Lie A=\langle\log_n A(k)\rangle$ it follows that $(\Lie A,A)$ is an acceptable pair. As stressed in \cite{Nori}, the above theorem sets up a bijection $L\to \exp L$ between nilpotently generated Lie subalgebras of $M_n(k)$ and exponentially generated subgroup schemes of $\GL_{n,k}$.

\subsection{Basics on $k$-saturation}
Let $1\le n\le M$ and $k$ be as in the previous subsection. We further assume that $\ell>c_1(M)\ge2M-1$. This subsection contains a detailed treatment of some aspects of saturation over perfect fields, building on \cite{Serre-Semisimplicite,Serre-Moursund,Seitz,McNinch-Abelian,Serre-CR,Deligne-Semisimplicite,BDP}, that mostly focus on $k=\overline k$. Throughout this subsection {\em we assume} that $k$ is perfect and we let $G\subset\GL_{n,k}$ be a closed smooth $k$-subgroup such that $G^o/R_u(G)$ is {\em quasi-split}, i.e., $G$ contains closed subgroup $B$ such that the base change $B_{\overline k}$ of $B$ to $\overline k$ is a Borel subgroup of $G_{\overline k}$. The latter condition seems natural since the groups we are interested in should have an ample supply of unipotent subgroups. But also, our main case of interest is that of finite fields, and, by a result of Lang, $G^o/R_u(G)$ is quasi-split if $k$ is finite; see~\cite[Exer.~16.2.9]{Springer}.

By $G_u$ we denote the closed subscheme of $G$ of unipotent elements. By the Jordan decomposition in linear algebraic groups, it is equal to $\GL_{n,u}\cap G$. Since $\GL_{n,u}$ is defined over $k$ -- it is defined by the closed condition that the characteristic polynomial be $(T-1)^n$ --, so is $G_u$. Because $n<\ell$, we have $u^\ell=1$ or equivalently $(u-1)^\ell=0$ for $u\in \GL^\unip_{n}(\overline k)$. For $u\in\GL_{n}^\unip(\overline k)$ and $t\in\overline k$, we define 
\begin{equation}\label{Eqn-DefUtForN}
u^t:=\sum_{i=0}^{\ell-1}  (u-1)^i \binom{t}i =\exp_n(t\log_n u)\in\GL_n(\overline k).
\end{equation}

\begin{Def}[{\cite[\S~4.2]{Serre-Semisimplicite}}]\label{Def-Saturation}
\begin{enumerate}
\item One calls $G$ {\em saturated (in $\GL_{n,k}$)} if for any $u\in G_u(\overline k)$ and $t\in \overline k$, the element $u^t\in\GL_n(\overline k)$ lies in $G(\overline k)$.
\item One defines the {\em saturation of $G$ (over $k$)} as the  intersection of all closed saturated subgroups of $\GL_{n,k}$ that contain~$G$, and denotes it by $G^\sat_k$.
\end{enumerate}
\end{Def}

\begin{Lem}\label{Lem-SaturationIsIntersection}
The intersection of all saturated subgroups of $\GL_{n,\overline k}$ that contain $G$ is defined over $k$ and hence equal to $G_k^\sat\otimes_k\overline k$.
\end{Lem}
\begin{proof}
Let $H\subset \GL_{n,\overline k}$ be a saturated subgroup that contains $G$. Let $A_k$ be the Hopf algebra of $\GL_{n,k}$,  $A_{\overline k}$ its base change to $\overline k$, $I_{H,\overline k}$ the Hopf ideal that defines $H$, and $I_{G,k}$ the Hopf ideal that defines the $k$-subgroup $G$. Note that $\overline k=k^\sep$ since $k$ is perfect. Let $\sigma\in\Gal(\overline k/k)$. Then $\sigma (I_{H,\overline k})$ is an ideal of $A_{\overline k}$, and because the Hopf structure of $A_k$ is defined over $k$, $\sigma I_{H,\overline k}$ is a Hopf ideal.\footnote{A Hopf ideal of $A$ is an ideal $I$ such that $\Delta(I)\subset A\otimes I+I\otimes A$ for the comultiplication $\Delta\colon A\to A\otimes A$, such that $S(I)\subset I$ for the inversion $S\colon A\to A$, and such that $\eps(I)=0$ for the counit $\eps\colon A\to k$.} The ideal $\sigma (I_{H,\overline k})$ defines the subgroup $\sigma H\subset \GL_{n,\overline k}$. The group $\sigma H$ is saturated, because $G_u$ is defined over $k$ and for $u\in G(\overline k)$ and $t\in\overline k$ one has $\sigma(u)^t=\sigma(u^{\sigma^{-1}(t)})\in \sigma H$. The intersection of the $\Gal(\overline k/k)$-orbit of $H$ is defined by the ideal of $A_{\overline k}$ spanned by $J_H:=\bigcup_{\sigma\in\Gal(\overline k/k)}\sigma I_{H,\overline k}$. Since $G$ is contained in $H$ we have $I_{G,k}\otimes_k\overline k\subset J_H$. Moreover $J_H$ is a Hopf ideal, and by its very definition $J_H$ is invariant under $\Gal(\overline k/k)$. It follows that $\bigcap_{\sigma\in\Gal(\overline k/k)} \sigma(H)$ is saturated, contained in $H$, defined over $k$ and contains $G$. The lemma follows; cf.~\cite[I.14]{Borel}.
\end{proof}

\begin{Cor}\label{Cor-SatBaseChange}
For any algebraic extension $k\to k'$, the base change $(G^\sat_k)\otimes_kk'\subset\GL_{n,k}\otimes_kk'\cong\GL_{n,k'}$ is the saturation of $G\otimes_kk'$ in $\GL_{n,k'}$.
\end{Cor}
\begin{proof}
By the previous lemma, the group $G^\sat_k$ is the intersection of all saturated subgroups $H$ of $\GL_{n,\overline k}$ that contain $G$ -- because this intersection is defined over $k$. But the same intersection also defines $(G\otimes_kk')^\sat_{k'}$. Hence the result is clear.
\end{proof}
\begin{Cor}
Suppose $H$ is a closed subgroup of $\GL_{n,k}$ such that $H\otimes_kk'$ is saturated for some finite extension $k\to k'$. Then $H$ is saturated.
\end{Cor}
\begin{proof}
Consider $H\subset H^\sat_k$. By base change this becomes $H\otimes_kk'\cong (H\otimes_kk')^\sat_{k'}\cong H^\sat_k\otimes_kk'$. But then we must have $H\cong H^\sat_k$.
\end{proof}

The following result and its proof are inspired by \cite{Serre-Semisimplicite}, and in particular \S~4 of op.cit.
\begin{Prop}\label{Prop-SatAndOtimes}
Define $\otimes'\colon \GL_{n_1,k}\times\GL_{n_2,k}\to\GL_{n,k},(A,B)\mapsto A\otimes B$ for $n_i\in\BN$ such that $n=n_1n_2$. Let $\Gamma$ be a finite group with homomorphisms $r_i\colon\Gamma\to \GL_{n_i,k}$, $i=1,2$, and set $r=\otimes'(r_1\times r_2)\colon \Gamma\to\GL_{n}$. Let $G_i\subset\GL_{n_i,k}$ be closed saturated subgroups for $i=1,2$. Then
\begin{enumerate}
\item $\otimes'(G_1\times G_2)$ is saturated in $\GL_{n,k}$.
\item  If $G_i$ contains $r_i(\Gamma)$, $i=1,2$, and $r(\Gamma)_k^\sat=\otimes'(G_1\times G_2)$, then $r_i(\Gamma)^\sat_k=G_i$ for $i=1,2$.
\end{enumerate}
\end{Prop}
\begin{proof}
Note first that $\otimes'(G_1\times\{1\})$ and $\otimes'(\{1\}\times G_2)$ are commuting subgroups in $\GL_n$. Hence $\otimes'(G_1\times G_2)$ is a subgroup of $\GL_{n,k}$. It is easy to see that it is closed. Part (a) now follows from the following two simple facts whose proof we leave to the reader: (i) If $u=u_1\otimes u_2\in \otimes'(G_1\times G_2)\subset\GL_n$ is unipotent in $\GL_n$, then so are $u_i\in \GL_{n_i}$; (ii) one has $\exp_n(A\otimes 1_{n_2})=\exp_{n_1}(A)\otimes 1_{n_2}$ and $\exp_n(1_{n_1}\otimes B)=1_{n_1}\otimes \exp_{n_2}(B)$ for $A\in M_{n_1}^\nilp(\overline k)$ and $B\in M_{n_2}^\nilp(\overline k)$.

To see (b), note that by (a) and our hypotheses we have inclusions 
\[\otimes'(G_1\times G_2)=r(\Gamma)_k^\sat\subset \otimes'(r_1(\Gamma)_k^\sat\times r_2(\Gamma)_k^\sat)\subset \otimes'(G_1\times G_2).\]
It follows that we have equality everywhere, and (b) follows from the injectivity of $\otimes'$.
\end{proof}

In the following we write $\wt G$ for $R_u(G) G^\der\unlhd G$ and $G^+$ for the subgroup of $G$ generated by $G_u$; recall that by our convention $G^\der=(G^o)^\der$. As observed above, $G_u$ is defined over $k$, and, because $k$ is perfect, $R_u(G)$ is defined over $k$ (see.~\cite[12.1.7]{Springer}). Thus $\wt G$ and $G^+$ are defined over $k$. We now clarify some further properties of $\wt G$ and~$G^+$.
\begin{Lem}\label{Lem-OnQuasiSplit}
\begin{enumerate}
\item The group $\wt G$ contains a Borel subgroup $\wt B$ defined over $k$.
\item Let $T$ be a maximal torus of a $\wt B$ as in (a). Then $\wt G$ is generated by the~subgroup schemes $R_u(B)$, where $B$ ranges over the Borel subgroups of $\wt G$ defined over $k$ and containing~$T$.\footnote{The result only needs $\ell>\max\{3, m+1\mid G^o/R_u(G)\hbox{ has a quotient of type }A_m\}$.}
\end{enumerate}
\end{Lem}
\begin{proof}
Since $R_u(G)$ is closed normal connected unipotent in $\wt G$ and defined over the perfect field $k$, and since the inverse image of a $k$-Borel subgroup under the canonical map $\wt G\to \wt G/R_u(G)$ is a $k$-Borel subgroup, it suffices to prove the lemma for the connected semisimple group $H=\wt G/R_u(G)$ which satisfies $\wt H=H$. By hypothesis $H$ is quasi-split, and hence (a) is clear. 

Let $\wt B$ be a $k$-Borel subgroup of $H$ and $T\subset \wt B$ a maximal $k$-torus. Let $H'$ be the closed (reduced) smooth subgroup scheme of $H$ generated by the unipotent subgroups $R_u(B)$ where $B$ runs through all $k$-Borel subgroups of $H$ containing $T$. By \cite[5.2.12]{ConradSGA3}, there exists a unique $k$-Borel subgroup $B'$ of $H$ such that $\wt B\cap B'=T$, and over $\overline k$ the group $B'$ is the opposite Borel subgroup to $\wt B$. It follows that $\Lie H'\otimes_k\overline k$ contains all subspaces of $\Lie H\otimes_k\overline k$ on which $T$ acts non-trivially. In particular, it implies that $H'(\overline k)$ contains the root subgroups of $H(\overline k)$. Since $H$ is semisimple, we deduce $H'=H$ which proves~(b).
\end{proof}
\begin{Cor}\label{Cor-OnG+}
\begin{enumerate}
\item The group  $\wt G$ is normal in $G$. It is contained in $G^+$.
\item The groups $G^+/\wt G$ and $G^o/\wt G$ commute in $G/\wt G$. 
\item The homomorphism $G^+/\wt G\to G/G^o$ is well-defined, its image is the subgroup generated by all element of $G/G^o$ of order~$\ell$, and it sets up a bijection $(G^+/\wt G)_u\to (G/G^o)_u$. 
\end{enumerate}
\end{Cor}
\begin{proof}
To see (a) note that $R_u(G)$ and $G^\der$ are characteristic subgroups of $G$ and hence so is $\wt G=R_u(G)G^\der$. The assertion $\wt G\subseteq G^+$ follows from \autoref{Lem-OnQuasiSplit}. To prove (b), observe first that $G/\wt G$ is an extension of the finite group $G/G^o$ by the torus $\overline T:=G^o/\wt G$. It will suffice to show that any element $\bar u$ of $G^+/\wt G$ acts trivially on $\overline T$. If not, passing to $\overline k$ we get a non-trivial homomorphisms $\BZ/(\ell)\to \Aut(\overline T_{\overline k})\cong\GL_m(\BZ)$ for $m=\rank \overline T\le n$. The image of a generator of $\BZ/(\ell)$ is a matrix $A\in\GL_m(\BZ)$ of exact order $\ell$. The cyclotomic polynomial $\Phi_\ell=(x^\ell-1)/(x-1)$ is then irreducible and by the Cayley-Hamilton Theorem it must divide the characteristic polynomial of $A$. This implies $\ell-1\le m\le n\le M$, contradicting our hypothesis $\ell\ge2M-1(\ge M+2)$ since $M\ge 3$. This completes (b) and shows the first part of~(c).

For the second part of (c) recall first that under homomorphisms of algebraic groups unipotent elements are mapped again to such (see \cite[Thm.~2.4.8]{Springer}). To complete the second assertion of (c), let $\bar u\in G/G^o$ be unipotent, i.e., of order $\ell^m$ for some $m\ge1$. Let $u\in G$ be a preimage, and split $u=u_su_u$ in its (commuting) unipotent and semisimple parts. We have $u_u^\ell=1$ because $\ell>n$ and the image of $u_s$ is semisimple and hence of order prime to $\ell$. It follows that $\bar u$ is the image of $u_u$, and that its order is $\ell$ (cf.~\cite[Lect.~3, Prop.~3]{Serre-Moursund}).

It now remains to show that $(G^+/\wt G)_u\to (G/G^o)_u$ is injective. So suppose that $u,u'$ are elements in $(G^+/\wt G)_u$ that have the same image in $(G/G^o)_u$. It follows that there exists $t$ in the torus $G^o/\wt G$ such that $u'=ut$. Clearly the elements of the torus are semisimple, and we deduce $u=u'$ from the uniqueness in the Jordan decomposition of linear algebraic groups. 
\end{proof}
\begin{Ques}
We do not know if the kernel of $G^+/\wt G\to G/G^o$ is always finite, or if $(G^+\cap G^o)/\wt G$ can contain a torus.
\end{Ques}
\begin{Lem}\label{Lem-BasicsOnSat} 
\begin{enumerate}
\item If $T$ is a torus, $\phi\colon G\to T$ is a homomorphism of algebraic groups, then $G$ is saturated if and only if $\kernel\phi$ is so. 
\item $G$ is saturated if and only if $G^o$ is saturated and $G/G^o$ has order prime to $\ell$; cf.~\cite[Prop.~11]{Serre-Semisimplicite}.
\item $G$ is saturated, if and only if $\wt G$ is saturated and $\wt G=G^+$.
\item If $G\subset \GL_{n,k}$ is exponentially generated, then $G$ is saturated and $G=G^+=\wt G$.
\end{enumerate}
\end{Lem}
\begin{proof}
Concerning (a) and (b) let us first show that if $G$ is saturated, then so are $\kernel \phi$ and $G^o$. This follows from the fact that any homomorphism from $\BG_a$ to a finite group scheme or a multiplicative group scheme is trivial. In (b) note also that the subgroup of $G/G^o$ generated by elements of order $\ell$ coincides with the image of $G^+$ by \autoref{Cor-OnG+}(c). However because $G$ is saturated, every element of $G^+$ is contained in a $1$-parameter subgroup and hence $G^+\subset G^o$ which implies that $G/G^o$ is of order prime to $\ell$.

For the `if part' of (a), suppose that $\kernel \phi$ is saturated. Now note that the inverse image of $T\setminus\{1_T\}$ under $\phi$ is disjoint from $G_u$, and hence $G$ is saturated if and only if $\kernel\phi$ is saturated. For the `if part' of (b) note that if $G/G^o$ is of order prime to $\ell$, then by  \autoref{Cor-OnG+}(c), we have $G^+\subset G^o$, and hence $G_u\subset G^o$, so that $G$ is saturated if and only if $G^o$ is so.

We now prove (c) and begin with the only if part. If $G$ is saturated, then by (b) $G^+\subseteq G^o$ and $G^o$ is saturated. Let $\psi\colon G^o\to T$ be the maximal torus quotient of $G^o$. It is easy to see that $\wt G$ is the kernel of $\psi$. Hence $\wt G$ is saturated by (a), and $G^+\subset \wt G$ because $G^+$ maps to the identity under $\psi$. Since we also have $\wt G\subseteq G^+$ by \autoref{Cor-OnG+}, we obtain $\wt G=G^+$.

For the if part of (c) assume that $\wt G=G^+$ and that $\wt G$ is saturated. The first condition together with part(b) and \autoref{Cor-OnG+}(c) shows that it suffices to prove that $G^o$ is saturated. For this we apply the `if part' of (a), which allows us to deduce the saturatedness of $G^o$ from that of $\wt G=\kernel \psi$ with $\psi\colon G^o\to T$ as in the previous paragraph.

The first assertion of (d), i.e., that $G$ is saturated, follows from $\ell>c_1(M)$, \autoref{Thm-Nori-A} and Definitions~\ref{Def-ExpoGen} and~\ref{Def-acceptable}: Let $u$ be in $G_u(\overline k)$. Then $\log_n u\in M_n^\nilp(\overline k)$, $\exp_n(\log_n u)=u$, and by \autoref{Def-acceptable}(c), we have $\log_nu\in \Lie G\otimes_k\overline k$. But then $t\log_nu\in  \Lie G\otimes_k\overline k$ for all $t\in\overline k$, and again by \autoref{Def-acceptable}(c) we have $u^t=\exp_n(t\log_nu)\in G(\overline k)$. The equality $G^+=\wt G$ now follows from~(c), and $G=G^+$ holds because $G$ is exponentially generated. 
\end{proof}

\begin{Prop}\label{Prop-SatOfWtG}
$G$ is saturated if and only if $\wt G=G^+$ and $\wt G$ is exponentially generated. 
\end{Prop}

\begin{proof}
By \autoref{Lem-BasicsOnSat}(c) and (d), we may assume that $\wt G$ is saturated, and we need to show that it is exponentially generated. By \autoref{Lem-OnQuasiSplit}, $\wt G$ is generated by its subgroup schemes $R_u(B)$, where $B$ ranges over the Borel subgroups of $\wt G$ defined over $k$. Let now $U=R_u(B)$ for such a $B$. It suffices to show that any such $U$ is exponentially generated.

Because $k$ is perfect, the connected unipotent group $U$ is $k$-split, and hence trigonalizable inside $\GL_{n,k}$, see~\cite[Ex.~12.3.5(3), Rem.~14.1.5, Thm.~14.1.4]{Springer}. This means that there exists $g\in \GL_n(k)$ such that $gUg^{-1}$ is a closed subgroup of the subscheme $U_{n,k}\subset \GL_{n,k}$ of upper triangular unipotent matrices. Thus below we shall assume $U=U_{n,k}\cap\wt G$.

For $1\le c \le n$ define $U_{n,k}^{(c)}$ as the normal closed subscheme of $U_{n,k}$ of elements $u$ such that the matrix $(a_{i,j})_{i,j\in\{1,\ldots,n\}}:=u-1$ satisfies $a_{i,j}=0$ for all $i\ge j-c$. Clearly the quotient
$U_{n,k}^{(c)}/U_{n,k}^{(c+1)}$ is isomorphic to $\BG_a^{c}$. From \cite[Thm.~14.2.6]{Springer} one deduces an isomorphisms
\[ U\cong \prod_{c=1}^{n} (U\cap U_{n,k}^{(c)})/ (U\cap U_{n,k}^{(c+1)})\cong\prod_{c=1}^{n}\BG_a^{d_{U,c}} \]
of $k$-schemes with $d_{U,c}$ defined by any isomorphism $(U\cap U_{n,k}^{(c)})/ (U\cap U_{n,k}^{(c+1)})\cong\BG_a^{d_{U,c}}$. It follows that the $k$-span of $\langle\log_n U(k)\rangle\subset \Lie U_n(k)$ has dimension at least $\dim U$. If we take a basis $(X_i)_{i\in I}$ of this $k$-span it follows that the subgroup $U'$ of $U_{n,k}$ generated by the $f_{X_i}$, $i\in I$, has dimension at least $\#I\ge\dim U$. By the saturation property of $\wt G$ these subgroups lie in $\wt G$, and hence for dimension reasons $U'=U$. This proves that $U$ is exponentially generated. 
\end{proof}

\begin{Rem}\label{Rem-NoriEnvelope}
Suppose $k$ is finite and $G\subset\GL_{n,k}$ satisfies $G=G^+$. This holds for instance if $G$ is a finite group generated by unipotent elements. By \autoref{Prop-SatOfWtG} the group $G^\sat_k$ coincides with what is often called the Nori envelope of $G\subset\GL_{n,k}$, i.e.\ the smallest exponentially closed subgroup that contains~$G$ (which only exists if $G$ is generated by $G_u$, i.e., if $G=G^+$).
\end{Rem}
\begin{Cor}\label{Cor-SatAndWtG}
Suppose $H$ is a normal subgroup of $G$. Then
\begin{enumerate}
\item The group $G$ normalizes $H^\sat_k$, and $H^\sat_k$ is a normal subgroup of $H^\sat_kG$.
\item The group $G^+$ is a normal subgroup of $G$, and $G^\sat=(G^+)^\sat_k G$.
\end{enumerate}
\end{Cor}
\begin{proof}
For (a) note that by the normality of $H$ in $G$, the adjoint action of $G$ on $\GL_n$ permutes the saturated subschemes of $\GL_{n,k}$ that contain $H$, and so by \autoref{Lem-SaturationIsIntersection}, the group $G$ normalizes $H^\sat$. The second part of (a) is immediate from the first.

For the first part of (b) note that the adjoint action of $G$ on $\GL_{n,k}$ clearly preserves $G_u$ and thus $G^+$, and hence $G^+$ is normal in $G$. Regarding the second part of (b), observe that $\supseteq$ is clear from the definitions. To show $\subseteq$ it thus suffices to show that $(G^+)^\sat_k G$ is saturated. Now from the first paragraph it follows that $(G^+)^\sat_k$ is a normal subgroup of $(G^+)^\sat_k G$. Hence we obtain canonical epimorphisms
\[ (G^+)^\sat_k G\stackrel\phi \longrightarrow \bar G:=(G^+)^\sat_k G/(G^+)^\sat_k\cong G/((G^+)^\sat_k \cap G)\stackrel\psi \longleftarrow G/G^+.\]
From $\psi$, the definition of $G^+$ and \autoref{Cor-OnG+}(c), it follows that $\bar G^o$ is a torus and $\bar G/\bar G^o$ is of order prime to $\ell$. From \autoref{Lem-BasicsOnSat}(a) and (b) applied to $\phi$, and from the saturatedness of $(G^+)^\sat_k$ we deduce that $(G^+)^\sat_k G$ is saturated.
\end{proof}
\begin{Cor}\label{Cor-SatIsConn}
The group $(G^+)^\sat_k$ is connected. In particular, if $G$ is connected, then so is $G^\sat_k$, and if $G$ has no torus quotient, then neither does $G^\sat_k$.
\end{Cor}
\begin{proof}
Suppose first that $G=G^+$; i.e., $G^o/R_u(G)$ is connected semisimple and $G/G^o$ is a finite and generated by the elements of $(G/G^o)_u$ which have order dividing $\ell$ (see the proof of \autoref{Cor-OnG+}). It will suffice to show that $G^\sat_k$ is exponentially generated. 

In a first step, we choose for any unipotent element $\bar u\in G/G^o(\overline k)$ a unipotent element $u\in G(\overline k)$ mapping to it; see the proof of \autoref{Cor-OnG+}(c). Then $G^\sat$ has to contain the $1$-parameter subgroups $t\mapsto u^t$ for any such $t$. The group $H$ generated by these and $G$ is invariant under $\Gal(\overline k/k)$ and hence a $k$-group. Moreover $H=H^+=\wt H$, and by construction $H$ must be contained in $G^\sat_k$. Hence it suffices to show that $G^\sat_k$ is exponentially generated, assuming~that~$G^+=\wt G$. 

For this we proceed as in the proof of \autoref{Prop-SatOfWtG}. The group $G=\wt G$ is generated by the groups $U_B=R_u(B)$ where $B$ ranges over the Borel subgroups of $G$ defined over $k$ and containing a fixed maximal torus $T$, and for a fixed $U_B$ let $g_B\in\GL_n(k)$ be such that $g_B U_B g_B^{-1}=g_B Gg_B^{-1}\cap U_{n,k}$. Let $U'_B$ be the subgroup of $g_B^{-1}U_{n,k}g_B$ generated by $U_B$ together with the images of the $f_X$ for $X\in \langle\log_n U_B(k)\rangle$. If $U'_B=U_B$ then as in the proof of \autoref{Prop-SatOfWtG} we deduce that $U_B$ and hence $ G$ are exponentially generated. 

Suppose now that some $U'_B$ properly contains $U_B$, so that $\dim U'_B>\dim U_B$. Let $G'$ be the subgroup of $\GL_{n,k}$ generated by the $U'_B$. Because the $U'_B$ are connected unipotent, the group $G'$ is connected and satisfies $\wt G'=G'$. Moreover we have $\dim G'>\dim G$. We conclude the proof by induction on $\dim G$.

For the remaining assertions, note that by \autoref{Cor-SatAndWtG} we have $G^\sat_k=(G^+)^\sat_kG$. By what we just proved $(G^+)^\sat_k$ is connected and has no non-trivial torus quotient; cf.~\autoref{Lem-BasicsOnSat}(d). Thus if $G$ is connected or has no-non-trivial torus quotient, the same property will hold for~$G^\sat_k$.
\end{proof}
The following result is straightforward from \cite[4.3]{Serre-Semisimplicite}, \cite[Thm.~8]{Serre-Moursund} and  \autoref{Cor-SatIsConn}. We omit the proof.
\begin{Prop}[Serre]\label{Prop-GConnSemisimple}
Suppose the representation of $G$ on $V=k^n$ afforded by $G\subseteq \GL_{n,k}$ is a direct sum $V=\oplus_i V_i$ of absolutely irreducible representations. Then 
\begin{enumerate}
\item $G^\sat_k$ preserves each $V_i$ and acts absolutely irreducibly on it,
\item $(G^\sat_k)^o$ is reductive,
\item if $G=G^+$, then $G^\sat_k$ is semisimple connected.
\end{enumerate}
Conversely, if $(G^\sat_k)^o$ is reductive, then the action of $G$ on $V$ is semisimple.\footnote{The condition $n=h_{\GL_n}\le p$ needed in \cite[Thm.~8]{Serre-Moursund} follows from our standing assumption $\ell>2M-1$; see~\autoref{Rem-Coxeter} for the meaning and value of $h_{\GL_n}$.}
\end{Prop}
\begin{Prop}\label{Prop-GConnSemisimple2}
Suppose that  the representation of $G$ on $V=k^n$ afforded by $G\subseteq \GL_{n,k}$ is a direct sum $V=\oplus_i V_i$ of representations $V_i$ of $G$, and denote by $\pr_i\colon\prod_i \Aut_k(V_i)\to\Aut_k(V_i)$ the projection onto the $i$-th factor. Then $\pr_i$ maps $G^\sat_k$ onto the saturation $\pr_i(G)^\sat_k$ formed inside $\Aut_k(V_i)$.
\end{Prop}

\begin{proof}
Observe first that the group $\prod_i\Aut_k(V_i)$ is saturated. Also, since $G$ lies in $\prod_i \Aut_k(V_i)$, for elements of $G$ we have $\log_n=\bigtimes_i \log_{\dim V_i}$ and for elements of $\Lie G$ we have $\exp_n=\bigtimes_i \exp_{\dim V_i}$. Thus (i): if $H\subseteq\Aut_k(V_i)$ is saturated and contains $\pr_i(G)$, then $\pr_i^{-1}(H)$ is saturated and contains $G$. Moreover (ii): if $G\subseteq N\subseteq \prod_i\Aut_k(V_i)$ is saturated, then $\pr_i(N)$ is saturated and contains $\pr_i(G)$. Now from (i) we deduce the inclusion $G_k^\sat\subseteq\pr_i^{-1}(\pr_i(G)_k^\sat)$, and from (ii) that $\pr_i(G)_k^\sat\subseteq\pr_i(G^\sat_k)$, and the assertion follows.
\end{proof}

The following result will be useful when comparing geometric and arithmetic mod $\lambda$ monodromy.
\begin{Prop}\label{Prop-GeomArithSaturation}
Suppose that : (i) $G\subset\SL_n$,  (ii) $N$ is a closed normal subgroup of $G$ with finite cyclic quotient, (iii) $N(k)$ acts absolutely irreducibly on $V:=k^n$, (iv) $\ell>n$. Then
\begin{enumerate}
\item $(N^\sat_k)^o=(G^\sat_k)^o$ and $G^\sat_k=G N^\sat_k$, and if $N=N^+$, then also $(N^\sat_k)^o=N^\sat_k=(G^+)^\sat_k$;
\item the group $G^\sat_k/N^\sat_k$ is finite cyclic, and its order divides $n!$ and is bounded by $n(n-1)$.
\end{enumerate} 
\end{Prop}
\begin{proof}
Let us first prove the proposition under the further hypothesis $N=N^+$. By \autoref{Cor-SatAndWtG}, the group $N^\sat_k$ is normal in $N^\sat_kG$, and from (ii) it follows that $\bar G:=N^\sat_kG/N^\sat_k$ is finite cyclic. Consider the adjoint homomorphism $c\colon N^\sat_kG\to\Aut(N^\sat_k),g\mapsto c_g$ defined by conjugation, and denote by $c^\out\colon N^\sat_kG\to\Out(N^\sat_k)$ the induced homomorphism to the outer automorphisms. By~\cite[7.1]{ConradSGA3}, the group $\Out(N^\sat_k)$ injects into the group of automorphisms of a based root datum attached to $N^\sat$. Because of $N=N^+$ and (iii),  \autoref{Prop-GConnSemisimple} implies that $N^\sat_k$ is connected semisimple, and it follows that $\Out(N^\sat_k)$ is a subgroup of the automorphisms of the Dynkin diagram of $N^\sat_k$. The Dynkin diagram (of $N^\sat_k$ over $\overline k$) is the union of the diagrams of the simple factors of $N^\sat_k\otimes_k{\overline k}$, and the diagram contains of at most $n-1=\rank \SL_n\ge\rank N^\sat_k$ vertices. It follows that the maximal cyclic subgroup in $\Out(N^\sat_k)$ has order bounded by $n-1$. Define $H$ as the kernel of the restriction $c^\out|_G$. Then $N\subset H\subset G$, the group $G/H$ is cyclic of order bounded by $n-1$, and $\kernel c^\out=N^\sat_k H$.

Let now $h\in H$ be an element whose image is a generator of the cyclic group $H/N$, and let $h'\in N^\sat_k(\overline k)$ be an element such that $c_{h'}=c_h$. Then $h'h^{-1}$ maps to a generator of $N^\sat_kH/N^\sat_k$, and $c_{h'h^{-1}}$ is trivial, so that $h'h^{-1}$ commutes with $N_k^\sat$ and hence by (iii) lies in the scalars in $\SL_n$. It follows that the order of $h'h^{-1}$ divides $n$, and hence that $N^\sat_kH/N^\sat_k$ is cyclic of order dividing~$n$. We deduce that the order of $GN^\sat_k/N^\sat_k$ divides~$n!$ and is bounded by~$n(n-1)$.

From (iv) it now follows that $GN^\sat_k/N^\sat_k$ is finite of order prime to $\ell$. Then \autoref{Lem-BasicsOnSat}(b) implies that $GN^\sat_k$ is saturated, and hence equal to $G^\sat_k$. Thus (b) is proved. The equality $N^\sat_k=(G^\sat_k)^o$ in (a) follows from (b) because $N^\sat_k$ is connected and saturated. The first equality in (a) follows by applying the proposition to $N\subseteq G^+$, in which case one deduces the equality from~(b).

Let now $N$ be arbitrary. Choose $g\in G$ such that $g$ maps to a generator of $G/N$, and let $G'$ be the subgroup of $G$ generated by $N^+$ and $g$. Observe that $N^+$ is characteristic in $G$, and in particular a normal subgroup. Hence we can use (a) and (b) for the pair $(N^+,G')$. From $GN^\sat_k =G'N^\sat_k=G'(N^+)N^\sat_k$ we deduce that $GN^\sat_k/N^\sat_k\cong G'/(N^\sat_k\cap (G'N^+)^\sat_k)$ is a quotient of $G'/(N^+)^\sat_k$, and this proves (b). From (b) it follows that $GN^\sat_k$ is saturated since the quotient $GN^\sat_k/N^\sat_k$ is of order prime to $\ell$; cf. the argument proving~\autoref{Lem-BasicsOnSat}(b). From (b) we also deduce $(G^\sat_k)^o\subset N^\sat_k\subset G^\sat_k$, and this completes the proof of~(a).
\end{proof}

\subsection{On saturation within any reductive group~$G$}
\label{Subsec-GenlSat}
In this subsection, $G$ will denote a connected reductive group over $k$. We shall recall some results from \cite{Serre-Moursund} on saturation within arbitrary reductive groups $G$. Our presentation is also inspired by~\cite{BDP}. 

We begin with the definition of Dynkin height and Coxeter number.
\begin{Def}[{\cite[II.2]{Serre-Moursund} or \cite[Def.~4.3]{BDP}}] \label{Def-Dynkin} 
Let $G$ be a reductive group over $k$, $T$ a maximal torus of $G_{\overline k}$, and $R^+$ the set of positive roots for some Borel subgroup containing $T$. The {\em Dynkin height} of a representation $V$ of $G$ is 
\[\height_G(V):=\max\{\sum_{\alpha\in R^+}\langle\lambda,\alpha^\vee\rangle\mid \lambda \hbox{ a weight for the action of $T$ on $V\otimes_k\overline k$}\}.\]

A representation $\rho\colon G\to \GL(V)$ with $\ell>\height_G(V)$ is said to be of {\em low $\ell$-height for $G$}.
\end{Def}
\begin{Rem}\label{Rem-Dynkin}
Suppose $G$ has a root base $(\alpha_i)_{i\in I}$, and $V$ is an irreducible representation with dominant weight $\lambda^+$ and smallest weight $\lambda^-$, and write $\lambda^+-\lambda^-=\sum_im_i\alpha_i$ for $(m_i)\in\BQ_{\ge0}^I$. Then $\height_G( V)=\sum_i m_i$. For instance, if $G=\GL_n$ and $V$ is the standard representation, then $\height_G(\bigwedge^i V)=i(n-i)$ (cf.~\cite[Part II, Lect.~2]{Serre-Moursund}).
\end{Rem}
\begin{Ex}
Let $E$ be a number field with the ring of integers $\CO_E$ and $n\in\BN$. Let $\CG$ be a reductive group scheme over $\CO_E[1/n]$ and $\rho\colon \CG\to\GL(\CV)$ be representation over $\CO_E[1/n]$ with generic fiber $G\to\GL(V)$. Then for all primes $\Fp$ of $\CO_E[1/n]$ with residue field $k_\Fp$ such that $\Char k_\Fp>\height_G(V)$, the base change $\CG_{k_\Fp}\to \GL(\CV_{k_\Fp})$ is of low $\Char k_\Fp$-height. 
\end{Ex}

\begin{Lem}\label{Lem-HeightUnderIsog} 
\begin{enumerate}
\item If $G=\prod G_i$ for reductive groups $G_i$, and $V=\bigotimes_i V_i$ is an (external) tensor product of representations $V_i$ of $G_i$, then $\height_G(V)=\sum \height_{G_i}(V_i)$.
\item If $\pi\colon G\to G'$ is a central isogeny of reductive $k$-groups, $V'$ is a representation of $G'$, and $V$ the representation of $G$ on $V'$ via $\pi$, then $\height_G(V)=\height_{G'}(V')$.
\end{enumerate}
\end{Lem}
\begin{proof}
Part (a) is a simple exercise left to the reader. Concerning (b), we may assume that $k=\overline k$. Let $T'$ be a maximal torus of $G'$ and $T=\pr^{-1}(T')$, denote by $X$, $X'$ the character groups of $T$ and $T'$, by $X^\vee$ and $X^{\prime,\vee}$ their cocharacter groups, and let $R$ and $R'$ be roots of $G$ with respect to $T$ and $G'$ with respect to $T'$ respectively. Let further $f\colon X'\to X$ be the homomorphisms of character groups induced by $\pi|{T}\colon T\to T'$ and $f^\vee\colon X^\vee\to X^{\prime,\vee}$ be the dual map. Then by \cite[9.6, formula (44)]{Springer} we have $f^\vee(R^{+,\vee})=R^{\prime+,\vee}$ if we define $R'$ by a Borel subgroup $B'\subset G'$ containing $T'$ and $R$ by $B=\pi^{-1}(B')$. Moreover a weight $\lambda'\in X'$ of $V'$ becomes the weight $\lambda'\circ\pi$ of $V$. By \cite[9.6.1]{Springer} for $\alpha\in R$ one has $\langle \lambda',f^\vee(\alpha)\rangle=\langle\lambda\circ\pi,\alpha\rangle$. Summing over $\alpha\in R^+$ and maximizing over the weights $\lambda'$ of $V'$, part (b) of the lemma follows.
\end{proof}

\begin{Def}\label{Def-Coxeter}
If the root system of $G_{\overline k}$ is irreducible, the {\em Coxeter number $h_G$} is defined to be
\[h_G=1+ \sum n_i, \]
where $\sum n_i\alpha_i$ is the highest root of $\Lie G_{\overline k}$ (with $(\alpha_i)_{i\in I}$ a root base), cf.~\cite[VI.1.Prop.~31]{Bourbaki}. For a general reductive $G$ one defines \[h_G=\max\{h_H\mid H \hbox{ is a simple quotient of }G\}.\]
\end{Def}
\begin{Rem}\label{Rem-Coxeter}
If the root system of $G_{\overline k}$ is irreducible, one has $h_G=\langle\rho,\beta^\vee\rangle +1$  where $\rho$ is the half sum of the positive roots and $\beta^\vee$ is the highest short coroot, see \cite[\S~3.5]{Humphreys-Modular}. 

Note also that for $G=\GL(V)$ one has $h_G=\dim V$ for its Coxeter number and thus for $G=\prod_i\GL(V_i)$ one finds $h_G=\max \dim V_i$.
\end{Rem}

Let $\Fg^\nilp$ (resp.\ $G_u$) be the reduced closed subscheme of $\Fg=\Lie G$ (resp. $G$) with points the nilpotent (resp.\ unipotent) elements. Let $U$ be the unipotent radical of a parabolic subgroup $P$ of $G$ and $\Fu=\Lie U$. Suppose that $\ell \ge h_G$, so that for the $\ell$-th step of the descending central series of $\Fu$ one has $Z^\ell\Fu=0$. Then the Campbell- Hausdorff group law $\circ$ makes sense in characteristic $\ell$, and it turns $\Fu$ into an algebraic group over $k$. Furthermore there is an unique isomorphism
\begin{equation}\label{Eqn-ExpU}
\exp_U\colon (\Fu,\circ)\stackrel\simeq\longto U
\end{equation}
equivariant for the action of $P$ and whose differential at the origin is the identity; see \cite[Proof of Prop.~5.3]{Seitz} and \cite[Lect.~2]{Serre-Moursund}. Here $k$ can be any field.

If $\ell>h_G$, there is a unique $G$-equivariant isomorphism
\begin{equation}\label{Eqn-ExpG}
\exp_G\colon \Fg_{\overline k}^\nilp \stackrel\simeq\longto G_{\overline k,u},
\end{equation}
which induces (\ref{Eqn-ExpU}) on the unipotent radical of each parabolic subgroup; see \cite[Thm.~3, p.~21]{Serre-Moursund}, \cite[Proof of Prop.~5.3]{Seitz} and \cite[Prop.~11, p.~6]{McNinch-Optimal} for details. Let $\log_G\colon G_u\to \Fg^\nilp$ denote its inverse. The map $\exp_G$ is not necessarily defined over the field of definition $k$ of $G$. For $u\in G_u(k)$, one defines the ``t-power map'' 
\begin{equation}\label{Eqn-Tpower}
\BG_a\to G, t\mapsto u^t= \exp_G(t \log_G u).
\end{equation}
For $G=\GL(V)$ the maps $\exp_G$ and $\log_G$ are the maps in \autoref{Eqn-ExpLogNDef}, and $t\to u^t$ is the the map \autoref{Eqn-DefUtForN}.
\begin{Def}[{\cite[Lect.~3]{Serre-Moursund}}]
A closed $k$-subgroup $H\subset  G$ is called {\em saturated} if for all $u\in H_u(\overline k)$ and $t\in \overline k$ the elements $u^t=\exp_G(t\log_G u)$ lie in $H(\overline k)$. 

For any closed subgroup $\Gamma$ of $G$, its {\em $G$-saturation $\Gamma^{G\dash\sat}_k$} is defined as the intersection of all closed saturated $k$-subgroups of $G$ that contain~$\Gamma$. 
\end{Def}
\begin{Rem}
For $G=\GL_n$ and a closed subgroup $H\subset G$ we clearly have $H^\sat_k=H^{\GL_n\dash\sat}_k$.
\end{Rem}
\begin{Lem}\label{Lem-ExpFieldOfDef}
Suppose $G$ is quasi-split over $k$. Then $\exp_G$ is defined over~$k$.
\end{Lem}
\begin{proof}
Let $B$ be a closed subgroup of $G$, such that its base change $B_{\overline k}$ to $\overline k$ is a Borel subgroup of $G_{\overline k}$. Let $U$ be the unipotent radical of $B$, and so $U_{\overline k}$ is the unipotent radical of $B_{\overline k}$, and let $\exp_U\colon \Lie U\to U$ be the unique exponential map in \autoref{Eqn-ExpU} provided by \cite[\S~5]{Seitz} and defined over $k$. Consider the subvariety 
\[ Y=\{ (g,\Ad_{\Lie G}(g)n,\Ad_G(g)(\exp_U(n))\mid g\in G, n\in \Lie U\}\subset G\times\Lie G\times G\}. \]
The graph of $\exp_U$ is closed in $\Lie G\times G$ and defined over $k$, and since the adjoint action of $G$ on $\Lie G\times G$ (by $\Ad_{\Lie G}\times\Ad_G$) is also defined over $k$, its orbit $Y$ is defined over $k$, as well. Let $Z\subset \Lie G\times G$ be the image of $Y$ under the projection onto the last two factors of $G\times\Lie G\times G$. Note that $(\Lie G_{\overline k})^\nilp$ is the orbit of $\Lie U_{\overline k}$ under the adjoint action. Hence after base change from $k$ to $\overline k$, the uniqueness of $\exp_G$ in \autoref{Eqn-ExpG} implies that $Z_{\overline k}$ is the graph of $\exp_G$. It follows that $Z$ is the graph of a morphism defined over $k$, and hence that $\exp_G$ itself is defined over~$k$.
\end{proof}
\begin{Rem}
The main issue about $\exp_G$ being defined over $k$ is its domain. If, for instance, $G$ was the unit group of a quaternion algebra, then $G$ has no non-trivial parabolics defined over $k$ and hence the result of Seitz (or Serre) is void.
\end{Rem}
\begin{Cor}
If $G$ is quasi-split, then the intersection of all saturated subgroups of $G_{\overline k}$ that contain $H$ agrees with $H^{G\dash\sat}_k\otimes_k\overline k$.
\end{Cor}
\begin{proof}
Because of \autoref{Lem-ExpFieldOfDef} the isomorpism $\exp_G$ and its inverse $\log_G$ are defined over $k$. Now following the proof of \autoref{Lem-SaturationIsIntersection} implies the result of the corollary.
\end{proof}

The following are some important results from \cite[Lect.~4]{Serre-Moursund}:
\begin{Cor}[{\cite[Cor.~2 to Thm.~5]{Serre-Moursund}}]\label{Cor-ExpGVsExpN}
The representation $V$ is of low $\ell$-height for $G$ if and only if the following diagram is well-defined and it commutes
\[\xymatrix{
\Fg^\nilp\ar[d]_{\exp_G} \ar[r]^-{\mathrm{d}\rho_V} & \End(V)^\nilp\ar[d]^{\exp_{\GL(V)}}\\
G_u\ar[r]^-{\rho_V}& \GL(V)^\unip\rlap{.}\\
}\]
\end{Cor}
In other words, if $V$ has small $\ell$-height for $G$, then the saturation inside $G$ defined by Serre agrees with the saturation inside $\GL_n$ defined by Nori.
\begin{Cor}[{\cite[Cor.~1 to Thm.~5]{Serre-Moursund}}]\label{Cor-SaturatedSubgroupAndHeight}
If $H$ is a closed, reductive and saturated subgroup of $G$, one has $\height_H(V)\le\height_G(V)$.
\end{Cor}
The following corollary generalizes \cite[App.~B, Prop.~25, 26]{EHK} from $\BF_\ell$ to any finite field.
\begin{Cor}\label{Cor-LowHeightForSubgroups}
Let $H$ be a closed, connected reductive and saturated subgroup of $\GL_n$, and let $V$ be the standard representation of $\GL_n$. Denote by $\lambda^+$ and $\lambda^-$ the highest and  lowest weight of $V$ regarded as a representation of $H$ with respect to a root base $\Delta$. Let $\eps$ be in $\{\pm\}$. Then $\eps\lambda^\eps$ is dominant, and writing $\eps\lambda^\eps=\sum_{\alpha\in\Delta} m^\eps_\alpha\alpha$, one has $\max\{m^\eps_\alpha\mid \alpha\in\Delta,\eps\in\{\pm\}\}\le n-1$.
\end{Cor}
\begin{proof}
By \cite[VIII, Sec.~7.5, Prop.~11]{Bourbaki-Lie7-8} the lowest weight of an irreducible representation of $G$ is the negative of a dominant weight. This shows that $-\lambda^-$ is dominant, and for $\lambda^+$ this is clear. In particular the $m_\alpha^\eps$ lie in $\BQ_{\ge0}$. Now from \autoref{Cor-SaturatedSubgroupAndHeight} and \autoref{Rem-Dynkin} we deduce that $\max\{m_\alpha+m_\alpha^{-}\mid \alpha\in\Delta\}\le n-1$, and this completes the proof of the corollary.
\end{proof}
The following result and its proof are inspired by \cite[\S1, 2]{Friedlander} where similar but not quite the same results are developed.
\begin{Prop}\label{Prop-WeilResAndSat}
Let $G$ be quasi-split with $h_G<\ell$, let $k'$ be a subfield of $k$ such that $k/k'$ is finite separable. Then the following hold:
\begin{enumerate}
\item The Weil restriction $G':=\Res_{k/k'}G$ is quasi-split, $h_{G'}=h_G$, and $\exp_{G'}$ is defined~over~$k'$.
\item One has $\Res_{k/k'}\exp_U=\exp_{G'}|_{\Res_{k/k'}U}$ for any parabolic $P\subset G$ with unipotent radical $U$.
\item Let $H\subset G$ be a closed subgroup. If $H$ is saturated in $G$, then so is $\Res_{k/k'}H$ in $\Res_{k/k'}G$. For any $H$ we have an inclusion $(\Res_{k/k'}H)^\sat_{k'}\subset \Res_{k/k'}(H^\sat_{k})$.
\item Let $V$ be a vector space over $k$ and write $W$ for $V$ as a vector space over $k'$. Let $\iota\colon\Res_{k/k'}\GL(V_k)\to\GL(W_{k'})$ be the embedding defined by the natural transformation $\GL_{k\otimes_{k'}R}(V\otimes_{k'}R)\to \GL_R(W\otimes_{k'}R)$ on $k'$-algebras~$R$. Then:
\begin{enumerate}
\item We have $\height_{\Res_{k/k'}\GL(V)} W=[k:k'](\dim_k V-1)=\dim_{k'}W-[k:k']$.
\item  If $\ell>\dim_{k'}W-[k:k']$, then the following diagram commutes
\[\xymatrix{
(\Res_{k/k'}\End_{k}(V))^\nilp\ar[d]_{\exp_{\Res_{k/k'}\GL(V)}} \ar[r]^-{\mathrm{d}\iota} & \End_{k'}(W)^\nilp\ar[d]^{\exp_{\GL(W)}}\\
\Res_{k/k'}\GL(V_k)^\unip\ar[r]^-{\iota}& \GL(W_{k'})^\unip\rlap{.}\\
}\]
\item If $H$ is closed and saturated in $\GL(V)$ and if $\ell>\dim_{k'}W-[k:k']$, then $\Res_{k/k'}H$ is saturated inside $\GL(W_{k'})$ (embedded via $\iota$).
\end{enumerate}
\end{enumerate}
\end{Prop}
\begin{proof}
Note first that after base change to an algebraic closure $\overline k$ of $k$ we have an isomorphism
\begin{equation}\label{Eqn-WeilResIso}
\Res_{k/k'}G\otimes_k\overline k\cong \prod_{\sigma\in\Hom_{k'}(k,\overline k)} G\otimes_k^\sigma\overline k,
\end{equation}
where the map $k\to\overline k$ in the factor for $\sigma$ is given by $\sigma$. Using \autoref{Eqn-WeilResIso}, it is now easy to see that if $P$ is a $k$-parabolic subgroup of $G$, then $\Res_{k/k'}P$ is a $k'$-parabolic subgroup of $\Res_{k/k'}G$. This implies the first claim of (a). The second claim of (a) also follows from \autoref{Eqn-WeilResIso} and the definition of $h_G$. The last claim in (a) follows from the first two and \autoref{Lem-ExpFieldOfDef}. For (b) note that, again by \autoref{Eqn-WeilResIso}, $\Res_{k/k'}U$ is the unipotent radical of $\Res_{k/k'}P$. Now (b) follows from the characterization of $\exp_{U}$ given in the paragraph surrounding the equation \autoref{Eqn-ExpU} and the properties of the Weil restriction. By (b) the exponential $\exp_{\Res_{k/k'}G}$ after base change to $\overline k$ and restriction to the factor $\sigma$ in \autoref{Eqn-WeilResIso} is given by $\exp_G\otimes_k^\sigma\overline k$. Since clearly $H\otimes_k^\sigma\overline k$ is saturated in $G\otimes_k^\sigma\overline k$ under $\exp_G\otimes_k^\sigma\overline k$, this implies the first assertion of (c); the second follows from the first. Part (d)(i) follows from \autoref{Eqn-WeilResIso} and \autoref{Rem-Dynkin}. Part (d)(ii) is a consequence of (i)  and \autoref{Cor-ExpGVsExpN}. Part (d)(iii) is a consequence of (c) and~(d)(ii).
\end{proof}
\begin{Ex}\label{Ex-SatAndWeil}
It is in general not true that the operation of taking the saturation and of taking a Weil restriction commute. For this, suppose $k$ is a finite field with proper subfield $k'$, let $G=\SL_{n,k}$ and let $H$ be the finite group $\SL_n(k')$ considered as a reduced finite subgroup scheme. Then we have $H^\sat_k=G$, because the Nori envelope has to contain all root groups for the diagonal torus and so in particular the standard upper triangular Borel subgroup and its opposite. This implies that $\Res_{k/k'}(H^\sat_k)=\Res_{k/k'}G$. However $\Res_{k/k'}H=H$, and therefore $(\Res_{k/k'}H)_k^\sat=\SL_{n,k'}$ inside $\Res_{k/k'} \SL_{n,k'}$ via the the embedding $\SL_{n,k'}\into \Res_{k/k'} \SL_{n,k'}$ that is the adjunction map of the Weil restriction. We leave further details to the reader.
\end{Ex}

\subsection{Saturation and Lifting}\label{Subsec-OnLifting}
In this subsection, in \autoref{Prop-LiftingOfLowWeightReps}, we prove a result on lifting representations of low $\ell$-height of a saturated semisimple subgroup $G$ of $\GL_n$ from a finite field of characteristic $\ell>0$ to its ring of Witt vectors and thus to a finite extension of $\BQ_\ell$. 

We begin with some elementary preparations for the proof of \autoref{Prop-LiftingOfLowWeightReps}.
\begin{Lem}\label{Lem-TensorProdOfReps}
Let $k$ be a field. Let $H\subset\GL_{n,k}$ be a split semisimple connected closed subgroup. Denote by $V\cong k^n$ the resulting representation $\bar\rho\colon H\to\Aut_k(V)$ of the algebraic group $H$, and by $\pi_k\colon \wt H\to H$ the universal cover of $H$, so that $\wt H$ is semisimple connected and simply connected. Suppose that $V$ is irreducible. Then
\begin{enumerate}
\item $\wt H\cong\prod_{i\in I} \wt H_i$ with $\wt H_i$ split simple semisimple simply connected over~$k$.
\item Let $\wt T_i$ be a maximal split torus of $\wt H_i$. Then $\wt T=\prod_i\wt T_i$ is a maximal split torus of $\wt H$ and $T:=\pi_k(\wt T)$ of $H$. Let $\iota\colon \oplus_I X(\wt T_i)\!\stackrel\simeq\to \!X(\wt T)$ be the induced isomorphism on character~groups.
\item As a representation of $\wt H$ via $\bar\rho\circ\pi_k$ one has $V\cong\otimes_{i\in I} V_i$ (exterior tensor product) for irreducible representations $\rho_i\colon \wt H_i\to\Aut_k(V_i)$.
\item If $\mu_i\in X(\wt T_i)$ is the highest weight of $V_i$, then $\iota((\mu_i)_{i\in I})$ is the highest weight of $V$.
\item The order of $\kernel\pi_k$ is bounded by $n=\dim V$ and if $\Char k=0$ or $\Char k>n$, then $\kernel\pi_k$ is \'etale over~$k$.
\end{enumerate}
\end{Lem}
\begin{proof}
Part (a) follows from \cite[parag.~before 6.4.4 and Thm.~5.1.19]{ConradSGA3}, and (b) is clear. Parts (c) and (d) are standard facts. But we could not find a reference. Let $\mu$ be the highest weight of $V$ and define $(\mu_i)_{i\in I}:=\iota^{-1}(V)$ and $V_i$ as the irreducible $\wt H_i$ representation with highest weight $\mu_i$. Then (the external tensor product) $W:=\otimes_{i\in I} V_i$ has highest weight $\mu$ as a representation of $\wt H$. We show that $W$ is irreducible as a representation $\wt H(\overline k)$: Observe that by the Burnside theorem and the irreducibility of $V_i$ as a representation of $\wt H_i$ the homomorphism $\overline k[\wt H_i(\overline k)]\to \Aut_{\overline k}(V_i(\overline k))$ induced from the representation of $V_i$ is surjective. But from our product situation it is then clear that also the homomorphism $\overline k[\wt H(\overline k)]\to \Aut_{\overline k}(V(\overline k))$ induced from the representation of $V$ is surjective. Again from Burnside, we deduce that $W(\overline k)$ and hence $W$ is irreducible. By the classification of irreducible representations in terms of their highest weight for split reductive groups, see \cite[II.~Cor.~2.7]{Jantzen-Representations}, we find $V\cong W$, and hence (c) and (d) are shown.

To prove (e), we may pass to an algebraic closure $\overline k$ of $k$. Since $\kernel\pi_k$ is contained in the product of the centers of the $\wt H_i$, it suffices to show that the order of the center (as a group scheme) will be bounded by $\dim V_i$. Then the order of $\kernel\pi_k$ is bounded by $\prod_i\dim V_i=\dim V$. However for $G$ simply connected simple, over any field, the order of its center is always bounded by the smallest dimension of a non-trivial representation, as follows from the classification: The kernel has order at most $m$ if $G=\SL_m$, $2$ if $G$ is of type $B_m$, $C_m$ or $E_7$, $4$ if $G$ is of type $D_m$, $3$ if $G$ is of type $E_6$ and $1$ otherwise, e.g.~\cite[Planches, 250ff., entry (VIII)]{Bourbaki}. The minimal dimensions of non-trivial representations of $G$ are given in \cite[VIII, table 2, p.~214]{Bourbaki-Lie7-8}; they are $m$ for $\SL_m$, and $2m+1$, $2m$, $2m$, $27$, $56$, $248$, $26$, $7$ for $B_m$, $C_m$, $D_m$, $E_6$, $E_7$, $E_8$, $F_4$ and $G_2$, resp. The \'etaleness of $\kernel\pi_k$ now follows from $\Char k>n$ or $\Char k=0$.
\end{proof}
\begin{Rem}\label{Rem-OnTensorDecompForAdjoint}
The scheme theoretic intersection $Z_i:=\kernel\pi_k\times_{\wt H} \wt H_i$ is the kernel of $\wt H_i\to H$, so that $H_i:=\wt H_i/Z_i$ can be regarded as an almost simple subgroup of $H$. Moreover $Z_i$ acts trivially on $V_i$, so that $V_i$ can be regarded as a representation $V'_i$ of $H_i$. One can form $\bigotimes_iV'_i$ over $H':=\prod_i H_i$, and $\prod_i H_i\to H$ is a central isogeny whose kernel coincides with that of $\bigotimes_i V'_i$. Thus in the above lemma, one can replace, as often done in the literature, the group $\wt H$ by~$H'$.

Moreover if $H$ has trivial center, i.e. $H=H^\ad$, then $H'=H$. Then $H=\prod_i H_i$ for simple connected adjoint groups, and the decomposition of $V=\bigotimes_iV_i$ can be carried out with irreducible representations $V_i$ of~$H_i$. In this case (as in the case where $H$ is simply connected), the $V_i$ can also be regarded as representations of $H$ itself.
\end{Rem}
\begin{Rem}
By reduction to the simply connected case, the argument for (e) shows that for any connected semisimple group $G$ with center $Z$ over any field with an almost faithful representation $V$ one has $\dim V\ge\#Z$.
\end{Rem}
The next result makes use of the notion of Weyl module. Weyl modules are natural representations $V(\lambda)_\BZ$ attached to a fixed dominant weight $\lambda$ of a reductive group over $\BZ$. Via base change they give rise to representations $V(\lambda)_A$ of the base changed group over any base ring $A$. If $A$ is a field of characteristic zero, then $V(\lambda)_A$ is irreducible of highest weight $\lambda$. If $A$ is a field of positive characteristic, in general $V(\lambda)_A$ will be reducible. For more see \cite[II.8]{Jantzen-Representations}.
\begin{Prop}\label{Prop-LiftingOfLowWeightReps}
Let the hypotheses and the notation be as in \autoref{Lem-TensorProdOfReps}. Suppose further that $k$ is finite of characteristic $\ell>2(n-1)$ and that $H\subset\GL_n$ is saturated, and denote by $W(k)$ the ring of Witt vectors of $k$. Then
\begin{enumerate}
\item Each $V_i$ is a Weyl module for $\wt H_i$.
\item There exist semisimple $W(k)$-group schemes $\CH$, $\wt\CH$, $\wt\CH_i$ with special fiber $H$, $\wt H$ and $\wt H_i$, respectively; moreover $\wt\CH$ is simply connected, we have a universal cover $\pi\colon\wt\CH\to\CH$ over $W(k)$ whose reduction is $\pi_k$, and we have an isomorphism $\wt\CH\cong\prod_i\wt\CH_i$ over $W(k)$ whose reduction is the isomorphism in \autoref{Lem-TensorProdOfReps}(b).
\item Let $\CV_i$ be the Weyl module for $\wt \CH_i$ over $W(k)$ with the same highest weight than $V_i$. Then $\CV:=\otimes_{i\in I}\CV_i$ is a representation $\rho$ of $\wt\CH$ which reduces to $\bar\rho\circ\pi_k$. 
\item For $(\CV,\rho)$ from (c) we have $\kernel\rho=\kernel\pi$, and thus $\CV$ is a representation of $\CH$ over $W(k)$ that reduces to the representation $V$ of $H$.
\end{enumerate}
\end{Prop}
\begin{proof}
To prove (a) we shall verify the low alcove condition for $V_i$ with respect to $\wt H_i$ given in  \cite[Cor.~4.4.3]{McNinch-Dimensional} in order to deduce by loc.cit.\ that $V_i$ is a Weyl module. Let $\Phi_i\subset X(\wt T_i)$ be the root system of $\wt H_i$, choose a Borel subgroup $B_i\subset\wt H_i$ and let $\Delta_i$ be the corresponding root base for $\Phi_i$. Because $H$ is saturated in $\GL_n$ we deduce $\height_H(V)\le n-1$ from \autoref{Cor-LowHeightForSubgroups}, and from \autoref{Lem-HeightUnderIsog} we deduce $\height_{\wt H_i}(V_i)\le \height_{\wt H}(V)=\height_H(V)\le n-1$, using for the $\height_{\wt H_j}(V_j)$, $j\neq i$, that the Dynkin height is always positive. 

We need to collect some results about the paring $\langle\cdot,\cdot\rangle\colon X(\wt T_i)\times X^\vee(\wt T_i)\to\BZ$. Denote by $\sigma_i$  the half sum of positive roots of $\wt H_i$, and by $\alpha_{i,0}$ the highest short root of $\wt H_i$, see \cite[p.~251]{Carter-Affine}. Then from \autoref{Rem-Coxeter} we have $\langle\sigma_i,\alpha^\vee_{i,0}\rangle=h_{\wt H_i}-1$ for the Coxeter number of $\wt H_i$. For any simply connected simple group $G$ the value of $h_G$ is at most as large as the dimension of the smallest irreducible non-trivial representation of $G$ (see the proof of \autoref{Lem-TensorProdOfReps}(e)). Moreover from  \cite[p.~543ff.]{Carter-Affine} one can check that $\langle\alpha,\alpha^\vee_{i,0}\rangle\in\{0,1\}$ for any $\alpha\in\Delta_i$ and that the value $1$ is attained at least once and at most twice (this holds for any simply connected simple~$G$). 

Let $\mu_i$ be the highest weight of $V_i$ with respect to $\wt T_i$. Write  $\mu_i=\sum_{\alpha\in \Delta_i} m_{i,\alpha} \alpha$ for a unique tuple $(m_{i,\alpha})\in\BQ^{\Delta_i}_{\ge0}$. Then from $\height_{\wt H_i}(V_i)\le n-1$ and from \autoref{Cor-LowHeightForSubgroups} it follows that $\sum_{\alpha\in \Delta_i} m_{i,\alpha}\le n-1$. Then by the above we find
\[ \langle \sigma_{i}+\mu_{i},\alpha^\vee_{i,0} \rangle\le h_{\wt H_i}-1+\sum_{\alpha\in \Delta_{i}} m_{i,\alpha} \langle\alpha,\alpha^\vee_{i,0} \rangle 
\le \dim V_i-1+ \sum_{\alpha\in \Delta_i} m_{i,\alpha} \le 2(n-1).\]
Because $\ell>2(n-1)$, it follows now from \cite[Cor.~4.4.3]{McNinch-Dimensional}, that the Weyl module $\bar V(\mu_i)$ for the weight $\mu_i$ of $\wt H_i$ is simple, and hence isomorphic to $V_i$, which proves~(a).

Now we turn to (b). The existence of simply connected simple groups $\CH$, $\wt\CH$, $\wt\CH_i$ with special fiber $H$, $\wt H$ and $\wt H_i$, as well as the isogeny $\pr$ lifting $\pr_k$ follows from \cite[Thm.~6.1.16]{ConradSGA3}. Note also that the isomorphism $\wt H\cong\prod_{i\in I} \wt H_i$ defines an isomorphism of root data which yields an isomorphism between the root data for $\wt\CH$ and $\prod_i\wt\CH_i$, and then again by  \cite[Thm.~6.1.16]{ConradSGA3} we obtain an isomorphism $\wt\CH\cong\prod_i\wt\CH_i$ over $W(k)$ whose reduction is the isomorphism in  \autoref{Lem-TensorProdOfReps}(b).

Regarding (c), note first that as explained in \cite[II.8]{Jantzen-Representations} Weyl modules for split reductive connected groups are defined over $\BZ$, and hence over any ring. Thus we have Weyl modules $\CV_i$ for $\wt\CH_i$ defined over $W(k)$ whose reduction to $k$ is $V_i$. The highest weighs of $\CV_i$ and of $V_i$ are the same. We remark that Weyl modules are irreducible over any field of characteristic zero, and that their highest weight is independent of the characteristic. From their definition in loc.cit.\ and arguing as in the proof of \autoref{Lem-TensorProdOfReps}(c),(d), it follows that in fact $\CV=\bigotimes_i\CV_i$ is the Weyl module for $\wt\CH$ of highest weight $\iota((\mu_i)_{i\in I})$. Its reduction agrees with $V=\bigotimes_iV_i$, and this proves~(c). 

Regarding (d), let $\wt\CT$ denote a maximal split torus of $\wt \CH$ over $W(k)$ whose image under $\pi$ is a maximal split torus $\CT$ of $\CH$. Since the kernel of $\pi$ is equal to the kernel of $\pr\colon\wt\CT\to\CT$, it suffices to consider the restriction of $\rho$ to $\wt\CT$. This restriction can be regarded as an $n$-tuple of characters $(\chi_i)_{i=1,\ldots,n}$ in $\Hom(\CT,\BG_{m,W(k)})$, and thus $\kernel\pi=\bigcap_{i=1,\ldots,n}\kernel\chi_i$. Their reductions form the corresponding character tuple $(\bar\chi_i)_{i=1,\ldots,n}$ for $\bar\rho|_{\wt T}$. After identifying $\CT$ with $\BG_{m,W{k}}^r$ for some $r\ge1$, the tuples are represented by some $r\times n$ matrix over $\BZ$, and the order of the kernel is the product of the elementary divisors of this matrix. By \autoref{Lem-TensorProdOfReps}(e), over $k$, this order is bounded by $n<\ell$ (because $n>1$ and $\ell>2(n-1)$). In particular, this kernel is finite \'etale over $W(k)$. Since $\pr$ from  \cite[Thm.~6.1.16]{ConradSGA3} is constructed from the corresponding map of root data over $k$, and since over the special fiber $\kernel\rho$ agrees with $\kernel\pi_k$, we deduce $\kernel\pi=\kernel \rho$, and this completes~(d).
\end{proof}

\subsection{Saturatedness of some non-reductive groups}
\label{Subsec-NonReduSat}

We now assume that $k$ is finite, that $G$ is a connected smooth affine $k$-subgroup scheme of $\GL_{n,k}$ but not necessarily reductive, and that $G$ contains a maximal $k$-split torus $T$. Let $V=k^n$ be the standard representation of $\GL_{n,k}$. Let $\Lie G\cong\Lie T\oplus \bigoplus_{\alpha\in \Phi} \Fg_\alpha$ denote the eigenspace decomposition under the adjoint action of $T$ where $\Phi\subset X(T)\setminus\{0\}$ is the finite subset of characters $\alpha$ of $T$ whose weight space $\Fg_\alpha$ is non-zero. Throughout this subsection, we require the following conditions to hold: 
\begin{enumerate}
\item[(i)] the weight spaces $\Fg_\alpha$ are $1$-dimensional, 
\item[(ii)] for any $\alpha\in \Phi$ we have $\BZ \alpha\cap \Phi\subset\{\pm \alpha\}$ (i.e., $\Phi$ is reduced).
\item[(iii)] for any $\alpha\in \Phi$, we have
\[  \{i\in\BN\mid \alpha^i \hbox{ is a weight of the action of $T$ on $\End(V)$} \}\subseteq \{1,\ldots,\ell-1\}.\]
\end{enumerate}
Note that if $G\subset \GL_{n,k}$ is a reductive subgroup, then (i) and (ii) clearly hold and (iii) holds if $\End(V)$ is $\ell$-restricted.

\begin{Lem}
For $\alpha\in \Phi$, let $T_\alpha\subset T$ be the subtorus $\kernel \alpha$ and define $G_\alpha:= Z_G(T_\alpha)$. In the following we consider $(\Fg_\alpha,+)$ as a group scheme isomorphic to $\BG_a$. Then
\begin{enumerate}
\item $G_\alpha$ is a smooth subgroup scheme of $G$ and $\Lie(G_\alpha)=\Lie T\oplus \Fg_\alpha\oplus\Fg_{-\alpha}$.
\item If $\Fg_{\beta}\neq0$ for $\beta\in\{\pm\alpha\}$, then there exists a unique $k$-homomorphism $u_\beta\colon\Fg_\beta\to G_\alpha$, such that $\mathrm{d}u_\beta=\id_{\Fg_\beta}$ and $u_\beta(\lambda t)=\beta(\lambda)u_\beta(t)$ for all $t\in\overline k$. Let $U_\alpha$ be the image of $u_\alpha$.
\item If $\dim \Fg_\alpha\oplus\Fg_{-\alpha}=2$, then $G_\alpha$ is either reductive of rank $1$, or $R_u(G_\alpha)\cong U_\alpha\times U_{-\alpha}$, and $G_\alpha/R_u(G_\alpha)\cong T$. 
\item If $\dim \Fg_{-\alpha}=0$, then $G_\alpha$ is solvable and sits in an extension $1\to U_\alpha\to G_\alpha\to T\to1$. 
\end{enumerate}
\end{Lem}
\begin{proof}
Part (a) is a special case of \cite[Lem.~2.4.4]{ConradSGA3}, using condition (ii). If $G_\alpha$ is reductive, then assertions (b) and (c) are clear. In the other case $R_u(G_\alpha)$ is non-trivial, and $G_\alpha/R_u(G_\alpha)$ must be a reductive group of rank $0$, i.e., a torus. For dimension reasons, the induced map from $T$ to that torus is an isomorphism and thus $G_\alpha$ sits in a short exact sequence
$1\to R_u(G_\alpha)\to G_\alpha\to T\to1$. If $\dim R_u(G_\alpha)=1$, then we must have $\Lie R_u(G_\alpha)\cong \Fg_\alpha$ and hence the existence of $u_\alpha$ follows from \cite[14.3.11]{Springer}, and (d) follows. Because of hypotheses (i) and (ii) in the remaining case $\Lie R_u(G_\alpha)\cong \Fg_\alpha\oplus\Fg_{-\alpha}$ has dimension $2$. Since $k$ is perfect, we can regard $R_u(G_\alpha)$ as a $2$-dimensional connected subgroup of $U_n$. \cite[14.3.11]{Springer} guarantees the existence of subgroups isomorphic to $\BG_a$ on which $T$ acts via $\alpha$ and $-\alpha$, respectively. The existence of the $u_\alpha$ and hence $U_\alpha$ are now straightforward: One simply has to check that any action of $\BG_m$ on $\BG_a$  is of the form $\BG_m\times\BG_a\to\BG_a,(\lambda,t)\mapsto \lambda^it$ for some unique $i\in\BZ$ (and this generalizes in an obvious way to actions of $T$ on $\BG_a$). To see that $U_\alpha$ and $U_{-\alpha}$ commute, note that their commutator lies in $R_u(Z_G(T_\alpha))$ but has weight $0$ for the action of $T$, and hence is trivial.
\end{proof}

\begin{Lem}\label{Lem-OnExp}
Let $\alpha$ be in $\Phi$ and let $u_\alpha\colon\Fg_\alpha\to U_\alpha$ as in the previous lemma. Then for any $u\in U_\alpha(k)$ and corresponding $X=u_\alpha^{-1}(u)$ one has
\[\exp_n(t\log_n (u))=u_\alpha(t X).\]
 \end{Lem}
\begin{proof}
The following argument is inspired by \cite[proof of Thm.~5 and its Cor.~2]{Serre-Moursund}. For explicitness we may assume that the $k$-split solvable group $TU_\alpha$ is contained the the upper triangular Borel subgroup of $\GL_{n,k}$, and $U_\alpha$ in $U_{n,k}$. Consider the homomorphism $\phi_{\alpha,t}\colon\BG_a\to U_{n,k},t\mapsto u_\alpha(t)$. It can be given as a polynomial in $t$ of the form $1+\sum_{i\ge1}a_i(\alpha)t^i$, for suitable $a_i(\alpha)\in \End(V)$ ($V=k^n$). Consider the adjoint action of $\lambda\in T(\overline k)$:
\begin{eqnarray*}
1+\sum_{i\ge1}a_i(\alpha)\alpha^i(\lambda ) t^i& =&1+\sum_{i\ge1}a_i(\alpha)(\alpha(\lambda )t)^i \ = \ u_\alpha(\alpha(\lambda )t) \ = \ \lambda u_\alpha(t)\lambda^{-1}\\
&=&\lambda \Big(1+\sum_{i\ge1}a_i(\alpha)t^i\Big)\lambda^{-1} \ = \ 1+\sum_{i\ge1}\lambda a_i(\alpha) \lambda^{-1} t^i.
\end{eqnarray*}
It yields $a_i(\alpha)\alpha^i(\lambda )= \lambda a_i(\alpha) \lambda^{-1}$ for all $i\in\BN$, and thus either $a_i(\alpha)=0$ or $\alpha^i$ is a $T$-weight of $\End(V)$. Since by hypothesis (iii) the $T$-weights of $\End(V)$ are $\ell$-restricted, it follows that $a_i(\alpha)=0$ for $i\ge \ell$. Now consider
\[\exp_n(t\log_n (u))-u_\alpha(t X)\]
as a polynomial in $t$ with coefficients in $\End(V)$. By what we just proved and by the definition of $\exp_n$ it follows that this polynomial has degree at most $\ell-1$. As it vanishes for the $\ell$ values $t=0,\ldots, \ell-1$, it is identically zero.
\end{proof}
\begin{Rem}
It seems an interesting question whether one can use the $u_\alpha$ above (and their conjugates) together with the Baker-Campbell-Hausdorff formula to define a Springer type isomorphism between $G_u(\overline k)$ and $(\Lie G_{\overline k})^\nilp$, even if $G$ is not semisimple or reductive. If $G$ is semisimple such a Springer isomorphism was used by Serre in his approach to saturation for general reductive groups. A related question where $G$ is reductive but $\ell<h_G$ is asked in \cite[Rem.~27]{McNinch-Optimal}: does there exist a Springer isomorphism whose restriction to the unipotent radical of any parabolic is as in formula~\autoref{Eqn-ExpU}.
\end{Rem}
\begin{Thm}\label{Thm-GisSaturated}
The group $G$ is saturated.
\end{Thm}
\begin{proof}
Let $H$ be the subgroup of $G$ generated by the $U_\alpha$ (for one fixed $T$). Then by \autoref{Lem-OnQuasiSplit} we have $H=\wt G$. Because of \autoref{Lem-OnExp} each $U_\alpha$ and hence also $H$ is exponentially generated. Since $G$ is connected $G^+=\wt G$ by \autoref{Cor-OnG+}(c), and hence $G$ is saturated by \autoref{Prop-SatOfWtG}.
\end{proof}

\subsection{Reductions of reductive groups and saturation}
\label{Subsec-LPschemesAndSat}
The first aim here is to generalize a result of Larsen-Pink from $\BQ$- to $E$-rational compatible systems. In \autoref{Thm-RedIsSaturated} we deduce from this the saturatedness of an interesting class of a priori not necessarily reductive groups over a finite field. In order to carry out the first part, which is a straightforward extension of results of Larsen and Pink, we recall some notions and results from \cite[Sects.~4,6,7]{LarsenPink92}. We follow this reference rather closely.

Throughout this subsection, we let $K$ be any field of characteristic zero, and we assume that all algebraic groups and varieties will be defined over $K$. Denote by \[\ch\colon\GL_n\to \BG_m\times\BG_a^{n-1}\] the morphism associating to a matrix the coefficients of its characteristic polynomial. Let $D_n=\BG_m^n$ be the diagonal split maximal torus of $\GL_n$. The Weyl group of $\GL_n$ with respect to $D_n$ is the symmetric group $S_n$, acting on $D_n$ by permutation of factors. The restriction of $\ch$ to $D_n$ is a finite morphism which identifies $\BG_m\times\BG_a^{n-1}$ with the scheme-theoretic quotient $D_n/S_n$. 

Let $T_0\subset D_n$ be any closed subtorus. Then $T_0$ is also split over $K$, and as a map between split tori, the inclusion $T_0\subset D_n$ is defined over $\BQ$ (or even $\BZ$), and hence the variety $\ch(T_0)\subset \BG_m\times\BG_a^{n-1}$ is defined over $\BQ$. The variety $\ch(T_0)$ is geometrically connected, because $T_0$ is so. Let $\rho\colon T_0\to D_n$ denote the representation corresponding to $T_0\into D_n$, and consider on $D_n$ the action of $S_n=N_{\GL_n}(D_n)/C_{\GL_n}(D_n)$ that permutes its entries. Set $\Aut(T_0,\rho_0):=\{ \alpha\in\Aut(T_0)\mid \exists \sigma \in S_n: \rho\circ\alpha=\sigma\circ \rho\}$ -- the automorphism group of the pair $(T_0,\rho)$ (cf.~\cite{LarsenPink92}). One has an isomorphism $N_{S_n}(T_0)/C_{S_n}(T_0) \stackrel\simeq\to \Aut(T_0,\rho_0)$, and the pair $(T_0,\rho_0)$ is determined by $\ch(T_0)$ up to isomorphism. (cf.~\cite{LarsenPink92}).

Following \cite[4.4]{LarsenPink92}, for any $\sigma\in S_n\setminus C_{S_n}(T_0)$, we define a proper subgroup $H_\sigma$ of $T_0$: if $\sigma(T_0)=T_0$, then $H_\sigma=\{t\in T_0\mid \sigma(t)=t\}$; otherwise, set $H_\sigma:=T_0\cap \sigma(T_0)$. The intersection of different tori of the same dimension is a subgroup of smaller dimension. Also, the condition $\sigma(t)=t$ can be expressed in the vanishing of at least one character of $T_0$. Thus $\dim H_\sigma<\dim T_0$ for all $H_\sigma$. Since $\ch$ is finite on $D_n$ one deduces $\dim \ch(H_\sigma)=\dim H_\sigma<\dim \ch(T_0)=\dim T_0$. In particular $Y:=\bigcup_{\sigma\in S_n\setminus C_{S_n}(T_0)}\ch(H_\sigma)$ is a proper closed subset of $\ch(T_0)$.

For a connected reductive group $G\subset \GL_n$ there exists a torus $T_0\subset \GL_n$ as above such that every maximal torus of $G$ is, over an algebraic
closure $\overline K$ of $K$, conjugate to $T_0$. Over $\overline K$, the semisimple part of any $g\in G$ can be conjugated into $T_0$, so that $\ch(T_0)=\ch(G)$ pointwise. It follows that $\ch(G)=\ch(T_0)\subset \BG_m\times\BG_a^{n-1}$ is Zariski-closed and that $T_0$ is, up to conjugation by $S_n$, uniquely determined by $\ch(G)$. Hence $Y$ only depends on $\ch(G)=\ch(T_0)$ (cf.~\cite[4.3]{LarsenPink92}).

\begin{Def}[{\cite[\S~4]{LarsenPink92}}]\label{Def-GammaRegular}
Let $G\subset\GL_n$, be a connected reductive subgroup, and let $Y\subset \ch(G)$ be as above. An element $g\in G$ is {\em $\Gamma$-regular} if and only if $\ch(g)\notin Y$.
\end{Def}
The importance of $\Gamma$-regular elements $g\in G$ is explained in \cite[4.4-4.7]{LarsenPink92}: One has:
\begin{enumerate}
\item The element $g$ is regular semisimple in $G$, and it lies in a unique maximal torus $T_g$ of $G$; denote by $\rho_g\colon T_g\into \GL_n$ the induced representation of $T_g$ on $K^n$.
\item The element $g$, via the eigenvalues of $\rho_g(g)$ together with their multiplicities, determine uniquely the formal character of the action of $T_g$ on $K^n$ via $\rho_g$ from (a). 
\end{enumerate}

Let $\rho_\bullet$ be an $E$-rational $n$-dimensional compatible system. Following the proof of \cite[6.11]{LarsenPink92}, the following result is immediate.
\begin{Lem}[{\cite[6.11]{LarsenPink92}}]\label{Lem-CharIndepOfLambda}
Let $g\in\pi_1(X)$ and $\Gamma\subset \pi_1(X)$ be an open subgroup. Then the Zariski closure of $\ch(\rho_\lambda(g\Gamma))$ in $(\BG_m\times\BG_a^{n-1})_{E_\lambda}$ is defined over $E$ and independent of $\lambda\in\CP_E'$.
\end{Lem}
\begin{Rem}
\autoref{Lem-CharIndepOfLambda} also applies to compatible systems over number fields.
\end{Rem}
 
\begin{Con}\label{Conv-RhoConnected}
In the remainder of this subsection, $\rho_\bullet$ is assumed to be connected.
\end{Con}

\begin{Def}[{\cite[\S~1]{LarsenPink95}}]\label{Def-GoodPlace}
A place $x\in |X|$ is called {\em good (for $\rho_\bullet$ (and $\lambda$))} if $\rho_\lambda(\Frob_x)$ is $\Gamma$-regular.
\end{Def}
The following result is proved in \cite{LarsenPink92} for $\BQ$-rational compatible systems. The adaptations to the general case are minor, and we omit them.
\begin{Prop}[{\cite[7.2, 7.4]{LarsenPink92}}]\label{Prop-GammaRegular}
For any $\lambda\in\CP'_E$, the set of all $g\in\pi_1(X)$ such that $\rho_\lambda(g)$ is $\Gamma$-regular is open and dense in $\pi_1(X)$. In particular, the set of good places $x\in|X|$ has \v{C}ebotarov density $1$. Moreover for $x\in |X|$ the condition of $\Frob_x$ being $\Gamma$-regular is independent of~$\lambda$.
\end{Prop}

Suppose that for every $\lambda\in\CP_E'$ we choose an $\CO_\lambda$-lattice $\Lambda_\lambda\subset E_\lambda^n$ that is invariant under~the~action of $\pi_1(X)$ via $\rho_\lambda$. Denote by $\CG_\lambda$ the Zariski closure of $G_\lambda$ in $\Aut_{\CO_\lambda}(\Lambda_\lambda)$, endowed with the unique structure of reduced closed subscheme. This is a flat group scheme over $\CO_\lambda$. We have the following immediate generalization of a result of Larsen and Pink:
\begin{Prop}[{\cite[Prop.~1.3]{LarsenPink95}}]\label{Prop-LP-GisSmooth}
For all $\ell_\lambda\gg0$, the following hold: The group $\CG_\lambda$ is smooth. It contains a closed subtorus $\CT_\lambda$ defined over $\CO_\lambda$ whose generic fiber is a maximal torus $T_\lambda$ of $G_\lambda$, and whose special fiber $\CT_{k_\lambda}$ is a maximal torus of $\CG_{k_\lambda}=\CG_\lambda\times_{\CO_\lambda}k_\lambda$. Moreover 
\begin{enumerate}
\item the set of roots (and their multiplicities) of $\CG_{k_\lambda}$ with respect to $\CT_{k_\lambda}$ is the same as the set of roots of $G_\lambda$ with respect to $T_\lambda$, and 
\item the formal character for the action of $T_\lambda$ on $E_\lambda^n$ and for the action of $\CT_{k_\lambda}$ on $k_\lambda^n$ agree.
\end{enumerate}
\end{Prop}
\begin{proof}
There are some explanations in order. The proof in \cite{LarsenPink95} for $\BQ$-compatible systems consists of two parts. First a nice model $\CT_\lambda$ (in the sense of \cite{CHT}) in $\CG_\lambda$ of a maximal torus of $G_\lambda$ is constructed. This argument uses $\Gamma$-regular elements recalled earlier in this subsection, and it works for $E$-compatible systems with any number field $E$; see also \cite[9.2.1]{CHT}. In a second step, for the proof of smoothness, an $\CO_\lambda$-scheme map $\phi_\lambda$ from $\CX_\lambda:=\CT_\lambda\times_{\CO_\lambda}\prod_\alpha\CU_\alpha$ for models $\CU_\alpha$ of the root spaces of $G_\lambda$ to $\CG_\lambda$ is considered. According to \cite[footnote, p.~3]{Wintenberger} the argument in \cite{LarsenPink95} for $\phi_\lambda$ being an isomorphism near zero and $\CX_\lambda$ being $\CO_\lambda$-smooth seems to be incomplete. This is remedied in \cite[Thme.~1]{Wintenberger} for $G_\lambda$ reductive, and in \cite[9.1.1]{CHT} for $G_\lambda$ semisimple. Both references rely on \cite[\S\S~1.2 and 2.2]{Bruhat-Tits-SurUnCorpsLocal2}. Note also that the reduction of $\phi_\lambda$ proves assertion (a), while (b) follows from having constructed $\CT_\lambda$ via $\Gamma$-regular elements.

We also remark that \cite[\S~1]{LarsenPink95} assumes purity as part of the definition of a compatible system. However this hypothesis is not used in the proof of \cite[Prop.~1.3]{LarsenPink95}, nor in the preliminaries leading up to it, e.g not in \cite{LarsenPink92}, and it is also not used in the our supplementary references.
\end{proof}

\smallskip

We now come to the central result of this subsection:
\begin{Thm}\label{Thm-RedIsSaturated}
For all $\ell_\lambda\!\gg0$, the special fiber $\CG_{k_\lambda}$ is smooth, and its identity component is saturated in~$\GL_{n,k}$.
\end{Thm}
\begin{proof} 
By \autoref{Prop-LP-GisSmooth} the group scheme $\CG_\lambda$ is smooth over $\CO_\lambda$ for $\ell_\lambda\gg0$, and hence $\CG_{k_\lambda}$ is smooth. The theorem now follows from \autoref{Thm-GisSaturated}, once we verified hypotheses (i)--(iii) listed at the beginning of \autoref{Subsec-NonReduSat} for $\CG_{k_\lambda}^o\subset\GL_{n,k_\lambda}$.

For this let $R'$ be the set of weights of $\CT_{k_\lambda}$ acting on $\bar V_\lambda:=k_\lambda^n$ and let $\Phi$ be the root system of $\CG_{k_\lambda}$. By \autoref{Prop-LP-GisSmooth}, the set $\Phi$ is also the root system of the reductive group $G_\lambda$, and this shows that conditions (i) and (ii) hold. 

For (iii) note, again by \autoref{Prop-LP-GisSmooth}, that the set $R'$ is also the set of weights of $T_\lambda$ for its action on $E_\lambda^n$. The sets $\Phi$ and $R'$ are the sets of roots of the split motivic group and of the split motivic representation attached to $\rho_\bullet$ over some finite extension $F/E$, respectively, and hence by \autoref{Cor-GenChinTriple} these sets are independent of $\lambda$. Thus also $R=\{\alpha-\beta\mid\alpha,\beta\in R'\}$ is independent of $\lambda$, and it is the set of weights for the action of $T_{k_\lambda}$ on $\End(\bar V_\lambda)$. It is now clear that condition (iii) must hold for $\ell_\lambda\gg0$.
\end{proof}

\section{Residual saturation in the absolutely irreducible case}
\label{Sec-ResidualSaturation}

In this section we shall prove two main results on $E$-rational compatible systems $\rho_\bullet$ that are absolutely irreducible and connected, \autoref{Thm-OnIrred} and \autoref{Thm-AI-Reduction-IsSaturated}, and in \autoref{Thm-MainGenlCase} an important result on lifting (arbitrary) residually compatible systems. The first theorem shows that their reductions $\bar\rho_\lambda$ are absolutely irreducible for almost all $\lambda$, even after restriction to $\pi_1^\geo(X)$. For the second, suppose each $\rho_\lambda$ takes its image in $\GL_n(\CO_\lambda)$ and let $\CG^\geo_\lambda$ denote the Zariski closure of $\pi_1^\geo(X)$ therein. Then we show that for almost all $\lambda$ the $\CO_\lambda$-group scheme $\CG^\geo_\lambda$ is semisimple connected, and its special fiber is the saturation of $\bar\rho_\lambda(\pi_1^\geo(X))$. Our proofs combine L.\ Lafforgue's global Langlands conjecture, and a conjecture of de Jong, as proved by Gaitsgory. Whenever the image is smaller than expected, we construct a congruence between automorphic forms, to a given automorphic form. Since we manage to bound the candidates for such a congruence within a finite set, only finitely many congruences are allowed. The first author had learned the usefulness of similar congruence arguments from Ribet, e.g. \cite[Prop.~5.3]{Weston}. In fact, by our very approach, after some preliminary observations in \autoref{Subsec-CharAndLifting}, we shall in \autoref{SubSect-ModLambdaSystems} first prove that any (suitably defined) potentially tame residually compatible system is the reduction of a compatible system. From this we shall quickly deduce \autoref{Thm-OnIrred} in \autoref{Subsect-ModLambdaAI}. \autoref{Subsec-ModLmabdaSaturated} contains the proof of \autoref{Thm-AI-Reduction-IsSaturated}. A proof of \autoref{Thm-OnIrred} was outlined by us in \cite{Boeckle-DFG}. However we recently learned that in \cite[Prop.~E.10.1]{Drinfeld-ProSemisimple} Drinfeld had given such a proof by essentially the same method. Since we give more details and also want to stress the assertion of \autoref{Thm-MainGenlCase} we do include proofs.

\subsection{On central characters and lifting mod $\lambda$ representations}
\label{Subsec-CharAndLifting}

This subsection contains some auxiliary results. We begin by defining two natural reductions of a $\lambda$-adic representation $\rho_{\lambda}$ which is a member of an $E$-rational compatible system $\rho_\bullet$. Then we provide some results on $1$-dimensional compatible systems that will be relevant to the lifting techniques used in the later subsections, and we also provide two basic results on lifting.
\medskip

Let $\Lambda_\lambda$ be an $\CO_\lambda$-lattice in $E_\lambda^n$ that is stable under the action of $\pi_1(X)$ via $\rho_\lambda$; the existence of such a lattice can be deduced from $\pi_1(X)$ being a profinite group. Then the image of $\rho_\lambda$ lies in $\Aut_{\CO_\lambda}(\Lambda_\lambda)$, and the latter group has a reduction map to $\Aut_{\BF_\lambda}(\Lambda_\lambda/\varpi_\lambda\Lambda_\lambda)\cong\GL_n(k_\lambda)$. We define $\bar\rho_\lambda^{\Lambda}:=\bar\rho_\lambda^{\Lambda_\lambda}$ as the composite of $\rho_\lambda$ with the just defined reduction map, i.e., as the action of $\pi_1(X)$ on $\Lambda_\lambda/\varpi_\lambda\Lambda_\lambda$. The semisimplification of $\bar\rho_\lambda^{\Lambda}$ we simply denote by $\bar\rho_\lambda$. 

The theorem of Brauer and Nesbitt (see \cite[30.16]{CurtisReinerOld}) shows that $\bar\rho_\lambda$ is independent of any choices. In general $\bar\rho_\lambda^{\Lambda}$ depends on $\Lambda_\lambda$. If $\bar\rho_\lambda$ is irreducible, then $\bar\rho_\lambda\cong\bar\rho_\lambda^{\Lambda}$, and so then $\bar\rho_\lambda^{\Lambda}$ is independent of the choice of $\Lambda_\lambda$. In this section we mainly deal with $\bar\rho_\lambda$, in \autoref{Sec-IndepOfLat} with~$\bar\rho_\lambda^{\Lambda}$.

\begin{Rem}\label{Rem-RedAndCoeffChange}
Let $E'$ be a finite extension of $E$ and let $\rho'_\bullet$ be the compatible system $\rho_\bullet\otimes_EE'$. Let $\lambda'\in\CP_{E'}'$ have contraction $\lambda$ to $\CP_E'$ and define $\Lambda'_{\lambda'}:=\Lambda_\lambda\otimes_{\CO_\lambda}\CO'_{\lambda'}$. Then one has a functorial isomorphisms for the respective reductions
\[(\bar\rho')^{\Lambda'}_{\lambda'}\cong \bar\rho^\Lambda_\lambda\otimes_{k_\lambda}k_{\lambda'}\hbox{ \ \ and \ \ } \bar\rho'_{\lambda'}\cong \bar\rho_\lambda\otimes_{k_\lambda}k_{\lambda'},\]
where $k_\lambda\to k_{\lambda'}$ is the map of the residue field induced from $E\to E'$ at $\lambda'$. Both isomorphisms are straightforward. The second uses (again) the theorem of Brauer and Nesbitt and the semisimplicity of the reductions. It follows  that if $\bar\rho_\lambda$ is absolutely irreducible, then so is $\rho_\lambda$. 
\end{Rem}

We begin with some preparations for \autoref{Thm-MainGenlCase} in the following subsection.
\begin{Lem}\label{Lem-OnCentralChars1}
Let $\alpha_1,\ldots,\alpha_n$ be plain of characteristic $p$, and let $x$ be in $|X|$. For any partition $\wp=(\wp_1,\ldots,\wp_t)$ of $n$ and $j\in \{1,\ldots,t\}$ define $\alpha_{\wp,j}$ as a $[\kappa_x:\kappa]$-th root of $\prod_{i\in \wp_j}\alpha_i$ and define $\rho_{\wp,j,\bullet}$ as the $1$-dimensional compatible system $\rho_{\alpha_{\wp,j},\bullet}$. Let $X'\to X$ be a finite cover, and let $\Xi$ be the set of continuous characters $\tau\colon\pi_1(X)\to \GL_1(\overline E)$ that are trivial when restricted to $\pi_1(x)$ and factor via $\pi_1(X')^{\ab,p}$ when restricted to $\pi_1(X')$. Then 
\begin{enumerate}
\item each $\alpha_{\wp,j}$ is plain of characteristic $p$, and hence $\rho_{\wp,j,\bullet}$ is well-defined; 
\item the set $\Xi$ is finite;
\item each representation $\bar\delta_\lambda\colon \pi_1(X)\to \GL_1(k'_\lambda)$ such that $\bar\delta_\lambda|_{\pi_1(X')}$ factors via $\pi_1^{\ab,p}(X')$ and $\bar\delta_\lambda(\Frob_x)=\prod_{i\in \wp_j}(\alpha_i \!\!\pmod\lambda)$ is congruent modulo $\lambda$ to $\rho_{\tau,\lambda}\otimes\rho_{\wp,j,\lambda}$ for some $\tau\in \Xi$.
\item if all $\alpha_i$ are $q$-Weil numbers of weight $w$, then all $\alpha_{\wp,j}$ are $q$-Weil numbers of weight~$w\cdot\#\wp_j$.
\end{enumerate}
\end{Lem}
\begin{proof} For (a) observe that plainness of characteristic $p$ is preserved under the formation of products and roots. Part (b) follows from \autoref{Thm-KatzLang}. Regarding (c) note that the $1$-dimensional representation $\bar\delta_\lambda\otimes(\rho_{\wp,j,\lambda}^{-1}\!\!\pmod \lambda)$ of $\pi_1(X)$ is trivial when restricted to $\pi_1(x)$ and factors via $\pi_1(X')^{\ab,p}$, and thus it is of the form $\rho_{\tau,\lambda}\!\!\pmod \lambda$ for some $\tau\in \Xi$; note that any representation $\pi_1(X)\to\GL_1(k'_\lambda)$ that is trivial on $\pi_1(x)$ and factors via $\pi_1(X')^{\ab,p}$, has Teichm\"uller lift in~$\Xi$. Part (d) is clear from the definition of $\alpha_{\wp,j}$ using that $x$ is of degree $[\kappa_x:\kappa]$ over~$\kappa$. 
\end{proof}

\begin{Lem}\label{Lem-OnCentralChars2}
Suppose that $\rho_\bullet$ is $n$-dimensional. Fix $x\in |X|$ and denote by $\alpha_1,\ldots,\alpha_n$ the roots of $P_x(T)\in E[T]$ in a fixed algebraic closure $\overline E$ of $E$. For any partition $\wp=(\wp_1,\ldots,\wp_t)$ of $n$ and $j\in \{1,\ldots,t\}$ define $\alpha_{\wp,j}$ as in \autoref{Lem-OnCentralChars1}. For each $\lambda\in\CP_{E}'$ and embedding $k_\lambda\to \overline k_\lambda$ choose an isomorphism
\[\bar\rho_\lambda\otimes_{k_\lambda}\overline k_\lambda\cong \bigoplus_{j=1}^t \bar\rho_{\lambda,j}\]
with each $\bar\rho_{\lambda,j}$ irreducible. Then there exists a finite set $\Xi$ of characters $\tau\colon\pi_1(X)\to \GL_1(\overline E)$ such that for each $\lambda\in\CP_{E}'$ there exists a partition $\wp=(\wp_1,\ldots,\wp_t)$ of $n$ and a character tuple $\utau=(\tau_1,\ldots,\tau_t)$ in $\Xi^t$, with $(\wp,\utau)$ depending on $\lambda$, such that for all $j=1,\ldots,t$ we have
\[ \det\bar\rho_{\lambda,j} = \bar\rho_{\tau_j,\lambda}\otimes\bar\rho_{\wp,j,\lambda}. \]
\end{Lem}
\begin{proof}
Fix a finite cover $X'\to X$ such that for all $\lambda\in\CP_E'$ the ramification of $\rho_\lambda|_{\pi_1(X')}$ is pro-$\ell_\lambda$, cf.~\autoref{Lem-OnTamingCover}, and define $\Xi$ as in \autoref{Lem-OnCentralChars1} for this cover, upon noting that $\alpha_1,\ldots,\alpha_n$ are plain of characteristic~$p$ by \autoref{Rem-OnNotation}. Now the reduction $P_x(T)\mod \lambda$ is the product of the characteristic polynomials of the $\bar\rho_{\lambda,j}(\Frob_x)$, $j=1,\ldots,t$. Hence there is a partition $\wp=(\wp_1,\ldots,\wp_t)$ of $n$ such that $\charpol_{\bar\rho_{\lambda,j}(\Frob_x)}=\prod_{i\in \wp_j}(T-\alpha_i\!\!\!\pmod{\lambda})$, and thus 
\[\det\bar\rho_{\lambda,j}(\Frob_x)=\alpha_{\wp,j}^{[\kappa_x:\kappa]}\!\!\!\pmod\lambda=\prod_{i\in \wp_j}\alpha_i\!\!\!\pmod\lambda.\]
The proof of the lemma is now implied by \autoref{Lem-OnCentralChars1}(c).
\end{proof}

The following result is based on de Jong's conjecture from \cite{deJong} and its solution by Gaitsgory for $\ell>2$ in \cite{Gaitsgory}, and of de Jong for $n\le 2$. We recall it in the form needed below.
\begin{Lem}[de Jong and Gaitsgory]\label{Lem-UDRingGLn}
Let $X$ be a smooth projective geometrically irreducible curve over $\kappa$ and let $S\subset|X|$ be a finite subset. Let $k$ be a field of characteristic $\ell>n$ with ring of Witt vectors $W(k)$ and $K=W(k)[1/\ell]$. Let $\bar\rho\colon \pi_1(X\setminus S)\to \GL_n(k)$ be a continuous absolutely irreducible representation, and let $\eta\colon \pi_1(X\setminus S)\to\GL_1(\CO_L)$ be a character for some finite extension $L$ of $K$ such that the reduction of $\eta$ agrees with~$\det\bar\rho$.

Then there exists a continuous representation $\rho\colon\pi_1(X)\to \GL_n(\CO_{L'})$ over a finite extension $L'$ of $L$ with $\det\rho=\eta$ and reduction isomorphic to $\bar\rho$. If moreover $\bar\rho$ is at most tamely ramified (along $S$), then the conductor of $\rho$ is bounded by $nS$.
\end{Lem}
\begin{proof} 
Denote by $\CR$ the universal ring for deformations $\rho\colon \pi_1(X)\to\GL_n(R)$ of $\bar\rho$ with $\det\rho=\eta$, where $R$ ranges over all complete noetherian local $\CO_L$-algebras with maximal ideal $\Fm_{R}$ and residue field $k=\CO_L/\Fm_{\CO_L}$; see~\cite{Mazur} and~\cite[3.1-3.3]{deJong} for the precise definition of the deformation problem and the existence of a universal deformation $\wt\rho\colon\pi_1(X)\to\GL_n(\CR)$ using the absolute irreducibility of $\bar\rho$. 

Now de Jong's conjecture \cite[Conj.~2.3]{deJong}  for $\ell>2$ is known by work of Gaitsgory; see \cite[proof of Conj.~1.3]{Gaitsgory}. Therefore by \cite[Thm.~3.5]{deJong}, the resulting homomorphism $\CO_L\to\CR$ is a finite flat complete intersection morphism. In particular $\CR\otimes_{\CO_L}L$ is a finite non-empty product of finite extensions $L'$ of $L$. Fix one $L'$. The corresponding map $p\colon\CR\to L'$ takes its image in the ring of integers $\CO_{L'}$ of $L'$, and defining the representation $\rho:=p\circ \wt\rho\colon \pi_1(X)\to\GL_n(\CO_{L'})$ proves the first assertion of the lemma.

Regarding the second assertion, observe that if $\bar\rho$ is at most tamely ramified, the same holds for $\rho$, since the kernel of $\rho(\pi_1(X))\to\bar\rho(\pi_1(X))$ is a pro-$\ell$ group. But the conductor of a tame representation, at any place of $S$, of an $n$-dimensional representation is bounded by~$n$.
\end{proof}

For later use, we also state the following variant of \autoref{Lem-UDRingGLn}.
\begin{Lem}[{\cite[Thm.~5.14]{BHKT}}]\label{Lem-UDRingSemisimple}
Let $X$ be a smooth proper geometrically connected curve over $\kappa$, and $S\subset |X|$ be a finite set. Let $k$ be a finite field of characteristic $\ell\neq p$ with $\ell > 2(n-1)$. Let $\CG$ be semisimple connected and split affine algebraic group over $W(k)$ that has an almost faithful representation $i \colon \CG\times_{W(k)}k\to \GL_{n,k}$. Let $\bar\rho\colon \pi_1(X\setminus S)\to \CG(k)$ be a representation such that $i\circ \bar\rho$ is absolutely irreducible. Then there exists a finite extension $E$ of $W(k)[1/\ell]$ and a continuous homomorphism $\rho : \pi_1(X\setminus S) \to \CG(\CO_E)$ such that $\rho \pmod {\Fm_{\CO_E}} = \bar{\rho}$.
\end{Lem}

\subsection{Potentially tame mod $\lambda$ compatible systems lift}

\label{SubSect-ModLambdaSystems}

\begin{Def}\label{Def-ModLambdaCS}
An $n$-dimensional {\em $E$-rational mod $\lambda$ compatible system over $\CL\subset\CP_E'$} is a tuple \[((\bar\rho_\lambda)_{\lambda\in\CL},(P_x)_{x\in|X|}),\] or for short $\bar\rho_\bullet$, that consists of a homomorphism $\bar\rho_\lambda\colon \pi_1(X)\to\GL_n(k_\lambda)$ for every $\lambda\in\CL$, and a monic degree $n$ polynomial $P_x\in E[T]$ for each $x\in|X|$, such that
\begin{enumerate}
\item $\bar\rho_\lambda$ is semisimple for all $\lambda\in\CL$,
\item $P_x$ is plain of characteristic $p$ for all $x\in |X|$, and
\item $\charpol_{\bar\rho_\lambda(\Frob_x)}(T)=P_x(T)\pmod \lambda$ for all $(\lambda,x)\in \CL\times |X|$
\end{enumerate}
The system is called {\em pure of weight $w$}, if $P_x$ is pure of weight $w$ for $\#\kappa_x$ for all $x\in |X|$.

The notions {\em tame} and {\em potentially tame} are defined as for compatible systems, using $\bar\rho_\lambda$ in place of~$\rho_\lambda$. 

For $E\to E'$ a finite extension, we define the mod $\lambda$ compatible system $\bar\rho_\bullet\otimes_{\CO_E}\CO_{E'}$ over $\CL'$, the subset of $\CP_{E'}'$ of places above $\CL$, by setting $\bar\rho_\mu:=\bar\rho_\lambda\otimes_{k_\lambda}k_\mu$ for $\mu\in\CL'$ and $\lambda\in\CL$ its contraction, and by considering $P_x\in E'[T]$ via $E\to E'$; 
the $\bar\rho_\mu$ are again semisimple by \cite[69.8]{CurtisReinerOld}.
\begin{Rem}\label{Rem-ModLambdaCompatSyst-BasicFacts}
Note the following:
\begin{enumerate}
\item If $\rho_\bullet$ is an $n$-dimensional $E$-rational compatible system, then $(\bar\rho_\lambda)_{\lambda\in\CP_E'}$, defined as at the beginning of \autoref{Subsec-CharAndLifting}, is an $n$-dimensional $E$-rational mod $\lambda$ compatible system over $\CP_E'$. We then call $\bar\rho_\bullet$ the reduction of $\rho_\bullet$.  If $\rho_\bullet$ is pure, tame, potentially tame, then the same property holds for $\bar\rho_\bullet$
\item Given an $E$-rational mod $\lambda$ compatible system over some infinite subset $\CL\subset \CP_E'$, the polynomials $(P_x)_{x\in |X|}$ are uniquely determined by the representations $(\bar\rho_\lambda)_{\lambda\in\CL}$ and condition~(c). In particular, by \autoref{Lem-IsomOfCS}, there exists up to isomorphism at most one $E$-rational compatible system $\rho_\bullet$ whose reduction restricted to $\CL$ agrees with $\bar\rho_\bullet$.
\end{enumerate}
\end{Rem}
\end{Def}
\begin{Thm}\label{Thm-MainGenlCase}
Suppose $X$ is a geometrically connected curve over $\kappa$ and $\bar\rho_\bullet$ is a potentially tame $E$-rational $n$-dimensional mod $\lambda$ compatible system over an infinite set $\CL$. Then there exists a finite extension $E'$ of $E$ and an $E'$-rational compatible system $\rho_\bullet$ (over $\CP_{E'}'$) whose semisimplified reduction is $\bar\rho_\bullet\otimes_{\CO_E}\CO_{E'}$ over $\CL':=\{\lambda'\in\CP_{E'}'\mid \lambda' \hbox{ lies above some }\lambda\in\CL\}$.

If $\bar\rho_\bullet $ is absolutely irreducible for some $\lambda_0$, then $\rho_\bullet$ is absolutely irreducible. If $\bar\rho_\lambda\otimes_{k_\lambda}\bar k_\lambda$ is reducible for infinitely many $\lambda\in \CL$, then $\rho_\bullet$ is reducible.
\end{Thm}
\begin{proof}
Let $X'\to X$ be a finite Galois cover over which all $\bar\rho_\bullet$ are tame. Fix $x\in |X|$ and denote by $\alpha_1,\ldots,\alpha_n$ the roots of $P_x(T)\in E[T]$ in a fixed algebraic closure $\overline E$ of $E$. They are plain of characteristic $p$. For any partition $\wp=(\wp_1,\ldots,\wp_t)$ of $n$ and $j\in \{1,\ldots,t\}$ define $\alpha_{\wp,j}$ as in \autoref{Lem-OnCentralChars1}. 

For each $\lambda\in\CP_{E}'$ choose a finite extension $k_\lambda\to k'_\lambda$ over which there is an isomorphism
\[\bar\rho_\lambda\otimes_{k_\lambda} k'_\lambda\cong \bigoplus_{j=1}^{t_\lambda} \bar\rho_{\lambda,j}\]
with each $\bar\rho_{\lambda,j}$ absolutely irreducible. Then following the proof of \autoref{Lem-OnCentralChars2} there exists a finite set $\Xi$ of characters $\tau\colon\pi_1(X)\to \GL_1(\overline E)$ such that for each 
$\lambda\in\CP_{E}'$ there exists a partition $\wp_\lambda=(\wp_{\lambda,1},\ldots,\wp_{\lambda, t_\lambda})$ of $n$ and a character tuple $\utau_\lambda=(\tau_{\lambda,1},\ldots,\tau_{\lambda,t_\lambda})$ in $\Xi^{t_\lambda}$, such that for all $j=1,\ldots,t_\lambda$ we have
\[ \det\bar\rho_{\lambda,j} = \bar\rho_{\tau_{\lambda,j},\lambda}\otimes\bar\rho_{\alpha_{\wp_\lambda,j}}. \]

By \autoref{Lem-UDRingGLn}, we obtain a finite extension $L$ of $E_\lambda$ and for each $j\in\{1,\ldots,t_\lambda\}$ a representation
\[\rho_{\lambda,j}\colon\pi_1(X)\to \GL_{n_{\lambda,j}}(\CO_L)\]
that is a lift of $\bar\rho_{\lambda,j}$ with $\det\rho_{\lambda,j}=\rho_{\tau_{\lambda,j},\lambda}\otimes\rho_{\alpha_{{\wp_\lambda},j},\lambda}$. Note also that $\rho_{\lambda,j}$ is absolutely irreducible, because $\bar\rho_{\lambda,j}$ is so. 

Let $S$ be the set of places of $\kappa(X)$ not in $|X|$, and let $S'$ be the places of $\kappa(X')$ above $S$. Let $K_s$ be the local field at $s$, $K'_{s'}$ be the extension corresponding to $s'\in S'$ above $s$. Then any $\rho_{\lambda,j}$ restricted to $\Gamma_{K'_{s'}}$ is tame. Considering the higher ramification filtration, it follows that there is a uniform bound $M$ on the wild conductors at any $s\in S$ of any $\rho_{\lambda,j}$, supported on $S$. The tame conductor of $\rho_{\lambda,j}$ is bounded by $\wp_{\lambda,j}S$. By \autoref{Thm-Lafforgue} there is an automorphic representation $\Pi_{\lambda,j}$ of $\GL_{n/\BA_{\kappa(X)}}$ with central character $\tau_j$ and conductor bounded by $n_{\lambda,j}S+M$ with $\rho_{\Pi_{\lambda,j},\lambda}\cong \rho_{\lambda,j}\otimes\rho_{\wp,j,\lambda}^{-1}$. Define $\CS$ as the set of tuples $(\wp,\utau,\Pi_1,\ldots,\Pi_t)$ where $\wp=(\wp_1,\ldots,\wp_t)$ is a partition of $n$, each $\utau$ lies in $\Xi^t$, each $\Pi_j$ is a cuspidal automorphic representation for $\GL_{\#\wp_j}$ over $\BA_{\kappa(X)}$ of conductor bounded by $\wp_j S+M$ and central character that agrees with $\tau_j$. Now by \cite[Cor.~5.3]{Borel-Jacquet} the set $\CS$ is finite. 

From the finiteness of $\CS$ it follows that there exists a tuple $(\wp,\utau,\Pi_1,\ldots,\Pi_t)$ in $\CS$ and an infinite subset $\CL_0$ of $\CL$ such that for all $\lambda\in\CL_0$ we have $\rho_{\lambda,j}\cong \rho_{\Pi_j,\lambda}\otimes\rho_{\wp,j,\lambda}$. Let $E'$ be a finite extension of $E$ over which the compatible system $\rho'_\bullet:=\oplus_{j=1}^t\rho_{\Pi_j,\bullet}\otimes\rho_{\wp,j,\bullet}$ is defined. By the choice of $\CL_0$ we have for all $\lambda\in\CL_0$ an isomorphism $\rho'_\mu\!\pmod \mu \cong \bar\rho_\lambda\otimes_{k_\lambda}k_\mu$ for any $\mu\in\CP_{E'}'$ above $\lambda$. This implies that the polynomials $P_x$ for $\bar\rho_\bullet$ and for $\rho'_\bullet$ agree in $E'[T]$ for all $x\in|X|$. But then it follows that 
\[\charpol_{\rho'_\mu(\Frob_x)} = \charpol_{\bar\rho_\lambda(\Frob_x)}\otimes_{k_\lambda}k_\mu \]
for all $x\in|X|$ and $\mu\in\CL'$. The \v{C}ebotarov density theorem now implies that $\bar\rho_\bullet\otimes_{\CO_E}\CO_{E'}$ is the reduction of $\rho_\bullet$, and this concludes the proof of the main part. 

The first part of the final assertion follows from \autoref{Rem-RedAndCoeffChange} and \autoref{Lem-OnAbsIrredAndTwists}. To see the second part, let $\CL^\red$ be the set of places in $\CL$ at which $\bar\rho_\lambda\otimes_{k_\lambda}\bar k_\lambda$ is reducible, and let $\CS^\red$ be the set of tuples $(\wp,\utau,\Pi_1,\ldots,\Pi_t)$ in $\CS$  such that $t>1$. Then each $\bar\rho_\lambda$, $\lambda\in\CL^\red$ arises via reduction from some tuple in $\CS^\red$. Repeating the argument from the previous paragraph, we find that $(\bar\rho_\lambda)_{\lambda\in\CL^\red}$ is the reduction of a reducible compatible system $\rho_\bullet'$ over some~$E'$, where in fact $\rho_\bullet'$ is defined for all $\lambda\in\CP_{E'}'$. By the uniqueness of the compatible system for $\CL$ and for $\CL^\red$, we deduce that $\rho_\bullet$ constructed for $\CL$ is reducible over~$E'$.
\end{proof}

\begin{Rem}\label{Rem-DrinfeldCurveToVariety}
If $X$ is not a curve, the above result can be applied to the restriction of the system to any curve $C\into X$. Following \cite[\S~4]{Drinfeld12} this allows one to extend \autoref{Thm-MainGenlCase} to arbitrary normal, geometrically irreducible finite type varieties $X$ over $\kappa$ under a potential tameness hypothesis.
\end{Rem}
\begin{Ques}
We do not know if the hypothesis of potential tameness is needed in \autoref{Thm-MainGenlCase}. Can one construct counterexamples if $\CL$ is infinite, or even if $\CL$ has density one? 
\end{Ques}

\subsection{Mod $\lambda$ absolute irreducibility for almost all $\lambda$.}
\label{Subsect-ModLambdaAI}
      
\begin{Thm}[{cf.~\cite[Prop.~E.10.1]{Drinfeld-ProSemisimple}}]\label{Thm-OnIrred}
Suppose $\rho_\bullet$ is absolutely irreducible and $E$-rational. Then for almost all $\lambda\in\CP_E'$ the representation $\bar\rho_\lambda$ is absolutely irreducible.
\end{Thm}
\begin{proof}
By \autoref{Lem-OnAbsIrredAndTwists}(c) and \autoref{Lem-IsomOfCS}, there is a smooth geometrically irreducible curve $C\subset X_\reg$ such that $\rho_\bullet|_{\pi_1(C)}$ is absolutely irreducible. If now $\bar\rho_\mu|_{\pi_1(C)}$ is absolutely irreducible for almost all $\mu$, then the same will hold for $\bar\rho_\mu$; hence from now on we assume that $X$ is a smooth and geometrically irreducible curve. It follows from \autoref{Lem-OnTamingCover} that $\rho_\bullet$ is potentially~tame.

Let $\CL$ be the set of places of $\CP_E'$ at which $\bar\rho_\lambda\otimes_{k_\lambda}\bar k_\lambda$ is reducible. Assume that $\CL$ is infinite. Then by \autoref{Thm-MainGenlCase}, there exists a finite extension $E'$ of $E$ and an $E'$-rational compatible system $\rho'_\bullet$ (over $\CP_{E'}'$) whose semisimplified reduction is $\bar\rho_\bullet\otimes_{\CO_E}\CO_{E'}$ over $\CL':=\{\lambda'\in\CP_{E'}'\mid \lambda' \hbox{ lies above some }\lambda\in\CL\}$, and moreover this system is reducible. Now by \autoref{Rem-ModLambdaCompatSyst-BasicFacts}(b) we must have $\rho_\bullet\otimes_EE'\cong\rho'_\bullet$, which is a contradiction.
\end{proof}

\begin{Cor}\label{Cor-OnGeomIrred}
Suppose $\rho_\bullet$ is absolutely irreducible, connected and $E$-rational. Then for almost all $\lambda\in\CP_E'$ the representation $\bar\rho_\lambda|_{\pi_1^\geo(X)}$ is absolutely irreducible.
\end{Cor}
\begin{proof}
Let $\CL$ be the set of places $\lambda$ of $\CP_E'$ with $\ell_\lambda>n$ for which $\bar\rho_\lambda$ is absolutely irreducible, but $\bar\rho_\lambda|_{\pi_1^\geo(X)}$ is not. Let $\lambda$ be in $\CL$. The short exact sequence $1\to \pi_1^\geo(X)\to \pi_1(X)\to \Gamma_\kappa\to 1$ induces a short exact sequence
\[ 1\to N_\lambda:=\bar\rho_\lambda(\pi_1^\geo(X)) \to H_\lambda:= \bar\rho_\lambda(\pi_1(X)) \to C_\lambda \to 1,\]
where $C_\lambda$ is defined to be the cokernel of the morphism on the left. Clearly $C_\lambda$ is cyclic of finite order. By hypothesis the action of $H_\lambda$ on $V_\lambda:=\overline k^n_\lambda$ defined by $\bar\rho_\lambda$ is irreducible. By Clifford theory, cf.~\cite[\S~49, 50]{CurtisReinerOld}, the module $V_\lambda$ is completely reducible for the action of $N_\lambda$, and if $W_\lambda$ is an irreducible $N_\lambda$-submodule of $V_\lambda$, and $N^*_\lambda\supseteq N_\lambda$ the stabilizer of $W_\lambda$ in $H_\lambda$, then $V_\lambda\cong \Ind_{N^*_\lambda}^{H_\lambda}W_\lambda$. It follows that the index $[H_\lambda:N_\lambda^*]$ is at most of size $n$. Let $\kappa'$ be the unique extension of $\kappa$ of degree $n!$ and $X'\to X$ the base change of $X$ along $\Spec\kappa'\to\Spec\kappa$. Then for all $\lambda\in\CL$, the restriction $\bar\rho_\lambda|_{\pi_1(X')}$ is not absolutely irreducible. Because $\rho_\bullet$ is absolutely irreducible and connected, the same holds for $\rho_\bullet|_{\pi_1(X')}$. Thus $\CL$ is finite by \autoref{Thm-OnIrred}.
\end{proof}

\begin{Cor}\label{Cor-SatIsSemisimple}
Suppose $\rho_\bullet$ is absolutely irreducible, connected and $E$-rational. Then for almost all $\lambda\in\CP_E'$ we have $(\bar\rho_\lambda(\pi_1(X))^\sat_{k_\lambda})^o=(\bar\rho_\lambda(\pi^\geo_1(X))^\sat_{k_\lambda})^o=(\bar\rho_\lambda(\pi^\geo_1(X))^+)^\sat_{k_\lambda}$ and this group is connected semisimple.
\end{Cor}
\begin{proof}
Let $G=\bar\rho_\lambda(\pi_1(X))$ and $N=\bar\rho_\lambda(\pi_1^\geo(X))$. Then $N$ and $N^+$ are normal subgroups of $G$, $G/N$ is finite cyclic, and $N/N^+$ of order prime to $\ell$. By first applying  \autoref{Thm-AlmostIndependence} and then \autoref{Cor-OnGeomIrred}, it follows that $N^+$ acts absolutely irreducibly on $k_\lambda^n$ for almost all $\lambda\in\CP_E'$. By \autoref{Prop-GConnSemisimple}(b) we deduce that $G^\sat_{k_\lambda}$, $N^\sat_{k_\lambda}$ and $(N^+)^\sat_{k_\lambda}$ have a reductive identity component for almost all $\lambda$. Now by \autoref{Prop-GeomArithSaturation} we have $(G^\sat_{k_\lambda})^o=(N^\sat_{k_\lambda})^o$, and because $N/N^+$ is finite of order prime to $\ell$, by \autoref{Cor-SatAndWtG} we have $((N^+)^\sat_{k_\lambda})^o=(N^\sat_{k_\lambda})^o$. Finally $((N^+)^\sat_{k_\lambda})$ is connected semisimple by \autoref{Prop-GConnSemisimple}(c).
\end{proof}

\subsection{Mod $\lambda$ saturatedness for almost all $\lambda$.}
\label{Subsec-ModLmabdaSaturated}

Let $\rho_\bullet$ be an $E$-rational compatible system. For each $\lambda\in\CP_E'$, let $\Lambda_\lambda\subset E_\lambda^n$ be a $\pi_1(X)$-stable $\CO_\lambda$ lattice. Define $\CG_\lambda$ (resp.\ $\CG_\lambda^\geo$) as the Zariski closure of $\rho_\lambda(\pi_1(X))$ (resp.\ $\rho_\lambda(\pi_1^\geo(X))$) in $\Aut_{\CO_\lambda}(\Lambda_\lambda)\cong\GL_n(\CO_\lambda)$. Define 
$\CG_{k_\lambda}:=\CG_\lambda\otimes_{\CO_\lambda}k_\lambda$ and $\CG^\geo_{k_\lambda}:=\CG^\geo_\lambda\otimes_{\CO_\lambda}k_\lambda$.

\begin{Thm}\label{Thm-AI-Reduction-IsSaturated}
Suppose that $\rho_\bullet$ is absolutely irreducible and connected, and that $\rho_\bullet|_{\pi_1^\geo(X)}$ is connected, as well.\footnote{See \autoref{Cor-Reduction-IsSaturated} for a formulation that does not require any connectedness hypothesis.} Then for almost all $\lambda\in\CP_E'$ the following hold:
\begin{enumerate}
\item $\CG^\geo_{k_\lambda}$ is connected, semisimple and saturated; thus $\CG^\geo_{\lambda}$ is connected semisimple over~$\CO_\lambda$.
\item We have $\bar\rho_\lambda(\pi_1^\geo(X))^\sat_{k_\lambda}=\CG^\geo_{k_\lambda}$.
\item If $G^\geo_{\lambda}=G_{\lambda}$, then $\bar\rho_\lambda(\pi_1^\geo(X))^\sat_{k_\lambda}=\CG^\geo_{k_\lambda}=\CG_{k_\lambda}=\bar\rho_\lambda(\pi_1(X))^\sat_{k_\lambda}$.
\item If $G^\geo_{\lambda}\subsetneq G_{\lambda}$, denoting by $\CZ_\lambda$ the central torus of $\GL_{n/\CO_\lambda}$, then we have
\begin{enumerate}
\item $\CG_\lambda=\CG_\lambda^\geo\CZ_\lambda$, the group $\CG_\lambda$ is connected reductive over $\CO_\lambda$, and
\item $\big(\bar\rho_\lambda(\pi_1(X))^\sat_{k_\lambda}\big)^o =\CG^\geo_{k_\lambda}$. 
\end{enumerate}
\end{enumerate}
\end{Thm}
\begin{proof}
We first assert that $\CG_\lambda^\geo$ is smooth over $\CO_\lambda$ for almost all $\lambda\in\CP_E'$: By \autoref{Lem-OnAbsIrredAndTwists}(c), there exists a $1$-dimensional compatible system $\delta_\bullet$ such that $G_{\rho_\bullet\otimes\delta_\bullet,\lambda}$ has finite order determinant and at the same time
$G_{\rho_\bullet|_{\pi_1^\geo(X)},\lambda}=G_{\rho_\bullet\otimes\delta_\bullet|_{\pi_1^\geo(X)},\lambda}$. This implies $G_{\rho_\bullet|_{\pi_1^\geo(X)},\lambda}=G^o_{\rho_\bullet\otimes\delta_\bullet,\lambda}$ by \autoref{Prop-WeightZeroArithMonodromy}. Clearly there is a finite cover $X'\to X$ such that $\rho_\bullet\otimes\delta_\bullet|_{\pi_1(X')}$ and $\rho_\bullet\otimes\delta_\bullet|_{\pi_1^\geo(X')}$ are connected. Because $\rho_\bullet|_{\pi_1^\geo(X)}$ is connected, we have $G_{\rho_\bullet|_{\pi_1^\geo(X)},\lambda}=G_{\rho_\bullet\otimes\delta_\bullet|_{\pi_1(X')},\lambda}$. Now our assertion follows from \autoref{Prop-LP-GisSmooth}.

For (a) let $\CL$ be the set of places $\lambda\in\CP_E'$ such that $\bar\rho_\bullet|_{\pi_1^\geo(X)}$ is absolutely irreducible and $\CG^\geo_\lambda$ is smooth over $\CO_\lambda$. By \autoref{Cor-OnGeomIrred} and the above assertion, the set $\CP_E'\setminus\CL$ is finite. For $\lambda\in\CL$ we have $\bar\rho_\lambda^{\Lambda_\lambda}=\bar\rho_\lambda$ and it follows that $\bar\rho_\lambda(\pi_1^\geo(X))\subset\CG^\geo_{k_\lambda}(k_\lambda)$ by the definition of $\CG^\geo_\lambda$. Hence the action of $\CG^\geo_{k_\lambda}$ on $k_\lambda^n$ via its embedding into $\GL_{n,k_\lambda}$ is irreducible. By the proof of \cite[Thme.~1]{Wintenberger} the identity component of $\CG^\geo_{k_\lambda}$ is reductive -- the argument deduces from the semisimplicity of the action that the unipotent radical of $(\CG^\geo_{k_\lambda})^o$, which is a characteristic subgroup, must vanish. It follows from \cite[Prop.~3.1.3]{ConradSGA3} that $(\CG^\geo_\lambda)^o$ is a reductive group over $\CO_\lambda$ that is open and closed in $\CG^\geo_\lambda$. Since $\CG_\lambda^\geo$ was defined as the closure of the connected group $G^\geo_\lambda$, we deduce $\CG^\geo_\lambda=(\CG^\geo_\lambda)^o$ and hence $\CG^\geo_{\lambda}$ and $\CG^\geo_{k_\lambda}$ are connected reductive; they are semisimple because $G_\lambda^\geo$ is semisimple. From \autoref{Thm-RedIsSaturated} (and the reduction in the first paragraph), it follows that $\CG^\geo_{k_\lambda}$ is also saturated.

To prove (b) note first that by (a) we have $\bar\rho_\lambda(\pi_1^\geo(X))^\sat_{k_\lambda}\subseteq\CG^\geo_{k_\lambda}$ for $\lambda\in\CL$, since the group on the right is saturated. Let $X'\to X$ be a finite cover such that $\bar\rho_\lambda(\pi_1^\geo(X'))$ is $\ell_\lambda$-generated for almost all $\lambda$; it exists by \autoref{Rem-AIndep-AndEtoQ}. Note that this property is preserved under any further cover $X''\to X'$ since the index $[\bar\rho_\lambda(\pi_1^\geo(X')):\bar\rho_\lambda(\pi_1^\geo(X''))]$ is smaller than $\ell_\lambda$ for almost all $\lambda$. Set $\CH_{k_\lambda}:=\bar\rho_\lambda(\pi_1(X))^\sat_{k_\lambda}$ and $\CH^\geo_{k_\lambda}:=\bar\rho_\lambda(\pi_1^\geo(X))^\sat_{k_\lambda}$, and observe by  \autoref{Cor-SatIsSemisimple} that $\CH^\geo_{k_\lambda}$ is semisimple and connected for almost all $\lambda$. Enlarging the finite base field $\kappa$ to its unique extension of degree $n!$, and using \autoref{Prop-GeomArithSaturation}, we may assume from now on that $\CH^\geo_{k_\lambda}=\CH_{k_\lambda}$ for almost all $\lambda$. After another cover from \autoref{Lem-OnTamingCover}, and then by applying  \autoref{Cor-ReductionToCurve}, we may further assume that each $\rho_\lambda|_{\pi_1(X')}$ is at most $\ell_\lambda$-ramified, and that there is a smooth curve $C$ on $X'$ such that $\rho_\lambda(\pi_1(X'))=\rho'_\lambda(\pi_1(C))$ for  $\rho'_\lambda:=\rho_\lambda|_{\pi_1(C)}$ for all $\lambda$; in particular $\rho'_\bullet$ is tame. To complete the proof of (b) it will suffice to prove that $\bar\rho_\lambda(\pi_1^\geo(C))^\sat_{k_\lambda}=\CG^\geo_{k_\lambda}$ for almost all $\lambda\in\CL$. Note also that $\CG_\lambda$ and $\CG^\geo_\lambda$ are unchanged if we replace $X$ by $X'$ by the connectivity hypotheses of the theorem, and that $\CH_{k_\lambda}=\bar\rho_\lambda(\pi_1^\geo(C))^\sat_{k_\lambda}$, for almost all $\lambda$.

Let $\CL'\subseteq\CL$ be the set of places $\lambda$ where $\bar\rho'_\lambda(\pi_1^\geo(C))$ is $\ell_\lambda$-generated, $\CH^\geo_{k_\lambda}=\CH_{k_\lambda}$, and $\CG_\lambda$, $\CG^\geo_\lambda$ are smooth over $\CO_\lambda$, and $\ell_\lambda>2(n-1)$, and $\bar\rho'_\lambda$ is absolutely irreducible. By the results quoted so far $\CP_E'\setminus\CL'$ is finite. By the reduction argument from the first paragraph, to prove (b) we may also assume that $G_\lambda=G_\lambda^\geo$. Let $\CL''\subset\CL'$ be the set of places $\lambda$ at which $\CH_{k_\lambda}\subsetneq\CG_{k_\lambda}$, and fix $\lambda\in\CL''$. Let $k_\lambda\to k'_\lambda$ be a finite extension such that $\CH_{k'_\lambda}:= \CH_{k_\lambda}\otimes_{k_\lambda}k'_\lambda$ is split, connected and semisimple. Denote by $V_{k'_\lambda}$ the representation of $\CH_{k'_\lambda}$ via its embedding into $\SL_n$. Then from \autoref{Lem-TensorProdOfReps} and \autoref{Prop-LiftingOfLowWeightReps}. we obtain a semisimple connected group scheme $\CH_\lambda$ over $W(k'_\lambda)$ and a representation $\CV_\lambda$ of $\CH_\lambda$ defined over $W(k_\lambda')$ whose reductions to $k_\lambda'$ agree with $\CH_{k'_\lambda}$ and $V_{k'_\lambda}$. Since $W(k_\lambda')$ is regular local, we may identify $\CV_\lambda$ with $W(k_\lambda')^n$. Because $\CH_\lambda$ is semisimple connected, its image under the representation $\iota_\lambda\colon \CH_\lambda\to\Aut_{W(k)}(\CV_\lambda)$ must then lie in $\SL_{n,W(k'_\lambda)}$.

At this point, \autoref{Lem-UDRingSemisimple} provides a representation $\rho''_\lambda\colon \pi_1(C)\to\CH_\lambda(\CO'_\lambda)$ where $\CO'_\lambda$ is the valuation ring of some finite extension $E'_\lambda$ of $W(k_\lambda')[1/\ell_\lambda]$. The composition $\iota_\lambda\circ\rho''_\lambda\colon \pi_1(C)\to \SL_n(\CO'_{\lambda})$ is then a lift of $\bar\rho'_\lambda$. Since $\bar\rho'_\lambda$ is absolutely irreducible, so is the lift, if regarded as a representation into $\SL_n(E'_\lambda)$. Because $\det\bar\rho'_\lambda$ is trivial and $\ell_\lambda>n$, $\det \iota_\lambda\circ\rho''_\lambda$ is trivial. Note also that the monodromy group of $\iota_\lambda\circ\rho''_\lambda$ has smaller dimension than $G_\lambda$, because $\CH_{k_\lambda}\subsetneq\CG_{k_\lambda}$ and both are connected. Since $\bar\rho'_\lambda$ is $\ell_\lambda$-tame, this also holds for $\rho''_\lambda$, and hence the conductor of $\iota_\lambda\circ\rho''_\lambda$ is at most $n S$ where $S$ is the set of places of $\kappa(C)$ not in $|C|$, and so this conductor is independent of $\lambda$. Let $\CS$ be the set of cuspidal automorphic representations $\Pi$ for $\GL_n$ over $\BA_{\kappa(C)}$ of conductor bounded by $n S$ and with trivial central character. It follows from \autoref{Thm-Lafforgue} that there is an automorphic representation $\Pi_{\lambda}$ in $\CS$ with $\rho_{\Pi_{\lambda},\lambda}\cong \iota_\lambda\circ\rho''_\lambda$. Arguing as at the end of the proof of \autoref{Thm-OnIrred}, it follows that first $\CS$ and then $\CL''$ is finite, proving~(b).~Part~(c)~is immediate from (a), (b) and  $\CG^\geo_{k_\lambda}=\bar\rho_\lambda(\pi_1^\geo(X))^\sat_{k_\lambda}\subseteq \bar\rho_\lambda(\pi_1(X))^\sat_{k_\lambda}\subseteq \CG_{k_\lambda}=\CG^\geo_{k_\lambda}$ for $\lambda\in\CL\setminus\CL''$. 

We now turn to (d), and thus assume $G_\lambda^\geo\subsetneq G_\lambda$. To prove (i), we let $X'\to X$ and $\delta_\bullet$ be as in the first paragraph. Then $\rho_\bullet\otimes\delta_\bullet|_{\pi_1(X')}$ has trivial determinant and $\CG^\geo_{\rho_\bullet,\lambda}=\CG_{\rho_\bullet\otimes\delta_\bullet,\lambda}$. It follows that $\CG^\geo_{\lambda} \subseteq\CG_{\lambda} \subseteq \CZ_\lambda\CG_\lambda^\geo$. As $\det\rho_\bullet$ cannot have finite order in (d), $\det$ must surject $\CG_\lambda$ onto $\BG_m$, mapping the semisimple group $\CG^\geo_\lambda$ to $1$. Since $\CZ_\lambda$ is $1$-dimensional connected, it follows that $ \CG_\lambda\subseteq\CZ_\lambda\CG_\lambda^\geo$ must be an equality. Also, because $\CG_\lambda^\geo$ is connected semisimple over $\CO_\lambda$ and $\CZ_\lambda$ is a torus, the group $\CG_\lambda$ is connected reductive over $\CO_\lambda$, and completes the proof of~(i). Finally (ii) follows from (a) and (i), by which $\CG_{k_\lambda}$ is saturated and $\CG_\lambda(k_\lambda)/\CG_\lambda^\geo(k_\lambda)$ is finite of order prime to~$\ell_\lambda$.
\end{proof}

\begin{Rem}
\autoref{Thm-AI-Reduction-IsSaturated} for $\BQ$-rational compatible systems but only for a set of primes of density one and not all but finitely many primes is due to Larsen; cf.~\cite{Larsen95}. The result for $\BQ$-rational compatible systems that are cohomological in the sense of the introduction is proved in \cite{CHT}.
\end{Rem}

\begin{Rem}\label{Rem-AbsIrrThmFromCHToverQ}
Given \cite[Thm.~1.2]{CHT}, the following argument, whose details remain to be completed, gives a proof of \autoref{Thm-AI-Reduction-IsSaturated} for arbitrary $\BQ$-rational absolutely irreducible connected compatible systems $\rho_\bullet$ such that $\rho_\bullet$ and $\rho_\bullet|_{\pi_1^\geo(X)}$ are connected: Lafforgue's proof \cite{Lafforgue} of the Langlands correspondence for global function fields, partially quoted in \autoref{Thm-Lafforgue}, shows that $\rho_\bullet$ occurs in the \'etale cohomology of a (not necessarily smooth) scheme over the function field. Presumably by an argument like \cite[Thm.~7.3]{BGP-AdelicOpen} there is a finite number of smooth projective varieties over the function field such that $\rho_\bullet$ occurs in the direct sum of  their  \'etale $\ell$-adic cohomologies (though for each $\ell$ it is not clear in which summand). Now  \cite[Thm.~1.2]{CHT} gives the semisimplicity and saturatedness of the reduction for the direct sum for all $\ell\gg0$. This in turn should allow one to deduce the same for the factor corresponding to $\rho_\ell$, and hence \autoref{Thm-AI-Reduction-IsSaturated} for $\rho_\bullet$.
\end{Rem}

\section{Independence of lattices and main result on saturation}
\label{Sec-IndepOfLat}

Suppose $\rho_\bullet$ is $E$-rational compatible system which is a direct sum of absolutely irreducible compatible systems and has semisimple monodromy. In the first subsection we shall explain that the same monodromy groups, even integrally, arise from an absolutely irreducible compatible system, up to a finite kernel. In the second subsection, this will allow us in \autoref{Cor-Reduction-IsSaturated} to generalize \autoref{Thm-AI-Reduction-IsSaturated} to any system $\rho_\bullet$, i.e., the integral monodromy groups are smooth semisimple and their special fibers are the saturation of $\bar\rho_\bullet$, for almost all $\lambda\in\CP_E'$.  This generalizes some of the main results from \cite{CHT} from $\BQ$- to $E$-rational compatible systems. Our \autoref{Cor-CHT} gives an alternative proof of \cite[Thm.~1.2]{CHT}. In \autoref{Rem-OurMainThmFromCHToverQ} we indicate on argument by which their main result implies \autoref{Cor-Reduction-IsSaturated} for $\BQ$-rational systems. Finally, in \autoref{Cor-MonodromyIsUnramified} we deduce that the groups $G_{\rho_\bullet,\lambda}$ are unramified for almost all $\lambda\in\CP_E'$, as conjectured by Larsen and Pink (cf.~\cite[Conj.~5.4]{LarsenPink95}).

\subsection{Saturation under changing representations}
In this subsection we prove a criterion for when two representations whose generic monodromy groups differ by a central isogeny have saturations of their reductions that also differ by such an isogeny. Throughout this subsection, we let $\ell$ be in $\CP_\BQ'$ with $\ell>2$ and we let $K/\BQ_\ell$ be a finite extension with the ring of integers $\CO$, residue field $k$ and maximal ideal~$\Fm$. 

\begin{Lem}\label{Lem-IndepOfLattice}
Let $M$ be a semisimple connected group over $K$. Let $\alpha'\colon M\to\GL_{n'}$ be a faithful, and $\alpha\colon M\to\GL_n$ be an almost faithful representation. Let $\phi\colon\pi_1(X)\to M(K)$ be a homomorphism with Zariski dense image. Define $\rho=\alpha\circ\phi$, $\rho'=\alpha'\circ\phi$ and $\rho''=\rho\oplus\rho'$, let $\Lambda\subset K^n$ and $\Lambda'\subset K^{n'}$ be $\pi_1(X)$-stable $\CO$-lattices, and set $\Lambda'':=\Lambda\oplus\Lambda'$. Denote by $\bar\rho^\Lambda$, $\bar\rho^{\prime,\Lambda'}$ and $\bar\rho^{\prime\prime,\Lambda''}$ the corresponding reductions modulo $\Fm$, cf.\ the beginning of \autoref{Subsec-CharAndLifting}. Denote by $\CM_K$, $\CM'_K$ and $\CM''_K$ the image of $M$ under $\alpha$, $\alpha'$ and $\alpha\oplus\alpha'$, respectively, and by $\CM$, $\CM'$, $\CM''$ the $\CO$-group schemes that are the Zariski closures of $\CM_K$ in $\Aut_\CO(\Lambda)$, of $\CM'_K$ in $\Aut_\CO(\Lambda')$, and of $\CM''_K$ in $\Aut_\CO(\Lambda)\!\times\!\Aut_\CO(\Lambda')\subset\Aut_\CO(\Lambda'')$, respectively. Suppose that 
\begin{enumerate}
\item[(i)] the group schemes $\CM$, $\CM'$, $\CM''$ are smooth over $\CO$;
\item[(ii)] the identity components of their special fibers $\CM_k$, $\CM'_k$ and $\CM''_k$ are saturated;
\item[(iii)] the kernel of $\alpha$ is of order strictly less than~$\ell$.
\item[(iv)] the saturation $\bar\rho^\Lambda(\pi^\geo_1(X))^\sat_k$ is semisimple connected and agrees with $\CM_k$.
\end{enumerate}
Denote by $\pr'\colon\Aut_\CO(\Lambda)\times\Aut_\CO(\Lambda')\to \Aut_\CO(\Lambda')$ and $\pr\colon\Aut_\CO(\Lambda)\times\Aut_\CO(\Lambda')\to \Aut_\CO(\Lambda)$ the natural projection homomorphisms. We add the index $\CO$ for their restrictions $\pr_\CO\colon \CM''\to\CM$ and $\pr'_\CO\colon\CM''\to\CM'$. Then the following hold:
\begin{enumerate}
\item $\CM'$ over $\CO$ and $\CM'_k$ are connected semisimple, and the analogous assertion holds for~$\CM''$.
\item $\CM_k$, $\CM'_k$ and $\CM''_k$ are saturated in $\GL_{n,k}$, $\GL_{n',k}$ and $\GL_{n+n',k}$, respectively.
\item $\pr'_\CO$ is an isomorphism; $\pr_\CO$ is a central \'etale isogeny, and an isomorphism if $\alpha$ is faithful.
\item The inclusions $\bar\rho^{\prime\prime,\Lambda''}\!(\pi^\geo_1(X))^\sat_k\subseteq\CM''_k$ and $\bar\rho^{\prime,\Lambda'}\!(\pi^\geo_1(X))^\sat_k\subseteq\CM'_k$ from (b) \hbox{are equalities.}
\end{enumerate} 
\end{Lem}
\begin{proof}
Denote by $\pr_k\colon\CM_k''\to\CM_k$ and $\pr'_k\colon\CM_k''\to\CM_k'$ the morphisms induced from $\pr_\CO$ and $\pr_\CO'$ on special fibers, and consider the diagrams
\[\xymatrix@R-1pc @C+1pc
{
&\CM\ar@{^{(}->}[r]&\Aut_\CO(\Lambda)
\\
\pi_1(X)\ar@/^/[ur]^-{\rho} \ar@/_/[dr]_{\rho'}  \ar[r]^-{\rho''}&\CM'' \ar[u]^-{\pr_\CO}\ar[d]_{\pr_\CO'} \ar@{^{(}->}[r]&\Aut_\CO(\Lambda)\times \Aut_\CO(\Lambda') 
\ar[u]^-{\pr}\ar[d]_{\pr'}  
\\
&\CM'\ar@{^{(}->}[r]&\Aut_\CO(\Lambda')
\rlap{,}\\
\\
} \]
\vspace*{-2em}
\[\xymatrix@R-1pc @C+1pc
{
&\bar\rho^\Lambda(\pi_1^\geo(X))^\sat_k\ar@{^{(}->}[r]^-\iota& \CM_k\\
\pi_1(X)\ar@/^/[ur]^-{\bar\rho^\Lambda} \ar@/_/[dr]_-{\bar\rho^{\prime,\Lambda'}}  \ar[r]^-{\bar\rho^{\prime\prime,\Lambda''}}&\bar\rho^{\prime\prime,\Lambda''}(\pi_1^\geo(X))^\sat_k \ar[u]^-{\pr_k^\sat}\ar[d]_{(\pr')_k^\sat}  \ar@{^{(}->}[r]^-{\iota''}& \CM''_k\ar[u]^-{\pr_k}\ar[d]_{\pr_k'} 
\\
&\bar\rho^{\prime,\Lambda'}(\pi_1^\geo(X))^\sat_k\ar@{^{(}->}[r]^-{\iota'}& \CM'_k\rlap{,}\\
\\
} \]
where on the bottom we display the relevant situation on special fibers, and on the top over $\CO$. The closed immersion $\iota$ exists by hypotheses, the existence of $\iota',\iota''$ still needs a proof.

Since $\alpha$ is almost faithful and $\alpha'$ is faithful, the maps $M\to \CM''_K$ and $M\to \CM'_K$ are isomorphisms and $M\to \CM_K$ is a central isogeny with finite kernel. In particular $\dim M=\dim \CM_K=\dim \CM'_K=\dim \CM''_K$, and by (i) the groups $\CM$, $\CM'$ and $\CM''$ all have the same relative dimension over $\CO$. Since the restriction of $\pr_\CO$ and $\pr_\CO'$ to the generic fibers are clearly surjective, the maps $\pr_\CO'$ and $\pr_\CO$ are surjective homomorphisms by the definition of $\CM$, $\CM'$ and $\CM''$. Hence $\pr_k$ and $\pr'_k$ are surjective, as well. For dimension reasons, their kernels must be finite. The unipotent radical of $(\CM_k'')^o$ maps to a normal unipotent subgroup of $\CM_k$, and the induced map has finite kernel. Since $\CM_k$ is connected semisimple, the unipotent radical of $(\CM_k'')^o$ must be trivial and hence $(\CM_k'')^o$ is reductive. An analogous argument applies to centers and it follows that $(\CM_k'')^o$ is semisimple. From the surjectivity of $\pr'_k$ it follows similarly that $(\CM_k')^o$ is semisimple. As in the proof of \autoref{Thm-AI-Reduction-IsSaturated}, by \cite[Prop.~3.1.3]{ConradSGA3} it follows that $(\CM')^o$ and $(\CM'')^o$ are semisimple, and open and closed in $\CM'$ a $\CM''$, respectively. Since by assumption their generic fibers are connected, this proves (a). Part (b) follows immediately from~(a) and~(ii). Note that (b) also gives the existence of $\iota'$ and $\iota''$ as well as the inclusions in~(e).

Concerning (c) note first that by \cite[Prop.~6.1.10]{ConradSGA3} the maps $\pr_\CO$ and $\pr_\CO'$ in (c) and (d) are isogenies. Hence their kernels are finite flat over $\CO$ by \cite[Def.~3.3.9]{ConradSGA3}. By flatness and considering generic fibers, it follows that the kernel of $\pr_\CO'$ is trivial so that $\pr_\CO'$ is an isomorphism and that the kernel of $\pr_\CO$ is flat over $\CO$ of order strictly less than $\ell$. In particular $\kernel\pr_\CO$ is finite \'etale over $\CO$, so that $\kernel \pr_k$ is separable and hence central. Now by \cite[Prop.~3.3.10]{ConradSGA3} $\pr$ is central. If moreover $\alpha$ is faithful, then $\pr_\CO$ is an isomorphism by the proof for $\pr_\CO'$. This proves~(c).

We finally prove (d). We first note that by \autoref{Prop-GConnSemisimple2}, the maps $\pr_k^\sat$ and $(\pr')_k^\sat$ are well-defined and surjective. By (iv) it follows that 
\[\dim\bar\rho^{\prime\prime,\Lambda''}(\pi_1^\geo(X))^\sat_k\ge\dim\CM_k\] 
and from the existence of $\iota''$ and from (c) the reverse equality follows. Since $\CM_k''$ is semsimple connected, every proper closed subscheme has smaller dimension, and we find that $\iota''$ is an isomorphism. We already observed that $(\pr')_k^\sat$ is surjective, and that $\pr'_k$ is an isomorphism. Hence we deduce $\dim\bar\rho^{\prime,\Lambda'}(\pi_1^\geo(X))^\sat_k\ge\dim\CM''_k =\dim \CM_k'$, where the equality is due to (c). Arguing as above, the inclusion $\iota'$ is the identity as well. This completes (d).
\end{proof}
\begin{Cor}
Let $\rho\colon \pi_1(X)\to \GL_n(K)$ be a continuous homomorphism, and let $\Lambda$ and $\Lambda'$ be $\CO$-lattices in $K^n$ that are stable under the action of $\pi_1(X)$ via $\rho$. Let $\CM$ and $\CM'$ be the $\CO$-group schemes that are the Zariski closures of $\rho(\pi_1(X))$ in $\Aut_\CO(\Lambda)$ and $\Aut_\CO(\Lambda')$, respectively, and let $\CM''$ be the Zariski closure of $(\rho\oplus\rho)(\pi_1(X))\subseteq \Aut_\CO(\Lambda)\times \Aut_\CO(\Lambda')$. Suppose that
\begin{enumerate}
\item[(i)] the generic fiber $M:=\CM\otimes_\CO K(\cong \CM'\otimes_\CO K)$ is semisimple and connected;
\item[(ii)] the $\CO$-group schemes $\CM$, $\CM'$ and $\CM''$ are smooth; 
\item[(iii)] the identity components of their special fibers $\CM_k$, $\CM'_k$ and $\CM''_k$ are saturated;
\item[(iv)] the saturation $\bar\rho^\Lambda(\pi^\geo_1(X))^\sat_k$ is semisimple and agrees with $\CM_k$.
\end{enumerate}
Then $\bar\rho^{\Lambda'}(\pi^\geo_1(X))^\sat_k$ is isomorphic to $\CM'_k$, $\CM\cong\CM'$ and both are semisimple connected over~$\CO$.
\end{Cor}
\begin{proof}
Let $g\in\GL_n(K)$ be such that $\Lambda'=g\Lambda$, so that $\Aut_\CO(\Lambda')=g \Aut_\CO(\Lambda')g^{-1}$. The result is now a direct consequence of \autoref{Lem-IndepOfLattice}  applied to $\rho$ and $\rho'=g\rho g^{-1}$.
\end{proof}

\subsection{Geometric mod $\lambda$ semisimplicity for almost all $\lambda$.}

Let $\rho_\bullet$ be an $E$-rational compatible system. For each $\lambda\in\CP_E'$, let $\Lambda_\lambda\subset E_\lambda^n$ be a $\pi_1(X)$-stable $\CO_\lambda$ lattice. 
Define $\CG_\lambda$ (resp.\ $\CG_\lambda^\geo$) as the Zariski closure of $G_\lambda$ (resp.\ $G^\geo_\lambda$) in $\Aut_{\CO_\lambda}(\Lambda_\lambda)\cong\GL_n(\CO_\lambda)$. Define $\CG_{k_\lambda}:=\CG_\lambda\otimes_{\CO_\lambda}k_\lambda$ and $\CG^\geo_{k_\lambda}:=\CG^\geo_\lambda\otimes_{\CO_\lambda}k_\lambda$. The following is the main result of this section. It generalizes the part of \autoref{Thm-AI-Reduction-IsSaturated} concerning geometric monodromy groups to compatible systems that are not necessarily absolutely irreducible. 
\begin{Thm}\label{Thm-Reduction-IsSaturated}
Suppose that $\rho_\bullet|_{\pi_1^\geo(X)}$ is connected.\footnote{Because of \autoref{Cor-SSandGeomMonodromy}, here we do not need to assume that $\rho_\bullet$ is semisimple.} Then for almost all $\lambda\in\CP_E'$ 
\begin{enumerate}
\item the group $\CG^\geo_\lambda$ is  semisimple and connected over $\CO_\lambda$, and $\bar\rho^\Lambda_\lambda(\pi_1^\geo(X))^\sat_{k_\lambda}=\CG^\geo_{k_\lambda}$, and
\item $\bar\rho^\Lambda_\lambda$ is a semisimple representation, and thus isomorphic to $\bar\rho_\lambda$.
\end{enumerate}
\end{Thm}
\begin{proof}
By \autoref{Conv-Semisimple} $\rho_\bullet$ is semisimple, and using the connectivity hypothesis, it is clear that we may pass from the base field $\kappa$ to a finite extension, so that we can assume that $\rho_\bullet$ is connected as well. We also observe that to prove the theorem, we may pass from $E$ to a finite extension $E'$: If $\lambda$ is unramified for $E'/\BQ$, then the formation of the Zariski closure commutes with the unramified base change $\CO_\lambda\to \CO'_{\lambda'}$, and being semisimple connected descends from the base change to the original situation. Moreover we are given closed immersions $\bar\rho_\lambda^\Lambda(\pi_1^\geo(X))^\sat_{k_\lambda}\into\GL_{n,k}$ and $\CG^\geo_{k_\lambda}\into\GL_{n,k}$, and if these immersions are the same under a base change $k\to k'$, then the given immersions are the same: if $I_{k'},J_{k'}$ are the kernels in the coordinate ring $k'[\GL_n]$ describing these immersions, then $I_{k'}=J_{k'}$, and their invariants under $\Gal(k'/k)$ are the original kernels and hence the same. 

Thus from now on we also assume that $\rho_\bullet$ is a direct sum of absolutely irreducible representations $\bigoplus_i\rho_{i,\bullet}$. Since $\rho_\bullet$ is connected, as in the proof of \autoref{Prop-WeightZeroArithMonodromy} one has surjections $G_{\rho_\bullet,\lambda}\to G_{\rho_{i,\bullet},\lambda}$, and so the $\rho_{i,\bullet}$ are connected. Hence by \autoref{Lem-OnAbsIrredAndTwists} we may twist each $\rho_{i,\bullet}$ by an unramified compatible system $\delta_{i,\bullet}$, such that each $\rho_{i,\bullet}\otimes\delta_{i,\bullet}$ has finite order determinant, and hence trivial determinant by \autoref{Lem-AIandFiniteDet}. Since the twist does not change $\rho_\bullet|_{\pi_1^\geo(X)}$, we assume from now on that each summand $\rho_{i,\bullet}$ has trivial determinant. 

If necessary, we now enlarge $E$ again, so that a split motivic triple is defined over $E$; denote it by $(E,M,\alpha)$. Then by \autoref{Cor-ExistsAIRep}(d) we have an almost faithful representation $\beta$ of $M$, and a  compatible system $\rho'_{\bullet}$ such that for every $\lambda\in\CP_E'$ we have a homomorphism $\phi_\lambda\colon \pi_1(X)\to M(E_\lambda)$ and a commuting diagram
\[\xymatrix@R-1pc@C+1pc{
&&\GL_n(E_\lambda)\\
\pi_1(X)\ar@/^/[urr]^{\rho_\lambda} \ar@/_/[drr]_{\rho'_\lambda}  \ar[r]^{\phi_\lambda}&M(E_\lambda)\ar[ur]_{\alpha\otimes_EE_\lambda}\ar[dr]^{\beta\otimes_EE_\lambda}&\\
&&\GL_{n'}(E_\lambda)\rlap{.}
} \]
Let $\Lambda'_\lambda\subset E_\lambda^{n'}$ be an $\CO_\lambda$-lattice stable under the action of $\pi_1(X)$ under $\beta\otimes_EE_\lambda$, so that $\Lambda''_\lambda:=\Lambda_\lambda\oplus\Lambda'_\lambda\subset E^{n+n'}_\lambda$ is stable for $\pi_1(X)$ acting via $(\alpha\oplus\beta)\otimes_EE_\lambda$. Define for $?\in \{\lambda,k_\lambda\}$ and $*\in\{\prime,\prime\prime\}$ the group schemes $\CG^*_?$ and $\CG^{*,\geo}_?$ starting from $\Lambda_\lambda^*$ in the same way as $\CG_?$ and $\CG^\geo_?$ was constructed. By \autoref{Prop-LP-GisSmooth}, the $\CO_\lambda$-group schemes $\CG_\lambda,\CG'_\lambda,\CG''_\lambda$ are smooth for almost all $\lambda\in\CP_E'$. By \autoref{Thm-RedIsSaturated}, the identity components of their special fibers $\CG_{k_\lambda},\CG'_{k_\lambda},\CG''_{k_\lambda}$ are saturated for almost all $\lambda\in\CP_E'$. By \autoref{Thm-AI-Reduction-IsSaturated}, for almost all $\lambda\in\CP_E'$ one has 
\[\bar\rho^{\prime,\Lambda'}_\lambda(\pi_1^\geo(X))^\sat_{k_\lambda}=\CG^{\prime,\geo}_{k_\lambda},\]
and both groups are semisimple connected over $k_\lambda$. Note also that the size of kernel of $\beta\otimes_EE_\lambda$ is clearly independent of $\lambda$. Thus by \autoref{Lem-IndepOfLattice} (with the roles of $\rho'$ and $\rho$ reversed), we have for almost all $\lambda\in\CP_E'$ that $\bar\rho^\Lambda_\lambda(\pi^\geo_1(X))^\sat_k$ is semisimple connected and equal to $\CG^{\geo}_{k_\lambda}$ and that the group scheme $\CG^\geo_\lambda$ is semisimple connected over~$\CO_\lambda$. This proves part (a). Part (b) is immediate from (a) and \autoref{Prop-GConnSemisimple}.
\end{proof}

The following result summarizes the $\BQ$-rational case. 
\begin{Cor}\label{Cor-OpennessAboveEll0}
Let $\rho_\bullet$ be $\BQ$-rational connected and with semisimple split motivic group. Then for almost all $\ell\in\CP_\BQ'$ the following hold:
\begin{enumerate}
\item The $\BZ_\ell$-group scheme $\CG_\ell$ is smooth semisimple connected.
\item One has $\big( \bar\rho^{\Lambda}_{\ell}(\pi_1^\geo(X))\big)^\sat_{\BF_\ell}= \CG_{\BF_\ell}$, and hence $\big( \bar\rho^{\Lambda}_{\ell}(\pi_1^\geo(X)) \big)^+=  \CG_\ell(\BF_\ell)^+$.
\item For almost all $\ell$, one has $\CG_\ell(\BZ_\ell)^+ \subseteq \rho_\ell (\pi_1^\geo(X)) \subseteq  \CG_\ell(\BZ_\ell).$
\item Assertions (b) and (c) hold with $\pi_1(X)$ in place of $\pi_1^\geo(X)$.
\item There exists a finite subset $\CL$ of $\CP_\BQ'$, a compact open subgroup $H^\geo_\CL$ of $\prod_{\ell\in\CL}G_\ell(\BQ_\ell)$ and a finite cover $X'\to X$ such that \[
\big(\prod_{\ell\in\CP_\BQ'}\rho_{\ell}\big)(\pi^\geo_1(X'))=\prod_{\ell\in\CP_\BQ'\setminus\CL} \CG_{\ell}(\BZ_\ell)^+
\times H^\geo_\CL.\]
\end{enumerate}
\end{Cor}
\begin{proof}
By our hypotheses $G_\ell=G_\ell^o=G^\geo_\ell$, and now part (a) is immediate from \autoref{Thm-Reduction-IsSaturated}, as is the first half of part~(b). To see the second assertion in (b), we may pass to a cover of $X$, so that by \autoref{Thm-AlmostIndependence}, we can assume that $\bar\rho^\Lambda_\ell(\pi_1^\geo(X))$ is $\ell$-generated for almost all $\ell$. In that case,  by \autoref{Rem-NoriEnvelope} the group $\bar\rho^\Lambda_\ell(\pi_1^\geo(X))^\sat_{\BF_\ell}$ is the Nori envelope of $\bar\rho^\Lambda_\ell(\pi_1^\geo(X))$. The equality now follows from \cite[Thm.~B]{Nori}, and this concludes the proof of~(b).

In (c), only the inclusion on the left is at stake. By (b), we know the left inclusion after reduction modulo $\ell$. Let $\pi_\ell\colon\wt\CG_\ell\to\CG_\ell$ be the universal cover of $\CG_\ell$. Then by \cite[Thms.~12.4 and 12.6]{Steinberg-Endomorphisms}, the image of $\wt\CG_\ell(\BF_\ell)$ in $\CG_\ell(\BF_\ell)$ is precisely $\CG_\ell(\BF_\ell)^+$. Let $\Gamma\subset \wt\CG_\ell(\BZ_\ell)$ be the inverse image of $\rho^\Lambda_\ell(\pi_1^\geo(X))$ under $\pi_\ell$. Then by (b), the result just quoted and \cite[Prop.~2.6]{Larsen95}, we have for all $\ell\gg0$ that $\Gamma=\wt\CG_\ell(\BZ_\ell)$, and this implies~(c). Assertion (d) is clear since the claims in (b) and (c) are weaker if one replaces $\pi^\geo_1(X)$ by $\pi_1(X)$.

For (e) let $X'\to X$ be a finite cover as in \autoref{Thm-AlmostIndependence} so that $\rho_\bullet|_{\pi_1(X')}$ is independent and almost all $\rho_\ell (\pi_1^\geo(X'))$ are $\ell$-generated. Going through the above arguments we see that there is a finite set $\CL$ of primes such that for all $\ell\notin\CL$ we now have $\bar\rho^{\Lambda}_{\ell}(\pi_1^\geo(X')) =  \CG_\ell(\BF_\ell)^+$ in (b) and hence $\rho_\ell (\pi_1^\geo(X'))=\CG_\ell(\BZ_\ell)^+ $ in (c). Define $H_\CL:=\prod_{\ell\in\CL}\rho_\ell(\pi_1^\geo(X))$. Then (e) follows using the independence of $\rho_\bullet$.
\end{proof}

\begin{Cor}\label{Cor-GeomImageOpenInArithImage0}
Suppose that $\rho_\bullet$ is $E$-rational with semisimple split motivic group. Then 
\[(\prod_{\lambda\in\CP_E'}\rho_\lambda\big)\big(\pi_1^\geo(X)\big)\subset (\prod_{\lambda\in\CP_E'}\rho_\lambda\big)\big(\pi_1(X)\big)\] is an open subgroup.
\end{Cor}
\begin{proof}
It suffices to consider $\rho'_\bullet:=\Res_{E/\BQ}\rho_\bullet$. By \autoref{Lem-FactsOnChin} its split motivic group is semisimple, and hence we shall assume $E=\BQ$ from now on and work with $\rho_\bullet$. For the proof of the corollary, we may pass to a finite cover, and hence by \autoref{Thm-SerreOnComponents} we may assume that $\rho_\bullet$ is also connected. By \autoref{Thm-AlmostIndependence} we may further assume that $\rho_\bullet$ is independent. Now from \autoref{Cor-OpennessAboveEll0}(d) and (e), it follows that there exists a finite subset $\CL$ of $\CP_\BQ'$ and compact open subgroups $H^\geo_\CL\subset H_\CL\subset \prod_{\ell\in\CL}G_\ell(\BQ_\ell)$ such that 
\[\prod_{\ell\in\CP_\BQ'\setminus\CL} \!\CG_{\ell}(\BZ_\ell)^+
\times H_\CL^\geo=\big(\prod_{\ell\in\CP_\BQ'}\big)(\pi^\geo_1(X))  \stackrel{(\ast)}\subseteq\big(\prod_{\ell\in\CP_\BQ'}\big)(\pi_1(X))\subseteq 
\prod_{\ell\in\CP_\BQ'\setminus\CL} \!\CG_{\ell}(\BZ_\ell) \times H_\CL.
\]
Because $\CG_\ell$ is semisimple for $\ell\in\CP_\BQ'\setminus\CL$ and contained in $\GL_{n,\BQ_\ell}$, the quotients $\CG_\ell(\BZ_\ell)/\CG_\ell(\BZ_\ell)^+$ are finite abelian of cardinality bounded independent of $\ell$, e.g.~\cite[Lem.~3.8]{BGP15}. Since $H_\CL^\geo\subset H_\CL$ is open its index is finite. Now $\pi_1(X)/\pi^\geo_1(X)\cong\wh\BZ$, and so there is a finite extension $\kappa'/\kappa$ after which the inclusion labelled $(\ast)$ becomes an equality for $X'=X\otimes_\kappa\kappa'$ in place of~$X$. 
\end{proof}

\begin{Cor}\label{Cor-Reduction-IsSaturated}
Suppose that $\rho_\bullet$ is an $E$-rational compatible system of representations of $\pi_1(X)$, and that the lattices $\Lambda_\lambda$, $\lambda\in\CP_E'$ are stable under $\pi^\geo_1(X)$ (but not necessarily $\pi_1(X)$). Then for almost all $\lambda\in\CP_E'$ the group $\CG^{\geo,o}_\lambda$ is  semisimple over $\CO_\lambda$, and $\bar\rho^\Lambda_\lambda(\pi_1^\geo(X))^{\sat,o}_{k_\lambda}=\CG^{\geo,o}_{k_\lambda}$.
\end{Cor}
Note that because of  \autoref{Cor-SSandGeomMonodromy} the corollary also holds if $\rho_\bullet$ is not semisimple.
\begin{proof}
Either by the above remark, or by \autoref{Conv-Semisimple}, we shall assume that $\rho_\bullet$ is semisimple. Clearly it suffices to prove the result after passing to a finite coefficient extension $E'/E$, since $E'/E$ is unramified for almost all $\lambda$: the formation of the reduction is clearly compatible with this (cf.~\autoref{Rem-RedAndCoeffChange}); but also the formation of $\CG_\lambda^\geo$ as a Zariski closure, since taking the Zariski closure commutes with base change under unramified extensions $\CO'_{\lambda'}/\CO_\lambda$; cf.\ the proof of \cite[Prop.~1.3]{LarsenPink95}; and finally the formation of saturation is compatible with base change (and descent) by \autoref{Cor-SatBaseChange}. Thus we assume from now on that $\rho_\bullet$ is a direct sum of absolutely irreducible compatible systems $\rho_{i,\bullet}$. We may also twist the $\rho_{i,\bullet}$ according to \autoref{Lem-OnAbsIrredAndTwists}(c) by characters of $\Gamma_\kappa$ so that the $\rho_{i,\bullet}$ have finite order determinant, since the assertion we prove only concerns~$\pi_1^\geo(X)$. 

Since the assertion of the corollary only concerns identity components, we may also freely replace $X$ by a finite cover. Thus by \autoref{Thm-SerreOnComponents} we may assume that $\rho_\bullet$ and the $\rho_{i,\bullet}$ are connected, and have trivial determinant. This might destroy the absolute irreducibility of the $\rho_{i,\bullet}$. But by passing through the above steps finitely many times (the degree of $\rho_\bullet$ is finite), we may finally assume that the $\rho_{i,\bullet}$ are absolutely irreducible, and then by \autoref{Prop-WeightZeroArithMonodromy} that $\rho_\bullet$ is connected and has semisimple split motivic group. Now by \autoref{Cor-GeomImageOpenInArithImage0}, and a further finite cover we may assume that $(\prod_{\lambda\in\CP_E'}\rho_\lambda\big)\big(\pi_1^\geo(X)\big)=(\prod_{\lambda\in\CP_E'}\rho_\lambda\big)\big(\pi_1(X)\big)$.~In~particular, the lattices $\Lambda_\lambda$ are stable under (the new) $\pi_1(X)$. Now the corollary follows from \autoref{Thm-Reduction-IsSaturated}. 
\end{proof}

\autoref{Thm-Reduction-IsSaturated} also yields a proof of a main result of \cite{CHT} by independent methods, as is shown by the following result and its proof. 
\begin{Cor}[{\cite[Thm.~1.2]{CHT}}]\label{Cor-CHT}
Let $k$ be an algebraically closed field of characteristic $p$, let $X$ be a smooth $k$-scheme and let $f\colon Y\to X$ be a smooth proper morphism. Then the action of $\pi_1(X)$ on $R^if_*\BF_\ell$ is semisimple for almost all $\ell\in\CP_\BQ'$. 
\end{Cor}
\begin{proof}
By a standard reduction argument, cf.~\cite[Prop.~4.1]{CHT}, it suffices to prove the result in the case when $k=\kappa$ for some finite extension $\kappa$ of $\BF_p$, and when $f$ is the base change under $\kappa\to k$ of a smooth morphism $f_0\colon Y_0\to X_0$ between smooth geometrically connected $\kappa$-schemes. In the latter case, by \cite[Cor.~3.3.9]{Deligne-Weil2}, for each $i$, the family $\rho_\bullet$ of Galois representations attached to the family of lisse sheaves $(R^if_{0*}\BQ_\ell)_{\ell\in\CP_\BQ'}$ on $X$ is a $\BQ$-rational compatible system, pure of weight~$i$. It carries a natural integral structure given by the family of $\BZ_\ell$-lattices $\Lambda_\ell$ defined by (the generic fiber of) $R^if_{0*}\BZ_\ell$, $\ell\in\CP_\BQ'$. Let $\bar\rho_\ell$ be the induced representation of $\pi_1(X)$ on $\Lambda_\ell/\ell\Lambda_\ell$. 

To prove the corollary, we may clearly pass to a finite extension of $X$, and hence by \autoref{Thm-AlmostIndependence} we may assume that the images $\rho_\ell(\pi^\geo_1(X))$ are pro-$\ell$ generated for almost all~$\ell$. From \autoref{Thm-Reduction-IsSaturated} we deduce that $\bar\rho_\ell^\Lambda(\pi_1^\geo(X))^\sat_{\BF_\ell}=\CG^\geo_{\BF_\ell}$ is semisimple for almost all $\ell$, where $\CG^\geo_\ell$ is the Zariski closure of $\rho_\ell(\pi^\geo_1(X))$ in $\Aut_{\BZ_\ell}(\Lambda_\ell)$. Now the action of $\pi^\geo_1(X)$ on $R^if_*\BF_\ell$ is via $\bar\rho_\ell^\Lambda$, and since its saturation is semisimple, it follows from the last assertion of \autoref{Prop-GConnSemisimple} that $\pi_1^\geo(X)$ acts semisimply via $\bar\rho_\ell^\Lambda$, as was to be shown.
\end{proof}

\begin{Rem}\label{Rem-OurMainThmFromCHToverQ}
Combining \autoref{Rem-AbsIrrThmFromCHToverQ} with the reduction argument given in \autoref{Thm-Reduction-IsSaturated} should prove \autoref{Thm-Reduction-IsSaturated} for arbitrary $\BQ$-rational compatible systems $\rho_\bullet$ with semisimple split motivic group from the results in \cite[Thm.~1.2]{CHT} by Cadoret, Hui and Tamagawa without invoking any results from \autoref{Sec-ResidualSaturation}.
\end{Rem}

The following result answers in the affirmative Conjecture~5.4 of \cite{LarsenPink95} for representations of $\pi_1(X)$ with $X$ a normal variety over $\BF_p$.
\begin{Cor}\label{Cor-MonodromyIsUnramified}
Suppose that $\rho_\bullet$ is an $E$-rational compatible system. Then for almost all  $\lambda\in\CP_E'$, the groups $G^{\geo,o}_{\rho_\bullet,\lambda}$ and $G^{o}_{\rho_\bullet,\lambda}$ are unramified, i.e., they are quasi-split and become split over an unramified extension.
\end{Cor}
\begin{proof}
The assertions are unaffected by replacing $X$ by a finite cover, and so by \autoref{Thm-SerreOnComponents} and \autoref{Thm-AlmostIndependence} we may assume that $\rho_\bullet$ and $\rho_\bullet|_{\pi^\geo_1(X)}$ are connected, and that $\rho_\lambda(\pi^\geo_1(X))$ is $\ell_\lambda$-generated for almost all $\lambda\in\CP_E'$. By \autoref{Thm-Reduction-IsSaturated}, for almost all $\lambda$ the group $G^{\geo}_{\rho_\bullet,\lambda}$ is the generic fiber over a connected semisimple group scheme over $\CO_\lambda$. It follows that group scheme as well as its generic fiber are quasi-split and become split over an unramified extension; see \cite[5.2.14]{ConradSGA3}. This proves the corollary for $G^{\geo,o}_{\rho_\bullet,\lambda}$.

Next observe that by \autoref{Lem-GDer=Ggeo}, we have $G^{\geo,o}_{\rho_\bullet,\lambda}=G^\der_{\rho_\bullet,\lambda}$, since in our situation both are connected. But now if $B$ is a Borel subgroup in $G^{\geo}_{\rho_\bullet,\lambda}$, then its image $B^\ad$ in the adjoint group $(G^o_{\rho_\bullet,\lambda})^\ad$ is a Borel subgroup, and hence the preimage of $B^\ad$ under $G^o_{\rho_\bullet,\lambda}\to (G^o_{\rho_\bullet,\lambda})^\ad$ is a Borel subgroup in $G^o_{\rho_\bullet,\lambda}$. Thus whenever $G^{\geo,o}_{\rho_\bullet,\lambda}$ is quasi-split, then so is $G^o_{\rho_\bullet,\lambda}$. That $G^o_{\rho_\bullet,\lambda}$ becomes split over an unramified extension of $E_\lambda$ for almost all $\lambda$ is a consequence of the existence of a split motivic triple over a finite extension $F$ of $E$, since $F/E$ is unramified for almost all $\lambda$.
\end{proof}

\section{Adelic openness for some $E$-rational compatible systems}
\label{Sec-OpenImage}

Let $\rho_\bullet$ be an absolutely irreducible connected $E$-rational with simple adjoint split motivic group $M$, cf.~\autoref{Def-ChinAdjoint}, and suppose that the split motivic representation $\alpha$ is the adjoint representation of $M$.  We shall deduce adelic openness in the sense of Hui-Larsen for the algebraic and geometric monodromy groups attached to such a $\rho_\bullet$. We shall make crucial use of \cite{Pink-Compact} by Pink. In \autoref{Rem-AdelicOpen} we shall also indicate with a sketch of proof how far one can push adelic openness results by the methods developed here and in \cite{Pink-Compact}. From the results in this section, starting from \autoref{Cor-OpennessAboveEll0}, only, we will also deduce a second proof of \autoref{Thm-AI-Reduction-IsSaturated}, and thus by the results of \autoref{Sec-IndepOfLat} of \autoref{Thm-Reduction-IsSaturated}. Hence by \autoref{Rem-OurMainThmFromCHToverQ} this proof builds, up to some unchecked compatibilities, directly on \cite{CHT} and avoids the use of any results from~\autoref{Sec-ResidualSaturation}.

\subsection{Summary of parts of Pink \cite{Pink-Compact}}
\label{Subsec-OnPink}
We recall some context from \cite{Pink-Compact} specialized to the case of local fields of characteristic zero. In the following by $L$ we denote a finite product $\prod_iL_i$ of local non-archimedean fields of characteristic zero. Different $L_i$ may or may not have the same residue characteristic. By a local non-archimedean field of characteristic zero we mean a finite extension of some field $\BQ_\ell$. We regard $L$ as a commutative semisimple ring. By $K$ we denote a (topologically) closed subring of $L$. It may have fewer components. By $G$ we denote a semisimple group over $L$ which is fiberwise absolutely simple adjoint; i.e., $G$ is a disjoint union of groups $\coprod_i G_i$ where each $G_i$ is absolutely simple adjoint over $L_i$. Note also that $G(L)=\prod_i G(L_i)$. Moreover by $\Gamma$ we denote a compact subgroup of $G(L)$ that is fiberwise Zariski dense.
\begin{Def}\label{Def-MinimalQM}
A {\em quasi-model} of $(L,G,\Gamma)$ is a triple $(K,H,\phi)$ consisting of
\begin{enumerate}
\item[(i)] a semisimple closed subring $K$ of $L$ over which $L$ is of finite type, 
\item[(ii)] a fiberwise absolutely simple adjoint group $H$ over $K$,
\item[(iii)] an isomorphism $\phi\colon H\otimes_KL\to G$,
\end{enumerate}
such that $\Gamma$ is contained in $\phi(H(K))$.\footnote{In \cite[\S~3]{Pink-Compact}, $\phi$ can be a totally inseparable isogeny and $\mathrm{d}\phi$ can vanish at some factors. Both are not possible in characteristic zero.}

The triple $(L,G,\Gamma)$ is {\em minimal} if for every quasi-model $(K,H,\phi)$ one has $L=K$.

A quasi-model $(K,H,\phi)$ is called {\em minimal}, if the triple $(K,H,\phi^{-1}(\Gamma))$ is minimal.
\end{Def}
Fix a triple $(L,G,\Gamma)$ as above. A main abstract result on minimal quasi-models is the following:
\begin{Thm}[{\cite[Thm.~3.6]{Pink-Compact}}]\label{Thm-Pink-MinQM}
\begin{enumerate}
\item For $(L,G,\Gamma)$ there exists a minimal quasi-model $(K,H,\phi)$.
\item The subring $K\subset L$ in (a) is unique, and $(H,\phi)$ are determined up to unique isomorphism.
\end{enumerate}
\end{Thm}
In particular, if $(K',H',\phi')$ is a quasi-model for $(L,G,\Gamma)$ and if $(K,H,\phi)$ is a minimal quasi-model for $(K',H',(\phi^{\prime})^{-1}(\Gamma))$, then $(K,H,\phi'\circ\phi)$ is a minimal quasi-model for $(L,G,\Gamma)$.

However \cite{Pink-Compact} gives also a fairly concrete description of a minimal quasi-model. In the following result, we only quote the part on the minimal field of definition. The proof follows immediately from \cite[Thms.~2.3 and 2.14]{Pink-Compact} and \cite[Prop.~3.14]{Pink-Compact}; we omit details.

\begin{Thm}[{cf.~\cite[Prop.~0.6]{Pink-Compact}}]\label{Thm-ExplicitMinQM}
Let $\Ad_G\colon G\to \Aut(\Lie G)$ be the fiberwise adjoint representation of $G$. Let $K_{\Ad,G/L}\subset L$ be the closure in $L$ of the $\BQ$-algebra generated by the set $\{\trace(\Ad_G(\gamma))\mid\gamma\in\Gamma\}$. Then there exists a minimal quasi-model of the form $(K_{\Ad,G/L},H,\phi)$.
\end{Thm}

Also relevant to us is the following result:
\begin{Prop}[{\cite[Cor.~3.8]{Pink-Compact}}]\label{Prop-QMforNormal}
Let $\Gamma'\subset\Gamma$ be a closed normal subgroup. Suppose that $\Gamma'$ is also fiberwise Zariski dense in $G$. If $(L,G,\Gamma)$ is minimal, then so is $(L,G,\Gamma')$.
\end{Prop}

The following is the main result of \cite{Pink-Compact} on open images; note that in characteristic zero isogenies are separable and $\phi$ is an isomorphism for any quasi-model $(K,H,\phi)$.
\begin{Thm}[{\cite[Thm.~0.2]{Pink-Compact}}]\label{Thm-PinkOpen}
If $(K,H,\phi)$ is a minimal quasi-model for $(L,G,\Gamma)$, then the inclusion $\Gamma\into G(L)$ factors via $H(K)$ and $\Gamma$ is open in~$H(K)$. Conversely, if $\Gamma\subset G(L)$ is open then $(L,G,\Gamma)$ is minimal.
\end{Thm}
The converse part is rather simple to prove by considering the $\BQ_\ell$ dimension of $\Lie G$.

\subsection{Simple adjoint compatible systems}
\label{Subsec-ChinAdjoint}
In this section, we shall deduce some adelic openness results for certain $E$-rational compatible systems which we call (simple) adjoint. Using the just quoted results from Pink, for such systems we can apply Weil restriction to get into the setting of $\BQ$-rational systems. To the latter we apply Corollaries~\ref{Cor-OpennessAboveEll0} and \ref{Cor-Reduction-IsSaturated}. We then also indicate how, starting from \autoref{Cor-OpennessAboveEll0}, one can prove \autoref{Thm-AI-Reduction-IsSaturated} for general $E$-rational systems from its analog for $\BQ$-rational systems.
\begin{Def}\label{Def-ChinAdjoint}
We call a compatible system {\em simple adjoint} if it is connected semisimple, if its split motivic group $M$ is simple adjoint and if its split motivic representation is $\Ad_M$. We call it {\em adjoint}, if its split motivic group $M$ is semisimple adjoint, say $M=\prod_i M_i$ with $M_i$ simple, and if its split motivic representation is $\bigotimes_i\Ad_{M_i}$.
\end{Def}

\begin{Ex}
Let $\rho_\bullet$ be an $E$-rational compatible system. Then after possibly enlarging $E$,
\begin{enumerate}
\item
there exists an $E$-rational compatible system $\wt\rho_\bullet$ whose monodromy groups $G_{\wt\rho_{\bullet},\lambda}$ are semisimple and of adjoint type, and such that for all $\lambda\in\CP_E'$ there exists a surjective homomorphism $G_{\rho_{\bullet},\lambda} \to G_{\wt\rho_{\bullet},\lambda}$ with central kernel, and
\item
there exists an adjoint compatible system whose split motivic group is the adjoint group of the split motivic group of $\rho_\bullet$.
\end{enumerate}

To see (a), by  \autoref{Prop-OnAbsoluteCR} we may assume that $\rho_\bullet$ is absolutely completely reducible  after enlarging $E$. Let $\rho_\bullet=\oplus_{i\in I}\rho_{i,\bullet}$ be the corresponding decomposition into absolutely irreducible compatible systems $\rho_{i,\bullet}$  of dimension $n_i$. Then $\wt\rho_\bullet:=\bigoplus\Ad_{\rho_{i,\bullet}}$ is an $E$-rational compatible system. Since the center of each $\rho_{i,\bullet}$ is the center of the ambient $\GL_{n_i}$ intersected with $G_{\rho_{i,\bullet},\lambda}$, the groups $G_{\wt\rho_{\bullet},\lambda}$ are semisimple and of adjoint type, and (a) follows readily.

For (b) let $\wt\rho_\bullet^M$ be the $M$-compatible system from \autoref{Cor-ChinRepOverM}, where again we possibly have to enlarge $E$ so that $\wt\rho_\bullet^M$ and a split motivic triple $(E,M,\alpha)$ for $\wt\rho_\bullet$ exist over $E$. Since $M$ is split simple adjoint, we have $M=\prod_i M_i$ with $M_i$ absolutely simple, split adjoint. Let $\beta:=\bigotimes_i \Ad_{M_i}$. Then $\beta\circ \wt\rho_\bullet^M$ is the searched-for $E$-rational adjoint compatible system for~(b).
\end{Ex}

\begin{Def}\label{Def-MinimalE}
If $\rho_\bullet$ is an $E$-rational compatible system, we call $E$ {\em minimal for $\rho_\bullet$}, if $E/\BQ$ is generated by the coefficients of the polynomials $P_x$, $x\in|X|$.
\end{Def}

We deduce several immediate consequences:
\begin{Cor}\label{Cor-FieldOfDefOfRhoBullet}
Suppose $\rho_\bullet$ is $E$-rational simple adjoint. Let $E^\prime\subset E$ be the field of definition of the polynomials $P_x$, $x\in|X|$, for $\rho_\bullet$. Then the following hold:
\begin{enumerate}
\item There exists an $E^\prime$-rational compatible system $\rho^\prime_\bullet$, such that $\rho_\bullet=\rho^\prime_\bullet\otimes_{E^\prime}E$, i.e., $E^\prime$ is minimal for $\rho^\prime_\bullet$.
\item For any $\ell\in\CP_\BQ'$, the images of $\pi_1(X)$ and $\pi_1^\geo(X)$ under $\bigoplus\limits_{\lambda'|\ell}\rho'_{\lambda'}$ are open in $\prod\limits_{\lambda'|\ell} G_{\rho'_{\bullet},\lambda'}(E'_{\lambda'})$.
\item The $\BQ$-rational compatible system $\rho^\BQ_\bullet:=\Res_{E'/\BQ}\rho'_{\bullet}$ is semisimple, and its monodromy groups are
\[G_{\rho^\BQ_{\bullet},\ell}=\prod_{\lambda'|\ell} \Res_{E'/\BQ_{\ell}} G_{\rho'_{\bullet},\lambda'}.\]
\end{enumerate}
\end{Cor}
\begin{proof}
Let $M$ be the split motivic group of $\rho_\bullet$ and let $n=\dim M$, so that $\rho_\bullet$ is $n$-dimensional. Let $E_0$ be the subfield of $E$ generated by the coefficients of $T^{n-1}$ of the polynomials $P_x\in E[T]$, $x\in|X|$. In the course of the proof we shall see that $E_0=E'$. For $\ell\in\CP_\BQ'$ consider $L=\oplus_{\lambda|\ell}E_\lambda=E\otimes_\BQ\BQ_\ell$, $G=\coprod_{\lambda'|\ell} G_{\rho_\bullet,\lambda}$ and $\rho_L=\oplus_{\lambda|\ell}\rho_\lambda=\Ad_G$ (by hypothesis). Let $K$ be the closure of $E_0$ in $L$. Then $K=E_0\otimes_\BQ\BQ_\ell$ and it is embedded into $L$ via the completion at $\ell$ of $E_0\into E$. Then for all $x\in |X|$ one has $\trace(\rho_L(\Frob_x))\in K=E_0\otimes_\BQ \BQ_\ell$, and from the definition of $E_0$ and using the \v{C}ebotarov density theorem and the continuity of $\trace\circ\rho_L$ it is clear that $K=K_{\Ad,G/L}$.

Let $(K,H_K,\phi)$ be a minimal quasi-model from \autoref{Thm-ExplicitMinQM} for $\Gamma=\rho_L(\pi_1(X))$. Define $\rho_K\colon\pi_1(X)\to \Aut(\Lie H_K)$ as the composite of $\pi_1(K)\to\Gamma$ with $\Gamma\subset H_K(K)$ and then with $\Ad_H$. Because $\Ad_K\otimes_KL=\Ad_G$, and so $\rho_K\otimes_KL=\rho_L$ we have 
\[\charpol_{\rho_K(\Frob_x)}(T)=\charpol_{\rho_L(\Frob_x)}(T)=P_x(T) \]
for all $x\in|X|$. It follows that $P_x(T)\in E_0\otimes_\BQ\BQ_\ell[T]$, and this for all $\ell\in\CP_\BQ'$, and thus $P_x(T)\in E_0[T]$. Since this holds for all $x\in |X|$, we find $E_0=E'$. Now define $\rho'_{\lambda'}:=\rho_K\otimes_KE'_{\lambda'}$ for $\lambda'\in\CP_{E'}'$ above $\ell$. Then it is clear that $(\rho'_{\lambda'})_{\lambda'\in\CP_{E'}'}$ forms a compatible system defined over~$E'$, and (a) is immediate. Part (b) follows from \autoref{Thm-PinkOpen} and \autoref{Prop-QMforNormal} (we have $G_{\rho'_{\bullet},\lambda'}(E'_{\lambda'})=H_K(E'_{\lambda'})$).

The image of $\pi_1(X)$ under $(\Res_{E'/\BQ}\rho'_{\bullet})_\ell$ lies in the $\BQ_\ell$-points of $\prod_{\lambda'|\ell} \Res_{E'/\BQ_{\ell}} \!G_{\rho'_{\bullet},\lambda'}$ which is semismiple and a closed subgroup of
$ \Res_{E^\prime/\BQ}\GL_n\subset \GL_{n[E^\prime:\BQ]}$. By (b) the image is open in this semisimple subgroup, and this implies both parts of~(c).
\end{proof}

\begin{Def}\label{Def-SigmaConj}
For an $E$-rational compatible system and $\sigma\in\Gal(\overline\BQ/\BQ)$, we define the  {\em $\sigma$-conjugate $\sigma\rho_\bullet$} as the $\sigma(E)$-rational compatible system $\rho_\bullet\otimes_E^\sigma \sigma(E)$.
\end{Def}
Alternatively, if $\rho_\bullet=(\rho_\bullet,(P_x)_{x\in|X|})$ is a $E$-rational compatible and $\sigma\in\Gal(\overline\BQ/\BQ)$, then $\sigma\rho_\bullet=(\sigma\rho_\bullet,(P^\sigma_x)_{x\in|X|})$ can be described as follows, where for $\lambda\in\CP_{\sigma(E)}$ we let $\lambda|\sigma\in\CP_E'$ be the place under $\lambda$ via $\sigma\colon E\to \sigma(E)$: one has $(\sigma\rho_\bullet)_{\lambda}=\rho_{\lambda|\sigma}\otimes_{E_{\lambda|\sigma}}E_\lambda$ and $P^\sigma_x=\sigma(P_x)$.
\begin{Rem}\label{Rem-ConjugacyForMinimal}
Suppose $E$ is minimal for the $E$-compatible system $\rho_\bullet$. Then $\sigma\rho_\bullet$ is isomorphic to $\rho_\bullet$ if and only if $\sigma\in\Gal(\overline\BQ/\BQ)$ fixes $E$. The `if'-direction is obvious. For the `only if'-direction note that if $\sigma\rho_\bullet$ is isomorphic to $\rho_\bullet$, then $P^\sigma_x=P_x\in \overline\BQ[T]$ for all $x\in|X|$. But by the minimality of $E$ this implies that $\sigma$ fixes~$E$.
\end{Rem}
\begin{Lem}\label{Lem-OpennessInOtimes}
Let $E$ be a number field, let $I$ be a finite set, for each $i\in I$, let $E_i$ be a subfield of $E$ and let $\rho_{i,\bullet}$ be an $E_i$-rational simple adjoint compatible system such that each $E_i$ is minimal for $\rho_{i,\bullet}$. Suppose that 
\begin{enumerate}
\item[(i)] the split motivic group of $\rho_\bullet:=\bigotimes_{i\in I} (\rho_{i,\bullet}\otimes_{E_i}E)$ is the product of the split motivic groups of the $\rho_{i,\bullet}$ (over a suitable extension $F$ of $E$),
\item[(ii)] no $\rho_{i,\bullet}$ is $\sigma$-conjugate to some $\rho_{j,\bullet}$ for $i\neq j$ for some $\sigma\in\Gal(\overline\BQ/\BQ)$.
\end{enumerate}
Then for any $\ell\in\CP_\BQ'$, the subgroup
\begin{equation}\label{Eqn-OpenInProduct}
\Big( \bigotimes_{i\in I} \big( \bigoplus_{\lambda_i|\ell}\rho_{i,\lambda_i}\big) \Big) \big(\pi_1^\geo(X)\big)\subset  \prod_{i\in I} \Big(\prod_{\lambda_i|\ell} G_{\rho_{i,\bullet},\lambda_i}(E_{i,\lambda_i}) \Big) 
\end{equation}
is open, where the $\lambda_i$ run through all places of $\CP_{E_i}'$ above $\ell$, and the same holds with $\pi_1(X)$ in place of $\pi_1^\geo(X)$. Moreover $\rho_\bullet:=\bigotimes_{i\in I} (\Res_{E_i/\BQ}\rho_{i,\bullet})_\ell $ is semisimple with monodromy group $G_{\rho_\bullet,\ell}=\prod_{i\in I} G_{\Res_{E_i/\BQ}\rho_{i,\bullet},\ell}$.
\end{Lem}
\begin{proof}
Observe first that we can express \autoref{Eqn-OpenInProduct} also by the use of Weil restrictions as
\[ 
\rho_\ell(\pi_1^\geo(X))\subset \prod_{i\in I} G_{\Res_{E_i/\BQ}\rho_{i,\bullet},\ell}
(\BQ_\ell)  .\]
where we use \autoref{Cor-FieldOfDefOfRhoBullet}(c) on the right. Also, by  \autoref{Cor-FieldOfDefOfRhoBullet}(c), the groups $G_{\Res_{E_i/\BQ}\rho_{i,\bullet},\ell}$ are semisimple connected, and hence so is $\prod_{i\in I} G_{\Res_{E_i/\BQ}\rho_{i,\bullet},\ell}$. We need to show that this product is equal to $G_\ell=G_{\otimes_{i\in I} (\Res_{E_i/\BQ}\rho_{i,\bullet})_\ell,\ell}$, since then the openness follows from \autoref{Cor-Q-CS-OpenImage}. Turning things around, it suffices to establish openness for a single $\ell\in\CP_\BQ'$.

Choose $\ell$ such that $\ell$ is split in $E/\BQ$, so that for all $i$ and $\lambda_i\in\CP_i$ above $\ell$ and $\lambda\in\CP_E$ above $\ell$ we have $E_{i,\lambda_i}=\BQ_\ell$ and $E_\lambda=\BQ_\ell$. To match our setting with \autoref{Thm-ExplicitMinQM} and \autoref{Thm-PinkOpen}, let $L=\oplus_iL_i$ with $L_i=\oplus_{\lambda_i|\ell}E_{i,\lambda_i}$, $G=\prod G_i$ with $G_i=\prod_{\lambda_i|\ell}G_{\rho_{i,\bullet},\lambda_i}$, and let $\Gamma$ be the image of $\pi_1^\geo(X)$ in $G(L)$ under $ \otimes_{i\in I}\big(\oplus_{\gamma_i|\ell} \rho_{i,\lambda_i}\big) $, and $\Gamma_i$ the image of $\Gamma$ under projection to $G_i(L_i)$. By our hypothesis on $E_i$ and $\rho_{i,\bullet}$, and by \autoref{Cor-FieldOfDefOfRhoBullet}(b), the triple $(L_i,G_i,\Gamma_i)$ is minimal. 

Let now $(K,H,\phi)$ be a minimal quasi-model of $(L,G,\Gamma)$. We assume that $K$ is strictly smaller than $L$, which is a finite sum of copies of $\BQ_\ell$. Then there exists a simple factor of $K$, which must be a copy of $\BQ_\ell$, which embeds diagonally into at least two factors of $L$. Since $(L_i,G_i,\Gamma_i)$ is minimal, no two factors can lie in the same $i$. Hence to derive a contradiction, it suffices to assume from now on that $I=\{1,2\}$. By the definition of $K$, and considering the factor just singled out, we have
\begin{enumerate}
\item[(i')] a diagonal embedding $\iota\colon\BQ_\ell\to E_{1,\lambda_1}\times E_{2,\lambda_2}$ for suitable places $\lambda_i$ of $\CP_{E_i}'$ above $\ell$,
\item[(ii')] a semisimple $\BQ_\ell$-group $H_0$ and an isomorphism $H_0\otimes^\iota_{\BQ_\ell} (E_{1,\lambda_1}\!\times E_{2,\lambda_2})\stackrel\simeq\longto G_{\rho_{1,\bullet},\lambda_1}\!\times G_{\rho_{2,\bullet},\lambda_2}$,
\item[(iii')] a homomorphism $\rho_0\colon \pi_1(X)\to H_0(\BQ_\ell)$,
\end{enumerate}
such that, if we denote by $\phi_0$ the diagonal embedding $H_0\to  G_{\rho_{1,\bullet},\lambda_1}\times G_{\rho_{2,\bullet},\lambda_2}$ induced from (ii'), we have $\phi_0\circ \rho_0=\rho_{1,\lambda_1}\times\rho_{2,\lambda_2}$. Since $\iota$ is continuous, it is the canonical isomorphism $\BQ_\ell\cong E_{i,\lambda_i}$ on the factors for $i=1,2$. Interpreted now as an isomorphism $\iota_\ell\colon E_{1,\lambda_1}\to E_{2,\lambda_2}$ the map $\iota$ identifies $\rho_{1,\lambda_1}$ with $\rho_{2,\lambda_2}$. Let $\iota_\lambda\colon E_i\to E_{i,\lambda}$ denote the completion homomorphism. Then for all $x\in |X|$ we have $\iota_\ell\circ\iota_{\lambda_1}(P_{1,x})=\iota_{\lambda_2}(P_{2,x})$. This defines an isomorphism $E_1\to E_2$. By $\sigma$ we denote an extension to $\Gal(\overline\BQ/\BQ)$. Then $(\sigma\rho_{1,\bullet})_{\lambda_2}=\rho_{2,\lambda_2}$, and from \autoref{Lem-IsomOfCS} it follows that $\sigma\rho_{1,\bullet}=\rho_{2,\bullet}$. This was excluded by hypothesis (ii) and thus yields a contradiction.
\end{proof}

The following result shows that Galois conjugacy of compatible systems as tensor factors is an obstruction to an open image theorem, i.e, that hypothesis (ii) in \autoref{Lem-OpennessInOtimes} cannot be omitted. It should be compared with \cite[Thm.~3.2.2]{Loeffler-BigImage} where, by different methods, Galois conjugacy of compatible systems was identified as an obstruction to big image in pairs.
\begin{Prop}\label{Prop-GaloisConjugates}
Suppose $\rho_\bullet$ is $E$-rational simple adjoint and $E$ is minimal for $\rho_\bullet$. Let $F/E$ be finite Galois over $\BQ$, let $I$ be a non-empty subset of $\Hom_\BQ(E,F)$, and fix $\sigma_0\in I$. Then
\begin{enumerate}
\item the split motivic group of $\rho'_\bullet:=\bigotimes_{\sigma\in I}(\sigma\rho_{\bullet}\otimes_{\sigma(E)}F)$ is the product of the split motivic groups of the $\sigma\rho_\bullet$,
\item for any $\ell\in\CP_\BQ'$, the subgroup 

\[\Big(\bigotimes_{\sigma\in I}\big( \bigoplus_{\lambda_\sigma|\ell}
(\sigma\rho_\bullet)_{\lambda_\sigma}\big) \Big) \big(\pi_1^\geo(X)\big)\subset 
 \prod_{\sigma\in I}\Big(\prod_{\lambda_\sigma|\ell}
G_{\sigma\rho_{\bullet},\lambda_\sigma}(\sigma(E)_{\lambda_\sigma}) \Big)\]
is open in the diagonal embedding of $\prod_{\lambda|\ell} G_{\sigma_0\rho_{\bullet},\lambda}(E_{\lambda})$ under $\prod_{\sigma\in I}\sigma\sigma_0^{-1}$, where the $\lambda_\sigma$ run over all places of $\sigma(E)$ above $\ell$.
\end{enumerate}
\end{Prop}
\begin{proof}
It suffices to prove the lemma for $I=\Emb(E,\overline\BQ)$, and $\sigma_0=\id$. That the inclusion in (b) factors via the diagonal embedding of $\prod_{\lambda|\ell} G_{\rho_{\bullet},\lambda}(E_{\lambda})$ under $\prod_{\sigma\in I}\sigma$ is clear. The open image statement follows from \autoref{Cor-FieldOfDefOfRhoBullet}(b) applied to the single factor for $\sigma_0$. This proves~(b).

Regarding (a), we deduce from \autoref{Lem-WeilResAndMonodromy} that $ \Res_{E/\BQ}\rho_\bullet\otimes_\BQ F=\bigotimes_{\sigma\in I}(\sigma\rho_{\bullet}\otimes_{\sigma(E)}F)$, and thus the split motivic group in question is that of $ \Res_{E/\BQ}\rho_\bullet$. Now it follows from  \autoref{Cor-FieldOfDefOfRhoBullet}(c) that $G_{\Res_{E/\BQ}\rho_\bullet,\ell}= \prod_{\lambda|\ell} (\Res_{E_\lambda/\BQ_\ell}G_{\rho_\bullet,\lambda})$. 
Moreover for $\mu_0\in\CP_F'$ above $\ell$ we have, by  \autoref{Lem-CoeffExtOfCompSys}, that $ \prod_{\lambda|\ell} (\Res_{E_\lambda/\BQ_\ell}G_{\rho_\bullet,\lambda}\otimes_{\BQ_\ell} F_{\mu_0})= \prod_{\sigma\in I} G_{\sigma\rho_\bullet,\mu_0} $.  This completes the proof of~(a).
\end{proof}
We now combine \autoref{Lem-OpennessInOtimes} with \autoref{Cor-OpennessAboveEll0} to deduce some further results.
\begin{Thm}\label{Thm-OpennessAboveEll2}
Let $E$, $I$, $E_i$, $\rho_{i,\bullet}$ be as in \autoref{Lem-OpennessInOtimes}, and assume that all hypothesis of that lemma are satisfied. Let furthermore for each $i$ and each $\lambda_i\in\CP_{E_i}'$ denote by  $\Lambda_{i,\lambda_i}\subset E_{i,\lambda_i}^n$ an $\CO_{i,\lambda_i}$-lattice stable under the action by $\pi_1(X)$ via $\rho_{i,\lambda_i}$, and attach to this an $\CO_{i,\lambda_i}$-group scheme $\CG_{i,\lambda_i}$ as above \autoref{Thm-AI-Reduction-IsSaturated}.\footnote{Note that under our hypothesis we have $\CG_{i,\lambda_i}=\CG_{i,\lambda_i}^\geo$.}  Then for almost all $\ell\in\CP_\BQ'$ the following hold:
\begin{enumerate}
\item The $\BZ_\ell$-group scheme $\CG_\ell:=\prod_{i\in I}\prod_{\lambda_i|\ell}\Res_{\CO_{i,\lambda_i}/\BZ_\ell}\CG_{i,\lambda_i}$ is smooth and semisimple.
\item One has 
\[
\Big(\big(\prod_{i\in I}\prod_{\lambda_i|\ell} \Res_{k_{i,\lambda_i}/\BF_\ell}\bar\rho^{\Lambda_i}_{i,\lambda_i}\big)\big(\pi_1^\geo(X)\big)\Big)^\sat_{\BF_\ell} 
=\prod_{i\in I}\prod_{\lambda_i|\ell} \Res_{k_{i,\lambda_i}/\BF_\ell}(\bar\rho^{\Lambda_i}_{i,\lambda_i}(\pi_1^\geo(X))^\sat_{k_{i,\lambda_i}})  =  
\CG_{\BF_\ell},\]
and so 
\[\Big(\prod\limits_{i\in I}\prod\limits_{\lambda_i|\ell} \bar\rho^{\Lambda_i}_{i,\lambda_i}(\pi_1^\geo(X))\Big)^+\!\!\!=\!
\prod\limits_{i\in I}\prod\limits_{\lambda_i|\ell} \bar\rho^{\Lambda_i}_{i,\lambda_i}\big(\pi_1^\geo(X)\big)^+\!\!=\!
\prod\limits_{i\in I}\prod\limits_{\lambda_i|\ell} \CG_{i,\lambda_i}(k_{i,\lambda_i})^+\!=\!
\CG_\ell(\BF_\ell)^+.\]
\item For almost all $\ell$, one has 
\[\prod\limits_{i\in I}\prod\limits_{\lambda_i|\ell} \CG_{i,\lambda_i}(\CO_{i,\lambda_i})^+\subseteq
\Big(\prod\limits_{i\in I}\prod\limits_{\lambda_i|\ell} \rho_{\lambda_i}\Big)\big(\pi_1^\geo(X)\big)\subseteq
\prod\limits_{i\in I}\prod\limits_{\lambda_i|\ell} \CG_{i,\lambda_i}(\CO_{i,\lambda_i})= \CG_\ell(\BZ_\ell).\]
\item Assertions (b) and (c) hold with $\pi_1(X)$ in place of $\pi_1^\geo(X)$.
\item There is a finite subset $\CL$ of $\CP_\BQ'$ and a compact open subgroup $H_\CL\subset \!\prod\limits_{\ell\in\CL}\prod\limits_{i\in I}\prod\limits_{\lambda_i|\ell} G_{\rho_{i,\bullet},\lambda_i}(E_{i,\lambda_i})$ such that 
\[\Big(\prod_{\ell\in\CP_\BQ'}\prod_{i\in I}\prod_{\lambda_i|\ell} \rho_{i,\lambda_i}\Big)(\pi^\geo_1(X))\supseteq 
\prod_{\ell\in\CP_\BQ'\setminus\CL} \Big(\prod_{i\in I}\prod_{\lambda_i|\ell} \CG_{i,\lambda_i}(\CO_{i,\lambda_i})^+\Big)
\times H_\CL.
\]
\end{enumerate}
\end{Thm}
\begin{proof}
Let us first consider the case $\#I=1$ and in this case omit the index $i$ from the notation. Suppose $\ell\gg0$. Then $\ell$ is unramified in $E/\BQ$, so that for all $\lambda\in\CP_E'$ above $\ell$ the ring $\CO_\lambda$ is unramified over $\BZ_\ell$ and hence \'etale. Moreover by \autoref{Prop-LP-GisSmooth} the schemes $\CG_\lambda$ are smooth with semisimple generic fiber. Then by \cite[\S 7.5, Prop.~5]{BGR} the Weil restriction $\Res_{\CO_\lambda/\BZ_\ell}\CG_\lambda$ is smooth connected with semisimple generic fiber, and hence so is $\CG'_\ell:=\prod_{\lambda|\ell} \Res_{\CO_\lambda/\BZ_\ell}\CG_\lambda$. 

Define $\rho_\bullet^\BQ:=\Res_{E/\BQ}\rho_\bullet$. Then the lattice $\Lambda_\ell:=\oplus_{\lambda|\ell}\Lambda_\lambda$ is preserved under the action of $\pi_1(X)$ under $\rho_\ell^\BQ$, and we define $\CH_\ell$ as the Zariski closure in $\Aut_{\pi_1(X)}(\Lambda_\ell)$ where $\pi_1(X)$ acts via $\rho_\bullet^\BQ$. Because $\ell\gg0$ the group scheme $\CH_\ell$ is smooth over $\BZ_\ell$ by \autoref{Prop-LP-GisSmooth}. By its definition, $\CH_\ell\subset\CG_\ell'$ is a closed immersion, and by \autoref{Cor-FieldOfDefOfRhoBullet} the inclusion is an isomorphism on generic fibers. Since both are smooth connected group schemes and of the same relative dimension over $\BZ_\ell$, they must be isomorphic. We deduce from \autoref{Cor-OpennessAboveEll0} that both are also semisimple and that their special fiber satisfies $\big( \bar\rho^{\Lambda}_{\ell}(\pi_1^\geo(X))\big)^\sat_{\BF_\ell}= \CH_{\BF_\ell}$.

It follows that $\Res_{\CO_\lambda/\BZ_\ell}\CG_\lambda$ is semisimple for each $\lambda$ above $\ell$ and hence similar to \autoref{Lem-WeilResAndMonodromy} that each $\CG_\lambda$ is smooth semisimple connected above $\CO_\lambda$. From \autoref{Prop-WeilResAndSat}(e) we furthermore deduce that
\[\big( \bar\rho^{\Lambda}_{\ell}(\pi_1^\geo(X))\big)^\sat_{\BF_\ell} \subset \prod_{\lambda|\ell} \Res_{k_\lambda/\BF_\ell} \big(\bar\rho^\Lambda_\lambda(\pi_1^\geo(X))\big)^\sat_{k_\lambda}. \]
By the previous paragraph, the group  $\CH_{\BF_\ell}$ is equal to the left hand side, but at the same time clearly contains the right hand side. By reducing $\CH_\ell=\CG'_\ell$, we deduce $\CH_{\BF_\ell}=\prod_{\lambda|\ell} \Res_{k_\lambda/\BF_\ell} \CG_{k_\lambda} $. Hence we must have $\big(\bar\rho^\Lambda_\lambda(\pi_1^\geo(X))\big)^\sat_{k_\lambda}=\CG_{k_\lambda}$, giving the first assertion of (b) for $\#I=1$. 

Regarding the second part of (b), from our definitions and from what we proved so far, we have
\[
\bar\rho^{\Lambda}_{\ell}(\pi_1^\geo(X))^+ \!=\big(\prod_{\lambda|\ell} \bar\rho^\Lambda_\lambda(\pi_1^\geo(X))\big)^+\!\subseteq\prod_{\lambda|\ell} \bar\rho^\Lambda_\lambda(\pi_1^\geo(X))^+\!\subseteq\prod_{\lambda|\ell} \CG_\lambda(k_\lambda)^+\!=\CG_\ell(\BF_\ell)^+. \]
The outer terms are equal by \autoref{Cor-OpennessAboveEll0}(b), and this completes the proof of~(b). Moreover part (c) is immediate from our notation and \autoref{Cor-OpennessAboveEll0}(c). Assertion (d) is now clear since the claims in (b) and (c) are weaker if one replaces $\pi^\geo_1(X)$ by $\pi_1(X)$. Part (e) follows from \autoref{Cor-OpennessAboveEll0}(e) applied to $\rho_\bullet^\BQ$. Finally, if $\#I>1$ then one has to also invoke \autoref{Lem-OpennessInOtimes}. We leave the details to the reader.
\end{proof}

We now follow \cite{Hui-Larsen}. Let $G$ be a semisimple group over a field $M$, and let $\pi_G\colon \wt G\to G$ be its universal cover. Using that commutators with central elements are trivial, the commutator map $[\cdot,\cdot]\colon\wt G\times\wt G\to \wt G$ is seen to factor via a map $[\cdot,\cdot]^{\scriptscriptstyle\sim}\colon G\times G\to \wt G$. For $i\ge1$ and $\Gamma$ a subgroup of $G(M)$, denote by $[\Gamma,\Gamma]^{{\scriptscriptstyle\sim} i}$ the subset of $\wt G(M)$ of $i$-fold products of the set $\{[g,h]^{\scriptscriptstyle\sim}\mid g,h\in \Gamma\}$. If for any $\ell\in\CP_\BQ'$ one has a semisimple group $G_\ell$ over $\BQ_\ell$, and if $\Gamma$ is a closed subgroup of $\prod_{\ell\in\CP_\BQ'}G(\BQ_\ell)$, then we also define in the same way $[\Gamma,\Gamma]^{{\scriptscriptstyle\sim} i}$ as a subset of $\prod_{\ell\in\CP_\BQ'}\wt G(\BQ_\ell)$. 

In \cite[\S~2]{Hui-Larsen}, Hui and Larsen introduce the notion of a special adelic group. We alter their~definition by having as an indexing set $\CP_\BQ'$ and not the set set of all rational primes. An example of~such~a group is $\prod_{\ell\in\CP_\BQ'\setminus\CL}\wt\CG_\ell(\BZ_\ell)\times H_\CL$ with $H_\CL$ open compact in $\prod_{\ell\in\CL} \CG_\ell(\BQ_\ell)$, where $\CG_\ell$ is as in \autoref{Thm-OpennessAboveEll2} or \autoref{Cor-OpennessAboveEll0}. We abbreviate $\wt G(\BA_\BQ^p):=\prod'_{\ell\in\CP_\BQ'} \wt\CG_\ell(\BQ_\ell)$, and similarly without the~tilde, where the restriction is with respect to $\wt\CG_\ell(\BZ_\ell)\subset\wt\CG_\ell(\BQ_\ell)$ for $\ell\in\CP_\BQ'\setminus\CL$.
\begin{Cor}[{cf.~\cite[Conj.~1.3]{Hui-Larsen}}]\label{Cor-HuiLarsen}
Suppose either of the following hypotheses:
\begin{enumerate}
\item[(i)] Let $\rho_\bullet$ be a $\BQ$-rational connected compatible system with semisimple split motivic group, and define $\Gamma:=\prod_{\ell\in\CP_\BQ'}\rho_\ell(\pi_1^\geo(X))\subset G(\BA_\BQ^p)$.
\item[(ii)] Let the notation and hypotheses be as in \autoref{Thm-OpennessAboveEll2} and define \[\Gamma:=\Big(\prod_{\ell\in\CP_\BQ'}\prod_{i\in I}\prod_{\lambda_i|\ell} \rho_{i,\lambda_i}\Big)(\pi^\geo_1(X))\subset 
G(\BA_\BQ^p).\]
\end{enumerate}
Then $[\Gamma,\Gamma]^{\scriptscriptstyle\sim}$ generates a special adelic subgroup of $\wt G(\BA_\BQ^p)$, which is equal to $[\Gamma, \Gamma]^{{\scriptscriptstyle\sim}t}$ for some $t$. Moreover $[\Gamma, \Gamma]^{{\scriptscriptstyle\sim}2}$ contains a special adelic subgroup of $\wt G(\BA_\BQ^p)$.
\end{Cor}
\begin{proof}
Using \cite[Thm. 3.8]{Hui-Larsen}, under (i) the assertion follows from Corollaries~\ref{Cor-OpennessAboveEll0} and \ref{Cor-Reduction-IsSaturated}, and under (ii) it is a consequence of the explicit description of $\Gamma$ given in \autoref{Thm-OpennessAboveEll2}(e). 
\end{proof}

\begin{Rem}\label{Rem-AdelicOpen}
We have singled out in this section compatible systems $\rho_\bullet$ with simple adjoint split motivic groups and their adjoint representation as split motivic representations to prove an adelic openness theorem for $E$-rational systems with $E\supsetneq\BQ$ using the results from \cite{Pink-Compact}. The results in op.cit.\ also allow to show the following:

Suppose $\rho_\bullet$ is an $E$-rational compatible system whose split motivic group is semisimple and simple. Let $E_0$ be the subring of $E$ generated over $\BQ$ by $\{\trace(\rho_\lambda(\Frob_x))\mid x\in|X|\}$. Then by \cite[Thm.~2.14]{Pink-Compact} for all $\ell\gg0$ one can define a simple semisimple group $G_\ell$ over $E_{0,\ell}:=(E_0)\otimes_\BQ\BQ_\ell$ such that $\oplus_{\lambda|\ell}\rho_\lambda$ can be obtain by base change $E_{0,\ell}\to E\otimes_\BQ\BQ_\ell$ from a representation $\rho^0_\ell\colon \pi_1(X)\to \GL_{n,E_{0,\ell}}$ with Zariski closure~$G_\ell$. 

For finitely many $\ell$ it may happen in loc.cit.\ that $\rho_\ell$ takes its image in the units of a non-split central simple algebra $D_\ell$ of dimension $n^2$ over $E_0$. Even then one can define a group $G_\ell$ that is a closed subgroup of the unit group of $D_\ell$. The exceptional primes $\ell$ are not controlled by op.cit. In some cases, for instance if the split motivic representation is the adjoint representation, the set of such exceptional places is empty.

Let further $F\supset E$ be a number field over which a split motivic triple $(F,M,\alpha)$ exists for $\rho_\bullet$ and over which the $M$-compatible system $\rho^M_\bullet$ from \autoref{Cor-ChinRepOverM} exists, such that $\rho\otimes_EF\cong \alpha\circ\rho^M_\bullet$. Define $\rho^{\Ad}_\bullet:=\Ad_M\circ \rho^M_\bullet$ and let $E_\ad$ be the subring of $F$ generated over $\BQ$ by $\{\trace(\rho^{\Ad}_\lambda(\Frob_x))\mid x\in|X|\}$. If $E_0=E_\ad$, then by \cite[Prop.~0.6(c)]{Pink-Compact} for every $\ell\gg0$ in $\CP_\BQ'$ the representation $\rho_\ell$ has open image in $G_\ell(E_{0,\ell})$.

A careful reworking of the results of this section will then give an adelic openness result of the image of $\pi_1^\geo(X)$ in a suitably defined restricted product $\prod'_{\ell\in\CP_\BQ'}G(E_\ell)$.
\end{Rem}

Using only results from \autoref{Sec-OpenImage} given above and results up now to and including \autoref{Sec-Saturation}, we shall prove the following corollary which by  \autoref{Rem-AbsIrrThmFromCHToverQ} allows one to give another proof of  \autoref{Thm-AI-Reduction-IsSaturated} based on~\cite{CHT}. It basically says that \autoref{Thm-AI-Reduction-IsSaturated} for $\BQ$-rational compatible systems implies \autoref{Thm-AI-Reduction-IsSaturated} for $E$-rational compatible systems.

\begin{Cor}\label{Cor-SaturationForQRatImpliesARat}
Consider the following statement $(*)_E$: For all $E$-rational absolutely irreducible compatible systems $\rho_\bullet$ such that $\rho_\bullet$ and $\rho_\bullet|_{\pi_1^\geo(X)}$ are connected and all choices of $\pi_1(X)$-stable lattices $\Lambda_\lambda$ of $\CO_\lambda^n$ under the action of $\rho_\lambda$, it holds that for all but finitely many $\lambda\in\CP_E'$ one has equality of semisimple groups
\[\bar\rho_\lambda(\pi_1^\geo(X))^{\sat}_{k_\lambda}=\bar\rho_\lambda(\pi_1(X))^{\sat,o}_{k_\lambda}=\CG^{\geo}_{k_\lambda},\]
where $\CG^\geo_\lambda$ is defined as the Zariski closure of $\rho_{\lambda}(\pi_1^\geo(X))$ in $\Aut_{\CO_{\lambda}}(\Lambda)$.

Then $(*)_\BQ$ implies $(*)_E$ for all number fields $E$.
\end{Cor}
\begin{proof}
To prove the corollary, we may apply the following reductions: It follows from \autoref{Cor-SatAndWtG} that $\bar\rho_\lambda(\pi_1^\geo(X))^{\sat}_{k_\lambda}=\bar\rho_\lambda(\pi_1(X))^{\sat,o}_{k_\lambda}$, and so it suffices to prove that $\bar\rho_\lambda(\pi_1^\geo(X))^{\sat}_{k_\lambda}=\CG^{\geo}_{k_\lambda}$, and that these groups are semisimple. By \autoref{Lem-OnAbsIrredAndTwists} and passing to a finite extension $X'\to X$, we may assume that the split motivic group of $\rho_\bullet$ is semisimple. By~\autoref{Lem-IndepOfLattice},  it suffices to prove the corollary in the case where the split motivic group is an adjoint group $M=\prod_{i\in I} M_i$ with $M_i$ simple adjoint for some finite set $I$. By \autoref{Lem-IndepOfLattice} and \autoref{Cor-ChinAndAI2-New} we may further assume that $\rho_\bullet=\otimes_{i\in I}\rho_{i,\bullet}$, where the split motivic representation of $\rho_{i,\bullet}$ is $\ad_{M_i}$ on $M_i$ and the trivial representation on the factors $M_{i'}$, $i'\neq i$ (cf.~\autoref{Rem-OnChoicesOfSubstitutes}(a)). By \autoref{Cor-FieldOfDefOfRhoBullet}, we can assume that each $\rho_{i,\bullet}$ is $E_i$-rational with $E_i$ minimal for $\rho_{i,\bullet}$. We partition $I$ into a disjoint union $I=\bigcup_{j\in J}I_j$, where $i$ and $i'$ belong to the same class if and only if $\rho_{i,\bullet}$ and $\rho_{i',\bullet}$ are $\sigma$-conjugate for some $\sigma\in\Gal(\overline\BQ/\BQ)$. Because of \autoref{Prop-SatAndOtimes}, we may assume for the proof that for each $i\in I$ and $\sigma\in\Gal(\overline\BQ/\BQ)$, the system $\sigma\rho_{i,\bullet}$ is isomorphic to some $\rho_{i'.\bullet}$ with $i'\in I$. 

For each $j\in J$ we choose one $i_j\in I_j$, and we define $\rho_\bullet^\BQ:=\otimes_{j\in J}\rho_{i_j,\bullet}$. By \autoref{Thm-OpennessAboveEll2}, which in fact only uses $(*)_\BQ$ and not the full strength of \autoref{Cor-OpennessAboveEll0}, we have for $\ell\gg0$
\[\Big(\big(\prod_{j\in J}\prod_{\lambda_j|\ell} \Res_{k_{i_j,\lambda_j}/\BF_\ell}\bar\rho^{\Lambda_{i_j}}_{i_j,\lambda_j}\big)\big(\pi_1^\geo(X)\big)\Big)^\sat_{\BF_\ell}=
\prod_{j\in J}\prod_{\lambda_j|\ell} \Res_{k_{i_j,\lambda_j}/\BF_\ell}\CG^\geo_{i_j,k_{i_j,\lambda_j}},\]
where the notation is as in \autoref{Thm-OpennessAboveEll2}. Let $F/\BQ$ be a finite Galois extension containing all $E_i$ (and $E$), and for each $\ell$ denote by $k_\ell$ the residue field of $F$ above $\ell$. Note also, that by $(*)_\BQ$ the groups on both sides are semisimple for all $\ell_\lambda\gg0$. 

Set $\Lambda'_i:=\sigma\Lambda_{i_j}$ for $i\in I$ with $\rho_{i,\bullet}\cong \sigma\rho_{i_j,\bullet}$ for some $\sigma\in\Gal(\overline \BQ/\BQ)$, and observe that the lattices $\Lambda_i'$ are independent of any choices by \autoref{Rem-ConjugacyForMinimal} since $E_i$ is minimal for $\rho_{i,\bullet}$. Define $\CG^{\prime,\geo}_{i,\lambda_i}$ as the Zariski closure of $\rho_{i,\lambda_i}(\pi_1^\geo(X))$ in $\Aut_{\CO_{\lambda_i}}(\Lambda_i')$. Then for $\ell\gg0$ base change of the above to $k_\ell$ gives
\[\Big(\big(\prod_{i\in I}\prod_{\lambda_i|\ell} \bar\rho^{\Lambda'_i}_{i,\lambda_i}\otimes_{k_{i,\lambda_i}}k_\ell\big)\big(\pi_1^\geo(X)\big)\Big)^\sat_{k_\ell}=
\prod_{i\in I}\prod_{\lambda_i|\ell} \CG^{\prime,\geo}_{i,k_{i,\lambda_i}}\otimes_{k_{\lambda_i}}k_\ell.\]

Now for $\ell\gg0$ the field $F/\BQ$ is unramified above $\ell$. Using the notation $\rho^F_\bullet:=\otimes_{i\in I}\rho_{i,\bullet}\otimes_{E_i}F$, $\Lambda^F_\mu:=\otimes_{i\in I} \Lambda_i'\otimes_{E_{i,\lambda_i}} F_\mu$, and $\CG^F_\mu=\prod_{i\in I,}\CG^{\prime,\geo}_{i,\lambda_i}\otimes_{E_{i,\lambda_i}}F_\mu$ where $\mu$ restricts to $\lambda_i$ on $E_i$, we regroup the previous formula to 
\[\Big(\big(\prod_{\mu|\ell,\mu\in\CP_F'} \bar\rho^{F,\Lambda^F}_{k_\mu}\big)\big(\pi_1^\geo(X)\big)\Big)^\sat_{k_\ell}=
\prod_{\mu|\ell,\mu\in\CP_F'} \CG^F_{k_\mu}.\]
This implies $\bar\rho^{F,\Lambda^F}_{k_\mu}\big(\pi_1^\geo(X)\big)^\sat_{k_\mu}= \CG^F_{k_\mu}$ for all $\mu\in\CP_F'$ with $\ell_\mu\gg0$, and we still have that the groups $\CG^F_{k_\mu}$ are semisimple over $k_\mu$. It follows that the smooth group schemes $\CG^F_\mu$ are smooth and split semisimple over $\CO_\mu$ for $\ell_\mu\gg0$, and that $\bar\rho^{F,\Lambda^F}_{k_\mu}\big(\pi_1^\geo(X)\big)^\sat_{k_\mu} =\bar\rho^{F}_{\mu}\big(\pi_1^\geo(X)\big)^\sat_{k_\mu} $. By our reduction to the case where $I$ is closed under Galois conjugation, we have in fact that $\rho^F_\bullet=\rho_\bullet\otimes_EF$. Hence the previous equality and compatibility of saturation with descent by \autoref{Cor-SatBaseChange} implies that $\bar\rho_{\lambda}\big(\pi_1^\geo(X)\big)^\sat_{k_\lambda}= \CG^{\prime,\geo}_{k_\lambda}$. 

From \autoref{Thm-RedIsSaturated} we have the inclusion $\bar\rho_{\lambda}^\Lambda\big(\pi_1^\geo(X)\big)^\sat_{k_\lambda}\subset\CG^\geo_{k_\lambda}$. An explicit computation invoking $\exp_n$ from \autoref{Subsection-Nori} shows that saturation commutes with semisimplification, i.e., that one has a surjection $\pr_\lambda\colon\bar\rho_{\lambda}^\Lambda\big(\pi_1^\geo(X)\big)^\sat_{k_\lambda}\to \bar\rho_{\lambda}\big(\pi_1^\geo(X)\big)^\sat_{k_\lambda}$. Invoking the previous paragraph, we have
\[\dim \bar\rho_{\lambda}\big(\pi_1^\geo(X)\big)^\sat_{k_\lambda}= \dim\CG^{\prime,\geo}_{k_\lambda} =\dim G_\lambda=\dim \CG^\geo_{k_\lambda}\ge \dim \bar\rho_{\lambda}^\Lambda\big(\pi_1^\geo(X)\big)^\sat_{k_\lambda}.\]
Since all groups involved are smooth, $\pr_\lambda$ must be an isomorphism. Therefore the above inclusion is an equality and $\CG^\geo_{k_\lambda}$ is semisimple.
\end{proof}
\begin{Rem}
By the following argument the semisimplicity assertion in \autoref{Cor-SaturationForQRatImpliesARat} implies that the $\bar\rho_\lambda$ are absolutely irreducible for almost all~$\lambda$:

Let $\CT_\lambda$ be a split torus of $\CG_\lambda$ over $\CO_\lambda$; cf.~\autoref{Prop-LP-GisSmooth}. By \autoref{Prop-LP-GisSmooth} and \autoref{Prop-GenChinTriple}, the weights for the action of $\CT_\lambda$ on $F_\lambda^n$ and on $k_\lambda^n$ via $\rho_\lambda$ and $\bar\rho_\lambda$ are the same and independent of $\lambda$ for $\ell_\lambda\gg0$. The representations of generic fiber of the split group $\CG_{\lambda}$ on $F_\lambda^n$ are irreducible and hence the representation of $\CG_{k_\lambda}$ on $k_\lambda^n$ is irreducible for $\ell_\lambda\gg0$, cf.~\autoref{Prop-LiftingOfLowWeightReps}. Now apply \autoref{Prop-GConnSemisimple} to deduce that $\bar\rho_{k_\lambda}\big(\pi_1^\geo(X)\big)$ acts semisimply on $k_\lambda^n$, and hence absolutely irreducible since this is not changed under saturation.
\end{Rem}

\bibliographystyle{alpha}
\addcontentsline{toc}{section}{References}
\bibliography{FF-CS}

{\small {\sc Gebhard B{\" o}ckle\\
Computational Arithmetic Geometry\\
IWR (Interdisciplinary Center for Scientific Computing)\\
University of Heidelberg\\
Im Neuenheimer Feld 368\\
69120 Heidelberg, Germany}\\
E-mail address: \texttt{\small gebhard.boeckle@iwr.uni-heidelberg.de} 
\par\medskip

{\sc Wojciech Gajda\\
Faculty of Mathematics and Computer Science\\
Adam Mickiewicz University\\
Umultowska 87\\
61614 Pozna\'{n}, Poland}\\
E-mail adress: \texttt{\small gajda@amu.edu.pl}
\par\medskip

{\sc Sebastian Petersen\\ 
Universit\"at Kassel\\
Fachbereich 10\\
Wilhelmsh\"oher Allee 73\\
34121 Kassel, Germany}\\
E-mail address: \texttt{\small petersen@mathematik.uni-kassel.de}}

\end{document}